\documentclass[10pt, a4paper, twoside, titlepage]{article}
\usepackage[a4paper, margin=1.25in, marginparwidth=1in, marginparsep=0.125in]{geometry}
\pdfoutput=1

\usepackage[ngerman, english]{babel} 
\usepackage{amsmath}
\usepackage{amssymb}
\usepackage{amsthm}
\usepackage{thmtools}
\usepackage{algorithm, algpseudocode}
\usepackage{bbm}
\usepackage{blkarray}
\usepackage{pgfplots}
\usetikzlibrary{graphs,arrows.meta,matrix,fit,patterns,shapes,positioning,decorations.pathmorphing,fadings}
\pgfplotsset{compat=1.8}
\usepackage{makecell}
\usepackage{svg}
\usepackage{colonequals}
\usepackage{soul}
\usepackage{color}
\usepackage{setspace}
\usepackage{mathtools}
\usepackage{hyperref}
\usepackage[nameinlink,capitalize]{cleveref}
\usepackage[toc,acronym,symbols]{glossaries}
\usepackage{stmaryrd}
\usepackage{eso-pic}
\usepackage{transparent}
\usepackage{graphicx}
\usepackage[framemethod=TikZ]{mdframed}
\usepackage[shortlabels]{enumitem}
\usepackage[nottoc,numbib]{tocbibind}
\usepackage{ifoddpage}
\usepackage{lmodern} 
\usepackage[LGR,T1]{fontenc}
\usepackage{mathpazo} 
\usepackage{marginnote}
\usepackage{comment}
\usepackage{multirow}

\usepackage{titlesec} 
\numberwithin{equation}{section} 

\usepackage{tabularx}
\usepackage{multicol}
\newcolumntype{Y}{>{\centering\arraybackslash}X}
\usepackage[column=C]{cellspace}
\addparagraphcolumntypes{X}
\definecolor{tumBlue}{RGB}{0,101,189} 
\definecolor{tumDarkBlue}{RGB}{0,82,147} 
\definecolor{tumLightBlue}{RGB}{100,160,200} 
\definecolor{tumLighterBlue}{RGB}{152,198,234} 
\definecolor{tumOrange}{RGB}{227,114,34} 
\definecolor{tumGreen}{RGB}{162,173,0} 
\definecolor{tumGray}{RGB}{153,153,153} 
\definecolor{tumLight}{RGB}{218,215,203} 

\definecolor{commentpurple}{RGB}{102, 0, 102} 
\colorlet{commentpurplemuted}[RGB]{commentpurple!20!white}

\colorlet{sectionblue}{tumBlue}
\definecolor{linkred}{RGB}{127,0,0} 
\definecolor{darklinkred}{RGB}{50,0,0} 
\colorlet{headingcolor}{sectionblue}
\colorlet{headingcolormuted}[RGB]{headingcolor!20!white}
\colorlet{linkcolor}{linkred}

\usepackage[labelfont={color=headingcolor,bf}]{caption}

\makeatletter
\def\moverlay{\mathpalette\mov@rlay}
\def\mov@rlay#1#2{\leavevmode\vtop{%
   \baselineskip\z@skip \lineskiplimit-\maxdimen
   \ialign{\hfil$\m@th#1##$\hfil\cr#2\crcr}}}
\newcommand{\charfusion}[3][\mathord]{
    #1{\ifx#1\mathop\vphantom{#2}\fi
        \mathpalette\mov@rlay{#2\cr#3}
      }
    \ifx#1\mathop\expandafter\displaylimits\fi}
\makeatother

\newcommand{\midd}[1]{\mathrel{}\middle#1\mathrel{}}

\mdfsetup{ skipbelow=-0.2em }

\newlength{\spaceblength}
\declaretheoremstyle[
    headfont=\bfseries,
    notefont=\bfseries,
    notebraces={}{\\[\parskip]}, 
    bodyfont=\normalfont\upshape,
    headpunct={},
    postheadspace=\spaceblength,
    spacebelow=\parskip,
    spaceabove=\parskip,
    headformat=\color{headingcolor}\NAME\ \NUMBER,
    qed={$\begingroup\color{headingcolormuted}\blacktriangleleft\endgroup$}
]{boxstyle}
\declaretheoremstyle[
    headfont=\bfseries\itshape,
    notefont=\normalfont\bfseries,
    notebraces={}{\\[\parskip]}, 
    bodyfont=\normalfont,
    headpunct={},
    postheadspace=\spaceblength,
    spacebelow=\parskip,
    spaceabove=\parskip,
    headformat=\color{headingcolor}\NAME,
    qed=\qedsymbol
]{proofstyle}
\declaretheoremstyle[
    headfont=\bfseries\itshape,
    notefont=\bfseries\itshape\hypersetup{hidelinks},
    notebraces={}{}, 
    bodyfont=\normalfont,
    headpunct={},
    postheadspace=\spaceblength,
    spacebelow=\parskip,
    spaceabove=\parskip,
    headformat={\color{headingcolor}\NAME\ \NOTE},
    qed=\qedsymbol
]{oproofstyle}
\declaretheoremstyle[
    headfont=\bfseries,
    notefont=\bfseries,
    notebraces={}{\\[\parskip]}, 
    bodyfont=\normalfont,
    headpunct={},
    postheadspace=\spaceblength,
    spacebelow=\parskip,
    spaceabove=\parskip,
    headformat=\color{headingcolor}\NAME\ \NUMBER,
    qed={$\begingroup\color{headingcolormuted}\blacktriangleleft\endgroup$}
]{examplestyle}
\declaretheoremstyle[
    headfont=\normalfont,
    notefont=\bfseries,
    notebraces={}{\\[\parskip]}, 
    bodyfont=\normalfont,
    headpunct={},
    postheadspace=\spaceblength,
    spacebelow=\parskip,
    spaceabove=\parskip,
    headformat=\color{headingcolor}\NAME,
    qed={$\begingroup\color{headingcolormuted}\blacktriangleleft\endgroup$}
]{remarkstyle}

\declaretheorem[style=boxstyle]{definition}
\declaretheorem[style=boxstyle,sibling=definition]{theorem}
\declaretheorem[style=boxstyle,sibling=definition]{proposition}
\declaretheorem[style=boxstyle,sibling=definition]{lemma}
\declaretheorem[style=boxstyle,sibling=definition]{corollary}

\declaretheorem[style=proofstyle,numbered=no,name=Proof]{tproof}

\declaretheorem[style=remarkstyle,numbered=no]{remark}

\addto\extrasenglish{

}

\newlength{\secskip}
\titleformat*{\section}{\large\bfseries\color{headingcolor}}
\titleformat*{\subsection}{\bfseries\color{headingcolor}}
\titleformat*{\subsubsection}{\color{headingcolor}}

\hypersetup{
    colorlinks=true,
    allcolors=.,
    urlcolor=linkcolor,
    citecolor=linkcolor,
    linkcolor=linkcolor,
}

\newtagform{roman}[]()

\usepackage{marvosym}

\DeclareMathOperator{\DIST}{dist}

\newcommand{\hpipe}{\rotatebox[origin=c]{90}{$|$}}

\newcommand{\Prb}[1] {\mathbb{P}\left[#1\right]}
\newcommand{\Prbc}[2] {\mathbb{P}_{#1}\left[#2\right]}

\newcommand{\Exc}[2] {\mathbb{E}_{#1}\left[#2\right]}

\newcommand{\dx}[1] {\;\mathrm{d}#1}

\newcommand{\dmux}[2] {\;\mathrm{d}#1\left(#2\right)}

\newcommand{\dufrac}[1] {\frac{\textrm{d}}{\textrm{d}#1}}

\newcommand{\abs}[1] {\left|#1\right|}
\newcommand{\norm}[1] {\left\lVert #1\right\rVert}
\newcommand{\WD}[2] {\mathrm{W}_1\left(#1,#2\right)}

\newcommand{\bbR} {\mathbb{R}}

\newcommand{\bbN} {\mathbb{N}}
\newcommand{\bbZ} {\mathbb{Z}}

\newcommand{\calL} {\mathcal{L}}

\newcommand{\calV} {\mathcal{V}}

\newcommand{\1} {\mathbbm{1}}

\newcommand{\dist}[2] {\DIST\left(#1, #2\right)}

\DeclareMathOperator*{\argmax}{arg\,max}

\tikzset{
  diagonal fill a/.style n args=2{path picture={%
  \draw[fill=#1, draw=none] (path picture bounding box.south west) --
              (path picture bounding box.north east) -- (path picture bounding box.south east) -- cycle; %
  \draw[fill=#2, draw=none] (path picture bounding box.south west) --
              (path picture bounding box.north east) -- (path picture bounding box.north west) -- cycle;}},
  diagonal fill b/.style n args=2{path picture={%
  \draw[fill=#1, draw=none] (path picture bounding box.south west) --
              (path picture bounding box.north east) -- (path picture bounding box.south east) -- cycle; %
  \draw[pattern=north west lines, pattern color=#2, draw=none] (path picture bounding box.south west) --
              (path picture bounding box.north east) -- (path picture bounding box.north west) -- cycle;}},
  diagonal fill c/.style n args=2{path picture={%
  \draw[pattern=north west lines, pattern color=#1, draw=none] (path picture bounding box.south west) --
              (path picture bounding box.north east) -- (path picture bounding box.south east) -- cycle; %
  \draw[fill=#2, draw=none] (path picture bounding box.south west) --
              (path picture bounding box.north east) -- (path picture bounding box.north west) -- cycle;}},
  diagonal fill d/.style n args=2{path picture={%
  \draw[pattern=north west lines, pattern color=#1, draw=none] (path picture bounding box.south west) --
              (path picture bounding box.north east) -- (path picture bounding box.south east) -- cycle; %
  \draw[pattern=north west lines, pattern color=#2, draw=none] (path picture bounding box.south west) --
              (path picture bounding box.north east) -- (path picture bounding box.north west) -- cycle;}},
  table nodes/.style={
    rectangle,
    draw=none,
    align=center,
    minimum height=7mm,
    text depth=0.5ex,
    text height=2ex,
    inner xsep=0pt,
    outer sep=0pt
  },      
  table/.style={
    matrix of nodes,
    row sep=-\pgflinewidth,
    column sep=-\pgflinewidth,
    nodes={
        table nodes
    }
  }
}

\makeatletter
\newcommand{\oset}[3][0ex]{%
  \mathrel{\mathop{#3}\limits^{
    \vbox to#1{\kern-2\ex@
    \hbox{$\scriptstyle#2$}\vss}}}}
\makeatother

\newcommand{\transp} {\scriptscriptstyle \mathsf{T}}

\newcommand{\dabssep}{0.25em}

\newdimen\dabsOn   \newdimen\dabsOff
\newdimen\dabswidth
\newmuskip\dabsLskip
\newmuskip\dabsRskip
\dabsLskip=\thinmuskip   
\dabsRskip=\thinmuskip

\newcommand{\dabs}[1]{%
  \left.\mathchoice
    {\dabsaux\displaystyle{#1}}%
    {\dabsaux\textstyle{#1}}%
    {\dabsaux\scriptstyle{#1}}%
    {\dabsaux\scriptscriptstyle{#1}}%
  \right.%
}
\newcommand{\dabsaux}[2]{%
  \mathopen{}
  \tikz[baseline=(X.base)]{%
    \node[inner sep=0pt, outer sep=0pt] (X) {$#1#2$};
    \draw[line cap=round, line width=\dabswidth, dash pattern=on \dabsOn off \dabsOff]
      ([xshift=-\dabssep]X.north west) -- ([xshift=-\dabssep]X.south west);
    \draw[line cap=round, line width=\dabswidth, dash pattern=on \dabsOn off \dabsOff]
      ([xshift=\dabssep]X.north east) -- ([xshift=\dabssep]X.south east);
  }%
}

\usepackage{ifthen}
\newboolean{shorten}
\setboolean{shorten}{true}

\title{Wasserstein error bounds for aggregations of CTMCs}
\author{Fabian Michel} 
\date{08/10/2025}

\newenvironment{abstractblock}[1][0.8\textwidth]
  {\begin{center}\begin{minipage}{#1}\noindent}
  {\end{minipage}\end{center}}

\begin{document}

\begin{center}
  \textcolor{sectionblue}{
    \LARGE\bfseries
    Wasserstein error bounds for aggregations of continuous-time Markov chains
  }\\[1em]
  Fabian Michel\footnote{Universität der Bundeswehr München, Werner-Heisenberg-Weg 39, 85579 Neubiberg, Germany}
  (\Letter\ \href{mailto:fabian.michel@unibw.de}{fabian.michel@unibw.de})
\end{center}

\begin{abstractblock}
  \textcolor{sectionblue}{\textbf{Abstract.}}
  We study the approximation of a (finite) continuous-time Markov chain by a Markov chain
  on a reduced state space, and we provide formal error bounds for the approximated
  transient distributions in the Wasserstein distance. These bounds extend
  previous work on error bounds in the total variation distance, and
  are the first step towards a generalization to continuous-time
  Markov processes with continuous state spaces. A Wasserstein matrix norm
  is used to bound the error caused by the lower-dimensional approximation
  of the dynamics. In order to control the propagation of the accumulated
  error, we rely on the concept of coarse Ricci curvature of a Markov chain.
  The practical applicability of the presented bounds depends strongly on
  the curvature of the chain. Examples for CTMCs taken from the literature
  (where we added a metric on the state space) show that a negative curvature
  results in exponentially exploding bounds. On the other hand, certain
  CTMCs which we call translation-invariant always have non-negative curvature.
  When measuring the error in the total variation distance (a special case
  of the Wasserstein distance with the discrete metric), the curvature is also
  always non-negative. If it is strictly positive, the
  bounds presented in this paper are an improvement over previous work.

  ~

  \noindent Markov chains \textbullet\ State space reduction \textbullet\ Formal error bounds \textbullet\ Wasserstein distance \textbullet\ Aggregation \textbullet\ Coarse Ricci curvature
\end{abstractblock}

{
  \hypersetup{hidelinks}
  \tableofcontents
}

\section{Introduction}

State aggregation in dynamic systems has been studied extensively since
the 1960s (see \cite{ncd}). Due to the curse of dimensionality, continuous-time
Markov chains with large state spaces can quickly become computationally
intractable without state space reduction. One way to reduce computation time~--
or to turn the model into one which is easier to understand for humans~-- is
to approximate the original model with an aggregated model on a lower-dimensional state space.
Various cases where exact transient or stationary probabilities of
the original model can be derived from an aggregated model have been
identified and analysed (see, e.g.~\cite{exactordlump}).

Formal error bounds for the approximation error when exact
aggregation is not possible have received less attention. \cite{exactordlump}
already gave upper and lower bounds for the transient distribution of
a Markov chain which are derived from an aggregated model.
\cite{adaptformalagg} presented improved bounds for the transient
distribution of discrete-time Markov chains, which can also be applied
to continuous-time Markov chains via uniformisation. These bounds were
extended to a more general setting in \cite{formalbndsstatespaceredmc}.

Both \cite{adaptformalagg} and \cite{formalbndsstatespaceredmc} measured
the error (i.e., the difference between approximated and actual transient
distributions) in the total variation distance. However, this approach
is not suitable for an extension to general continuous-time Markov processes
with continuous state spaces: continuous movement as in the process $X_t = t$ cannot
be reproduced by an aggregated model on a finite, discrete state space, as required
for computation. Therefore, one must commonly allow the approximation
of the transient distribution of such a process to have probability mass
in slightly the wrong place, e.g.\ within a small interval instead of
concentrated on a single point. But the total variation distance is already
equal to the maximal value $1$ when comparing a uniform distribution on
a small interval to a Dirac measure. The Wasserstein distance is better
suited to measure the error in such settings.

This paper still focuses on continuous-time Markov chains with finite
state spaces, but it is intended as a step towards continuous-time Markov processes
with continuous state spaces. The discrete setting is simpler to analyze,
but we expect that many techniques carry over to the continuous setting.

Measuring the error in the Wasserstein distance instead of the total
variation distance introduces additional complications compared to \cite{formalbndsstatespaceredmc}.
While the error caused by the approximation of the dynamics of a continuous-time
Markov chain on a lower-dimensional state space can be controlled in
a similar way to \cite{formalbndsstatespaceredmc}, it is no longer true
that the error accumulated in the calculation up to a given time point
cannot blow up at a later stage. To deal with the accumulated error propagation, the concept of coarse Ricci curvature
of a Markov chain \cite{riccimarkovmetricspaces} turns out to be exactly the right
tool. Essentially, the
coarse Ricci curvature measures the rate at which two transient distributions
of a given Markov chain move toward or away from one another.

\subsection{Our contribution}

Our main contribution is a theory for calculating formal error bounds
for the difference between approximated and actual transient distributions
of a Markov chain in the Wasserstein distance. Such error bounds have
not appeared in the literature before. Our central result is the following:

Consider a continuous-time Markov chain on the finite state space
$S = \{1, \ldots, n\}$ equipped with some metric.
Let $p_0 \in \bbR^n$ be the initial distribution of the continuous-time Markov chain, denote the generator by $Q \in \bbR^{n \times n}$,
such that the transient distribution is $p_t^{\transp} = p_0^{\transp} e^{t Q}$,
and consider the approximation $\widetilde{p}_t^{\transp} = \pi_0^{\transp} e^{t \Theta} A$
with aggregated initial distribution $\pi_0 \in \bbR^m$, aggregated generator $\Theta \in \bbR^{m \times m}$,
and disaggregation matrix $A \in \bbR^{m \times n}$
(details in \autoref{ssec:prelim_mc}).
Similarly to what has been shown in \cite[Theorem 5]{formalbndsstatespaceredmc},
we can prove that (see \autoref{thm:wasserstein_error_growth_bound})
\begin{align}
  \dufrac{t^+} \WD{\widetilde{p}_{t}}{p_{t}}
  &\leq \norm{\Theta A - A Q}_{\mathrm{W}}
  + \WD{\widetilde{p}_{t}}{p_{t}} \cdot (-K)
  \label{eq:intro_central_result}
\end{align}
where $\WD{\cdot}{\cdot}$ is the Wasserstein distance, $\norm{\cdot}_{\mathrm{W}}$ is a Wasserstein
matrix norm on matrices with rows summing to $0$ (see \autoref{def:wasserstein_matrix_norm}),
and where $K$ is a lower bound on the coarse Ricci curvature $\underline{\kappa}(Q)$ of
the Markov chain (see \cite{riccimarkovmetricspaces} and \autoref{def:ricci_curvature}).
Error bounds for the transient distribution at a given time point can be
obtained by integrating \eqref{eq:intro_central_result}.

As a secondary contribution, we provide illustrating and more realistic examples
which show how the bounds behave in practice. Model properties which ensure
desirable error bound properties (non-explosion) are discussed, but the
examples also show where the limitations of the presented bounds are, in particular
for the practical applicability in the case of a negative Ricci curvature
which results in exponentially growing bounds.

\subsection{Paper structure}

\autoref{sec:prelim} introduces the basic concepts: Markov chains, the notion
of aggregation which is used in this paper, the Wasserstein distance and its relation
to linear programs, and finally coarse Ricci curvature as defined by \cite{riccimarkovmetricspaces}.
In \autoref{sec:w_error_bounds}, the central error bounds for the approximated
transient distributions in the Wasserstein distance are derived. The paper mainly
treats continuous-time Markov chains (in \autoref{ssec:w_error_bounds_ctmc}), but
their discrete-time counterpart is also briefly considered (in \autoref{ssec:w_error_bounds_dtmc}).
The propagation and growth of the error accumulated by the aggregation scheme
up to a given time point is bounded in \autoref{sssec:accum_error_growth} with
the help of the coarse Ricci curvature. \autoref{sssec:approx_error_growth}
then shows how the error growth contributed by the approximation on a lower-dimensional
state space can be bounded. The two bounds together yield the central result,
\autoref{thm:wasserstein_error_growth_bound}.

In \autoref{sssec:class_ctmc_nonneg_curv}, we show that Markov chains with
a translation-invariant structure have non-negative curvature which implies
non-exploding error bounds, and in \autoref{sssec:curv_improvement}, we show
how one error bound version given in \autoref{thm:wasserstein_error_growth_bound}
can be slightly improved.
A toy example to illustrate the theory is provided
in \autoref{sssec:toy_example}, followed by a more realistic example in
\autoref{sssec:rsvp}, which is analysed thoroughly and demonstrates some limitations of
the error bounds. More examples of models of a similar size are given in
\autoref{sssec:further_examples}. The conclusion can be found in \autoref{sec:conclusion}.

\section{Preliminaries}
\label{sec:prelim}

\subsection{Markov chains and state space reduction}
\label{ssec:prelim_mc}

We consider time-homogeneous discrete- and continuous-time
Markov chains (DTMCs and CTMCs) on the finite state space $S = \{1, \ldots, n\}$.
The dynamics are given by the stochastic transition matrix $P \in \bbR^{n \times n}$ for DTMCs,
where we have $P(i, j) = \Prb{X_{k+1} = j \midd| X_k = i}$ if $X_k$ denotes the
state of the DTMC at time $k$. For CTMCs, the dynamics are defined via the
generator matrix $Q \in \bbR^{n \times n}$, where $Q(i, j)$ is the transition rate from $i$ to $j$,
and $Q(i, i) = - \sum_{j \neq i} Q(i, j)$.
Given an initial distribution $p_0 \in \bbR^n$, the transient distribution
of a DTMC (respectively CTMC) is given by $p_k^{\transp} = p_0^{\transp} P^k$ (respectively $p_t^{\transp} = p_0^{\transp} e^{tQ}$).

We want to reduce the state space of the Markov chain to speed up computation
of various properties. We often refer to state space reduction as aggregation,
even though there are not necessarily groups of states which are aggregated into one single macro state.
Instead, we define the
aggregation of a Markov chain with an aggregated state space of dimension $m$ (where $m \leq n$) as follows: given a disaggregation matrix
$A \in \bbR^{m \times n}$ with non-negative entries ($A(i, j) \geq 0$) and rows summing to $1$ (i.e., a ``stochastic'', but non-quadratic
matrix with probability distributions in every row), an aggregated transition
matrix $\Pi \in \bbR^{m \times m}$ which is stochastic for DTMCs, an aggregated
generator matrix $\Theta \in \bbR^{m \times m}$ for CTMCs, and an aggregated initial probability distribution $\pi_0 \in \bbR^m$,
we approximate the dynamics of the original chain by setting
$\widetilde{p}_k^{\transp} := \pi_k^{\transp} A := \pi_0^{\transp} \Pi^k A$ and
$\widetilde{p}_t^{\transp} := \pi_t^{\transp} A := \pi_0^{\transp} e^{t \Theta} A$.
$\widetilde{p}_k$ and $\widetilde{p}_t$ are intended
to approximate the transient distributions of the original Markov chains, i.e.~$p_k$ and $p_t$.

We call $A$ disaggregation
matrix since $A$ describes how to blow up the aggregated transient
distribution $\pi_k$ to the full-state-space approximation $\widetilde{p}_k$
via the equation $\widetilde{p}_k^{\transp} = \pi_k^{\transp} A$, which
corresponds to disaggregating $\pi_k$.

The most commonly studied type of aggregation is more restrictive in
the possible choices for $\Pi$, $\Theta$, $A$ and $\pi_0$. In most published approaches, the
state space $S$ of the original chain is partitioned into aggregates
by some partition $\Omega = \{\Omega_1, \ldots, \Omega_m\}$ of $S$, with $\sigma \in \Omega$
being a subset of $S$ which represents all states belonging to one aggregate. The
aggregation function $\omega : S \to \Omega$ maps a state $s$ to
the aggregate to which $s$ belongs, i.e.~$s \in \omega(s)$. Instead
of an arbitrary stochastic disaggregation matrix $A$, one defines probability distributions $\alpha_\sigma \in \bbR^n$ with
support on $\sigma \in \Omega$. As a shorthand, we write $\alpha(s) := \alpha_{\omega(s)}(s)$.
The value $\alpha(s)$ should approximate the conditional probability of being in state $s$
when the chain is in the aggregate $\omega(s)$, i.e.~the probability
$\Prb{X_k = s \midd| X_k \in \omega(s)}$. This probability is in general
dependent on time, but commonly, only time-independent approximations $\alpha$ are considered.
$\alpha_\sigma$ can be thought of as a probability distribution
which splits the probability mass of the aggregate $\sigma$ among its constituting
states in the disaggregation phase, and can in general be chosen by the user.
One can then define the disaggregation matrix $A$ and the aggregation matrix $\Lambda$ as follows:
\begin{align*}
  \Lambda = \begin{pmatrix}
    | &  & | \\
    \1_{\Omega_1} & \hdots & \1_{\Omega_m} \\
    | &  & |
  \end{pmatrix} \in \bbR^{n \times m}, \;\;
  A = \begin{pmatrix}
    \hpipe \, \alpha_{\Omega_1}^{\transp} \, \hpipe \\
    \vdots \\
    \hpipe \, \alpha_{\Omega_m}^{\transp} \, \hpipe
  \end{pmatrix} \in \bbR^{m \times n}
  \qquad \textrm{(note: } A\Lambda = I \textrm{)}
\end{align*}
where $\1_{\sigma} \in \bbR^n$ is defined by
\begin{align*}
  \1_{\sigma}(s) = \begin{cases}
    1 & \textrm{ if } s \in \sigma \\
    0 & \textrm{ otherwise}
  \end{cases}
\end{align*}
A natural definition for $\Pi$ and $\Theta$ is then given by
$\Pi = A P \Lambda$ and $\Theta = A Q \Lambda$, which will ensure
that $\Pi$ is stochastic and that $\Theta$ is a generator.
In this case, $\Pi(\rho, \sigma)$ for $\rho,\sigma \in \Omega$ is
an approximation of the probability to transition from one aggregate state
into another, that is, an approximation of $\Prb{X_{k+1} \in \sigma \midd| X_k \in \rho}$.
Note that this probability may also depend on time (i.e.~on $k$) in general, in contrast to the
probability $\Prb{X_{k+1} = s \midd| X_k = r}$ for $r,s \in S$. However, we again consider only
time-independent approximations of $\Prb{X_{k+1} \in \sigma \midd| X_k \in \rho}$.
Simlarly, for CTMCs, we should have
\begin{align*}
  \Theta(\rho, \sigma) \approx \lim_{u \to 0} \frac{\Prb{X_{t+u} \in \sigma \midd| X_t \in \rho}}{u}
  \qquad \textrm{ for } \rho \neq \sigma
\end{align*}
if we aim at a faithful approximation of the dynamics.
Furthermore,
$\pi_0^{\transp} = p_0^{\transp} \Lambda$ is the natural choice for
the initial distribution when working with actual aggregates.

\subsection{Wasserstein distance}

We will measure the error caused by our aggregation scheme in the Wasserstein
distance \cite{wassersteindist,dobrushinwdist}, sometimes also called Kantorovich-Rubinstein distance \cite{kantrubinduality,krnorm}. Let us first
introduce the Wasserstein distance of two Borel probability measures
$\mu$ and $\nu$ on a general Polish space $S$. The Wasserstein distance depends
on a metric defined on the space $S$, which we will denote by $\DIST$,
and which we require to be lower semi-continuous
(this need not be a metric giving rise to the underlying topology of $S$).
\begin{definition}
  \label{def:wasserstein_polish}
  We define the Wasserstein distance between the two probability measures
  $\mu$ and $\nu$ as (cf.\ \cite[Theorem 1.14 on page 34]{topicsopttransportation}
  and \cite[Theorem 2.10]{lectopttransport} for the existence of the minimum)
  \begin{align}
    \WD{\mu}{\nu} &:= \min_{\gamma \in \Gamma(\mu, \nu)} \int_{S \times S} \dist{x}{y} \dmux{\gamma}{x,y} \label{eq:wasserstein_def_polish} \\
    \textrm{with } \Gamma(\mu, \nu) &:= \textrm{set of all probability measures on } S \times S \notag \\
    &\hphantom{:=\,\,} \textrm{ s.t.\ } \gamma(A \times S) = \mu(A) \textrm{ and } \gamma(S \times A) = \nu(A) \;\; \forall A \textrm{ measurable} \notag
  \end{align}
  $\Gamma(\mu, \nu)$ is the set of all couplings of the two measures $\mu$ and $\nu$.
\end{definition}
The Wasserstein distance measures the distance by which $\mu$'s mass has
to be moved to match $\nu$. The subscript $1$ in $\WD{\mu}{\nu}$ is the usual notation, and
distinguishes the above distance from Wasserstein distances where
$\dist{x}{y}$ is raised to some power within the integral above.

The Kantorovich-Rubinstein theorem \cite[Theorem 1.14 on page 34]{topicsopttransportation} gives
an alternative expression for \eqref{eq:wasserstein_def_polish}:
\begin{align}
  \WD{\mu}{\nu} &= \;\;\sup_{\substack{f : S \to \bbR \textrm{ bounded and }1\textrm{-Lipschitz w.r.t.\ }\DIST\\\abs{f} \textrm{ integrable w.r.t.\ }\abs{\mu - \nu}}} \left(\int_S f \dx{\mu} - \int_S f \dx{\nu}\right)
  \label{eq:wasserstein_dual_polish}
\end{align}
If $S = \{1, \ldots, n\}$, as in the finite-state Markov chain setting, then
by \cite[Remark 1.15 (i) on page 34]{topicsopttransportation} and \cite[Remark 1.4 (v) on page 20]{topicsopttransportation},
\eqref{eq:wasserstein_def_polish} and \eqref{eq:wasserstein_dual_polish} simplify to
\begin{align}
  \WD{p}{q} &= \min_{\gamma \in \Gamma(p, q)} \sum_{r, s \in S} \dist{r}{s} \cdot \gamma(r, s) \label{eq:wasserstein_def_finite}\\
  &= \;\;\max_{\substack{f \in \bbR^n \textrm{ is }1\textrm{-Lipschitz w.r.t.\ }\DIST\\\forall s \in S: 0 \leq f(s) \leq d_{\max}}}
  \left(\sum_{s \in S} f(s) \cdot p(s) - \sum_{s \in S} f(s) \cdot q(s)\right)
  \label{eq:wasserstein_dual_finite}
\end{align}
where $p, q \in \bbR^n$ are probability measures on $S = \{1, \ldots, n\}$ and
$d_{\max} := \max_{r, s \in S} \dist{r}{s}$. Note that $f \in \bbR^n$ being $1$-Lipschitz
simply means that $\abs{f(r) - f(s)} \leq \dist{r}{s}$ for $r, s \in \{1, \ldots, n\}$
in this context, where $f(s)$ is the $s$-th entry of the vector $f$.

\begin{remark}
  The restriction $\forall s \in S: 0 \leq f(s) \leq d_{\max}$ does not change
  the maximum in \eqref{eq:wasserstein_dual_finite}. This is due to two reasons:
  on the one hand, adding a constant to a function $f$ leaves the objective value
  over which we maximize unchanged. On the other hand, because $f$ needs to be $1$-Lipschitz,
  the difference between the maximum and the minimum of $f$ can be at most $d_{\max}$.
  Therefore, we can shift any $1$-Lipschitz $f$ (by adding the appropriate constant)
  such that it falls within the range $[0, d_{\max}]$ while keeping the objective
  value unchanged.

  Hence, we could also completely drop
  the restriction $\forall s \in S: 0 \leq f(s) \leq d_{\max}$, or restrict to non-negative $f$, etc.
\end{remark}

One important example for a metric on $S = \{1, \ldots, n\}$ is the so-called
discrete metric defined by
\begin{align*}
  \dist{r}{s} = \begin{cases}
    1 & \textrm{ if } r \neq s \\
    0 & \textrm{ otherwise}
  \end{cases} \quad \textrm{ for } r, s \in S
\end{align*}
For the discrete metric, we have
\begin{align}
  \begin{split}
    \WD{p}{q} &= \min_{\gamma \in \Gamma(p, q)} \;\;\; \sum_{s \in S} \sum_{r \in S, r \neq s} \gamma(s,r) \\
    &= \min_{\gamma \in \Gamma(p, q)} \;\;\; \sum_{s \in S} \big(p(s) - \gamma(s, s)\big)
    = \sum_{s \in S} \big(p(s) - \min\{p(s), q(s)\}\big) \\
    &\overset{\circledast}{=} \frac{1}{2} \sum_{s \in S} \abs{p(s) - q(s)}
    = \frac{1}{2} \norm{p - q}_1
    = \textrm{total variation distance between } p \textrm{ and } q
  \end{split}
  \label{eq:wd_eq_tovar_discrmet}
\end{align}
For $\circledast$, note that $\min\{p(s), q(s)\} = \frac{1}{2}\big(p(s) + q(s) - \abs{p(s) - q(s)}\big)$.
Hence, if we choose the discrete metric as our metric for the state space,
then we bound the error in the total variation distance and we recover the
setting that was treated in \cite{formalbndsstatespaceredmc}.
The dual expression \eqref{eq:wasserstein_dual_finite} can also be reduced
to a simplified version for the discrete metric:
\begin{align*}
  \WD{p}{q} &= \;\;\max_{f \in \bbR^n \textrm{ s.t.\ }\forall s \in S: 0 \leq f(s) \leq 1}
  \left(\sum_{s \in S} f(s) \cdot p(s) - \sum_{s \in S} f(s) \cdot q(s)\right)
\end{align*}

\begin{remark}
  On a finite state space, we can derive the following relation between the
  total variation distance and the Wasserstein distance for a general metric
  (not necessarily the discrete one):
  \begin{align}
    \begin{split}
      \WD{p}{q} &= \min_{\gamma \in \Gamma(p, q)} \;\;\; \sum_{s \in S} \sum_{r \in S, r \neq s} \dist{s}{r} \cdot \gamma(s,r) \\
      &\leq \min_{\gamma \in \Gamma(p, q)} \;\;\; \sum_{s \in S} \sum_{r \in S, r \neq s} d_{\max} \cdot \gamma(s,r)
      \;\overset{\textrm{\eqref{eq:wd_eq_tovar_discrmet}}}{=}\; \frac{d_{\max}}{2} \cdot \norm{p - q}_1
    \end{split}
    \label{eq:wd_tovar_bound}
  \end{align}
  where we write again $d_{\max} = \max_{r, s \in S} \dist{r}{s}$.
  That is, the Wasserstein distance is at most the diameter of the space
  times the total variation distance.
\end{remark}

\subsubsection{Wasserstein norm for matrices}

Next to probability measures, the Wasserstein distance can also be applied
to any two measures with equal total mass with the definition from
\eqref{eq:wasserstein_dual_polish}, or with \eqref{eq:wasserstein_def_polish}
where the coupled measure needs to have the same total mass as the
individual measures. We will use that extension for the error bounds
which we develop later for the aggregation scheme.
For these bounds, it will also be helpful to define a Wasserstein
norm for matrices.

\begin{definition}
  \label{def:row_wd_norm}
  Let $D \in \bbR^{m \times n}$ with rows summing to $0$ and assume
  that $\DIST$ is a metric on $S = \{1, \ldots, n\}$ with
  $d_{\max} = \max_{r, s \in S} \dist{r}{s}$. We define the column vector
  \begin{align*}
    \dabs{D}_{\mathrm{W}}
    := \begin{pmatrix}
      \max_{f \in \bbR^n \textrm{ is } 1\textrm{-Lip.\ w.r.t.\ }\DIST, \; \forall s \in S: \, 0 \leq f(s) \leq d_{\max}} \; D_1 f \\
      \max_{f \in \bbR^n \textrm{ is } 1\textrm{-Lip.\ w.r.t.\ }\DIST, \; \forall s \in S: \, 0 \leq f(s) \leq d_{\max}} \; D_2 f \\
      \vdots \\
      \max_{f \in \bbR^n \textrm{ is } 1\textrm{-Lip.\ w.r.t.\ }\DIST, \; \forall s \in S: \, 0 \leq f(s) \leq d_{\max}} \; D_m f
    \end{pmatrix}
    \in \bbR^m
  \end{align*}
  Here, $D_i$ denotes the $i$-th row of $D$.
\end{definition}

Note that the rows of both $\Theta A - A Q$ in the CTMC setting and of $\Pi A - A P$ in the DTMC setting
sum to $0$ so that \autoref{def:row_wd_norm} is applicable to these matrices.

\begin{remark}
  To clarify the relation to the Wasserstein distance, consider two matrices $B, C \in \bbR^{m \times n}$ with non-negative entries and
  rows summing to $1$. Then, every row of each matrix corresponds to a probability
  distribution, and we have
  \begin{align*}
    \dabs{B - C}_{\mathrm{W}}
    = \begin{pmatrix}
      \max_{f \in \bbR^n \textrm{ is } 1\textrm{-Lip.\ w.r.t.\ }\DIST, \, \forall s: 0 \leq f(s) \leq d_{\max}} \; (B_1 - C_1)f \\
      \max_{f \in \bbR^n \textrm{ is } 1\textrm{-Lip.\ w.r.t.\ }\DIST, \, \forall s: 0 \leq f(s) \leq d_{\max}} \; (B_2 - C_2)f \\
      \vdots \\
      \max_{f \in \bbR^n \textrm{ is } 1\textrm{-Lip.\ w.r.t.\ }\DIST, \, \forall s: 0 \leq f(s) \leq d_{\max}} \; (B_m - C_m)f
    \end{pmatrix} \overset{\textrm{\eqref{eq:wasserstein_dual_finite}}}{=}
    \begin{pmatrix}
      \WD{B_1}{C_1} \\
      \WD{B_2}{C_2} \\
      \vdots \\
      \WD{B_m}{C_m}
    \end{pmatrix}
  \end{align*}
  Hence, if $\dabs{\cdot}_{\mathrm{W}}$ is applied to the difference of two matrices $B$ and $C$
  which both contain probability measures as rows, then $\dabs{B - C}_{\mathrm{W}}$
  is a column vector with each entry corresponding to the Wasserstein distance between
  the two respective row measures in $B$ and $C$.

  In general, $\dabs{D}_{\mathrm{W}}$ measures, for every row $D_i$, the
  Wasserstein distance between the positive part of the row $D_i^+$ (the
  entry-wise maximum of $0$ and the respective row entries) and the negative
  part of the row $D_i^-$ (the negative of the entry-wise minimum of $0$ and the respective row entries).
  As each row $D_i$ is assumed to sum to $0$, $D_i^+$ and $D_i^-$ sum to the
  same total mass, so we can measure the Wasserstein distance between them
  (using the slightly extended definition mentioned at the beginning of this subsection).
\end{remark}

If $\DIST$ is the discrete metric, and $D$ a matrix with rows summing to $0$,
then $\dabs{D}_{\mathrm{W}} = \frac{1}{2} \abs{D} \cdot \mathbf{1}_n$
(here, $\abs{\cdot}$ is the element-wise absolute value and $\mathbf{1}_n \in \bbR^n$ is the
column vector consisting only of ones).

In a very similar way to the definition of $\dabs{\cdot}_{\mathrm{W}}$, we can define a Wasserstein norm
for matrices.

\begin{definition}
  \label{def:wasserstein_matrix_norm}
  Let $D \in \bbR^{m \times n}$. We define
  \begin{align*}
    \norm{D}_{\mathrm{W}} &:= 
    \begin{cases}
      \displaystyle \max_{i \in \{1, \ldots, m\}} \; \max_{\substack{f \in \bbR^n \textrm{ is } 1\textrm{-Lip.\ w.r.t.\ }\DIST,\\\forall s: 0 \leq f(s) \leq d_{\max}}} \; D_i f & \quad\textrm{if all rows of } D \textrm{ sum to } 0 \\
      \infty & \quad\textrm{otherwise}
    \end{cases} \qedhere
  \end{align*}
\end{definition}

$\norm{\cdot}_{\mathrm{W}}$ is a norm on the space of matrices with rows summing to $0$.
This can be seen by noting that
$\norm{\cdot}_{\mathrm{W}}$ is the maximum of the row-wise Kantorovich-Rubinstein norm (see \cite[Chapter VIII, §4, 4.3]{funcanakantorovich} or \cite{krnorm}, for example),
and therefore inherits the norm properties directly. $\norm{\cdot}_{\mathrm{W}}$ is not sub-multiplicative in general.
Furthermore, if $\DIST$ is the discrete metric, then (for a matrix $D$ with
rows summing to $0$) $\norm{D}_{\mathrm{W}} = \frac{1}{2}\norm{D}_{\infty}$, where $\norm{D}_{\infty}$
is the matrix norm given by
\begin{align*}
  \norm{D}_{\infty} = \max_{1 \leq i \leq m} \sum_{j=1}^n \abs{D(i, j)}
\end{align*}

\subsubsection{Linear programs and the Wasserstein distance}

In this subsection, we show alternative formulations for calculating the Wasserstein
distance and take a closer look at the two dual ways for its representation.
Consider the finite state space case $S = \{1, \ldots, n\}$ and the corresponding
forms for the Wasserstein distance in \eqref{eq:wasserstein_def_finite} and \eqref{eq:wasserstein_dual_finite}.
The duality between \eqref{eq:wasserstein_def_finite} and \eqref{eq:wasserstein_dual_finite}
follows directly from the duality in linear programming, as is shown in the
proof of the following proposition.

\begin{proposition}
  \label{prop:wasserstein_equiv}
  Let $p, q \in \bbR^n$ be probability measures
  on the state space $S = \{1, \ldots, n\}$ with metric $\DIST$.
  Then, we have
  \begin{align}
    \WD{p}{q} \textrm{ is the solution of}\;\;
    \max_{f \in \bbR^n, f \geq 0} (p^{\transp} - q^{\transp}) f
    \;\;\textrm{s.t.}\;\;
    \forall r, s \in S: f(r) - f(s) \leq \dist{r}{s}
    \label{eq:wasserstein_finite_primallp}
  \end{align}
  and, equivalently (by linear programming duality),
  \begin{align}
    &\WD{p}{q} \textrm{ is the solution of} \notag\\
    &\min_{\gamma \in \bbR^{n \times n}, \gamma \geq 0} \sum_{r,s \in S} \dist{r}{s} \gamma(r, s)
    \quad\textrm{s.t.}\quad
    \forall r \in S: \sum_{s \in S} \gamma(r, s) - \sum_{s \in S} \gamma(s, r) \geq p(r) - q(r)
    \label{eq:wasserstein_finite_duallp}
  \end{align}
  Furthermore, there is a pair of optimal solutions $f^\ast, \gamma^\ast$ of
  \eqref{eq:wasserstein_finite_primallp} and \eqref{eq:wasserstein_finite_duallp}
  which satisfies all of the following:
  \begin{enumerate}[(i)]
    \item \label{prop:optgammacoupling} $\displaystyle \gamma^\ast \in \Gamma(p, q) \textrm{, }\quad\textrm{i.e., } \gamma^\ast \textrm{ is a coupling of } p \textrm{ and } q$
    \item \label{prop:optgammaonlyincorout} $\displaystyle \forall r \in S: \quad \sum_{\substack{s \in S\\s \neq r}} \gamma^\ast(r, s) = 0 \quad \textrm{or} \quad \sum_{\substack{s \in S\\s \neq r}} \gamma^\ast(s, r) = 0$
    \item \label{prop:optfbounded} $\displaystyle \forall r \in S: \quad 0 \leq f^\ast(r) \leq d_{\max} \quad\textrm{with}\quad d_{\max} := \max_{r, s \in S} \dist{r}{s}$
    \item \label{prop:optfgammarelation} $\displaystyle \forall r, s \in S: \quad \gamma^\ast(r, s) > 0 \implies f^\ast(r) - f^\ast(s) = \dist{r}{s}$ \qedhere
  \end{enumerate}
\end{proposition}

\begin{tproof}
  The duality of \eqref{eq:wasserstein_finite_primallp} and
  \eqref{eq:wasserstein_finite_duallp} follows directly from standard
  linear programming duality, see e.g.\ \cite[Theorem 5.2]{linearprogramming}. As a corollary, we can show:

  \textbf{Proof of the duality of \eqref{eq:wasserstein_def_finite} and \eqref{eq:wasserstein_dual_finite}}:
  \eqref{eq:wasserstein_finite_primallp} clearly gives the same value
  as \eqref{eq:wasserstein_dual_finite} by the remark just after
  \eqref{eq:wasserstein_dual_finite}.
  To show that \eqref{eq:wasserstein_finite_duallp} has the same optimal value as \eqref{eq:wasserstein_def_finite}, we first
  note that the values of $\gamma(s, s)$ are irrelevant for the solution of
  \eqref{eq:wasserstein_finite_duallp}. It then suffices to show that at least one optimal $\gamma$ from \eqref{eq:wasserstein_finite_duallp} satisfies
  \begin{align}
    \forall r: \qquad \sum_{s \neq r} \gamma(r, s) \leq p(r), \qquad &\sum_{s \neq r} \gamma(s, r) \leq q(r)
    \label{eq:duallp_to_coupling}
  \end{align}
  which shows that one optimal $\gamma$ from \eqref{eq:wasserstein_finite_duallp} does indeed correspond
  to a coupling of $p$ and $q$.
  We can see that, in \eqref{eq:wasserstein_finite_duallp}, we must have $\forall r: \sum_s \gamma(r, s) - \sum_s \gamma(s, r) = p(r) - q(r)$
  because the left hand sides as well as the right hand sides of the inequalities sum up to $0$ when
  summing over $r$. In particular, as the common term $\gamma(r, r)$ in the two sums cancels, we have
  \begin{align}
    \forall r: \sum_{s \neq r} \gamma(r, s) - \sum_{s \neq r} \gamma(s, r) &= p(r) - q(r)
    \label{eq:duallp_equalities}
  \end{align}
  In order to show \eqref{eq:duallp_to_coupling},
  we will show below that (for at least one optimal $\gamma$ in \eqref{eq:wasserstein_finite_duallp})
  \begin{align}
    \forall r: \quad \sum_{s \neq r} \gamma(r, s) = 0 \quad \textrm{or} \quad \sum_{s \neq r} \gamma(s, r) = 0
    \label{eq:duallp_to_coupling_simplification}
  \end{align}
  As $\gamma \geq 0$ entry-wise, this implies, together with \eqref{eq:duallp_equalities},
  that
  \begin{align*}
    \sum_{s \neq r} \gamma(r, s) = \begin{cases}
      p(r) - q(r) & \textrm{ if } p(r) - q(r) \geq 0 \\
      0 & \textrm{ otherwise}
    \end{cases} \;\;\; \leq p(r)
  \end{align*}
  (and the same inequality for $\sum_{s \neq r} \gamma(s, r)$ and $q(r)$) as desired.

  To conclude, we now have to show \eqref{eq:duallp_to_coupling_simplification}. Assume for a contradiction
  that for all optimal $\gamma$ from \eqref{eq:wasserstein_finite_duallp}, there is some $r$ with $\sum_{s \neq r} \gamma(r, s) > 0$ and $\sum_{s \neq r} \gamma(s, r) > 0$. Then,
  there must be $u, s$ with $u \neq r$ and $s \neq r$ such that $\gamma(r, u) > 0$ and $\gamma(s, r) > 0$.
  We set $\varepsilon = \min\{\gamma(r, u), \gamma(s, r)\} > 0$. Then, we can define
  \begin{align*}
    \widetilde{\gamma}(r, u) &= \gamma(r, u) - \varepsilon \geq 0 \\
    \widetilde{\gamma}(s, r) &= \gamma(s, r) - \varepsilon \geq 0 \\
    \widetilde{\gamma}(s, u) &= \gamma(s, u) + \varepsilon \\
    \widetilde{\gamma}(\widetilde{r}, \widetilde{s}) &= \gamma(\widetilde{r}, \widetilde{s})
    \textrm{ for all other pairs }\widetilde{r}, \widetilde{s}
  \end{align*}
  Note that we still have
  \begin{align*}
    \sum_{\widetilde{s} \neq r} \widetilde{\gamma}(r, \widetilde{s})
    - \sum_{\widetilde{s} \neq r} \widetilde{\gamma}(\widetilde{s}, r)
    = \sum_{\widetilde{s} \neq r} \gamma(r, \widetilde{s})
    - \varepsilon
    - \sum_{\widetilde{s} \neq r} \gamma(\widetilde{s}, r) + \varepsilon
    = \sum_{\widetilde{s} \neq r} \gamma(r, \widetilde{s})
    - \sum_{\widetilde{s} \neq r} \gamma(\widetilde{s}, r)
  \end{align*}
  and equivalent equations for $u, s$ (as well as for all other states, where
  the value of $\widetilde{\gamma}$ remains unchanged from $\gamma$), so
  $\widetilde{\gamma}$ still satisfies \eqref{eq:duallp_equalities}, i.e.,
  $\widetilde{\gamma}$ is an admissible solution for the linear program
  \eqref{eq:wasserstein_finite_duallp}. However, we see that
  \begin{align*}
    \sum_{\widetilde{r}, \widetilde{s}} \dist{\widetilde{r}}{\widetilde{s}} \widetilde{\gamma}(\widetilde{r}, \widetilde{s})
    &= \sum_{\widetilde{r}, \widetilde{s}} \dist{\widetilde{r}}{\widetilde{s}} \gamma(\widetilde{r}, \widetilde{s}) + \underbrace{\varepsilon}_{> 0} \cdot \underbrace{\big(\dist{s}{u} - \dist{r}{u} - \dist{s}{r}\big)}_{\leq 0 \textrm{ by }\triangle\textrm{-inequ.}} \\
    &\leq \sum_{\widetilde{r}, \widetilde{s}} \dist{\widetilde{r}}{\widetilde{s}} \gamma(\widetilde{r}, \widetilde{s})
  \end{align*}
  Hence, the still admissible $\widetilde{\gamma}$ achieves an objective value
  smaller or equal than that achieved by $\gamma$. If the inequality is strict,
  we have a contradiction, if not, we can iterate the procedure until
  we reach a $\gamma$ of the desired form (this iteratrion must terminate
  because in every iteration, $\sum_{r \neq s} \gamma(r, s)$ is decreasing
  (the mass is actually moved to the diagonal, but this is hidden in our argument,
  because the diagonal entries of $\gamma$ are not relevant for the linear
  program in \eqref{eq:wasserstein_finite_duallp}), and if we go through
  all $r$ state by state to eliminate one of the two sums in \eqref{eq:duallp_to_coupling_simplification}, it is easy to check
  that for a later state $\widetilde{r}$ in the iteration, the sum which was
  set to $0$ for $r$ will remain unchanged).

  \textbf{Proof of \ref{prop:optgammacoupling}--\ref{prop:optfgammarelation}}:
  The existence of $\gamma^\ast$ satisfying \ref{prop:optgammacoupling} and
  \ref{prop:optgammaonlyincorout} follows from the previous part of the proof
  and in particular from \eqref{eq:duallp_to_coupling_simplification}.
  We now construct an optimal $f^\ast$ for \eqref{eq:wasserstein_finite_primallp} which
  satisfies \ref{prop:optfbounded} and \ref{prop:optfgammarelation}. Note that
  we can shift any admissible solution of \eqref{eq:wasserstein_finite_primallp}
  such that \ref{prop:optfbounded} is satisfied by the remark after \eqref{eq:wasserstein_dual_finite}. Hence, we only have to show that
  an optimal $f^\ast$ which satisfies \ref{prop:optfgammarelation} exists.

  In fact, we can choose any optimal $f^\ast$ for \eqref{eq:wasserstein_finite_primallp}
  and then invoke complementary slackness.
  As $\gamma^\ast$ is optimal for the dual \eqref{eq:wasserstein_finite_duallp},
  we have by \cite[Theorem 5.3]{linearprogramming} that
  \begin{align*}
    \forall r, s: \gamma^\ast(r, s) \cdot \underbrace{\bigg(\dist{r}{s} - \big(f^\ast(r) - f^\ast(s)\big)\bigg)}_{\textrm{primal slack}} = 0
  \end{align*}
  \ref{prop:optfgammarelation} follows immediately.
\end{tproof}

\subsubsection{Wasserstein distance in an example}

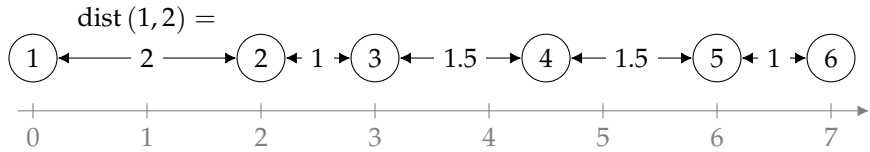
\begin{figure}[htb]
  \begin{center}
    \begin{tikzpicture}[>={Latex[length=1.5mm,width=1.5mm]}]
      \node[circle,draw=black] (A) at (0,0) {$1$};
      \node[circle,draw=black] (B) at (3,0) {$2$};
      \node[circle,draw=black] (C) at (4.5,0) {$3$};
      \node[circle,draw=black] (D) at (6.75,0) {$4$};
      \node[circle,draw=black] (E) at (9,0) {$5$};
      \node[circle,draw=black] (F) at (10.5,0) {$6$};
      \draw[<->] (A) -- node[fill=white] {$2$} node[above=2mm] {$\dist{1}{2}=$} (B);
      \draw[<->] (B) -- node[fill=white] {$1$} (C);
      \draw[<->] (C) -- node[fill=white] {$1.5$} (D);
      \draw[<->] (D) -- node[fill=white] {$1.5$} (E);
      \draw[<->] (E) -- node[fill=white] {$1$} (F);
      \draw[->,gray] (-0.2,-0.7) -- (11,-0.7);
      \draw[gray] (0,-0.6) -- (0,-0.8) node[below] {$0$};
      \draw[gray] (1.5,-0.6) -- (1.5,-0.8) node[below] {$1$};
      \draw[gray] (3,-0.6) -- (3,-0.8) node[below] {$2$};
      \draw[gray] (4.5,-0.6) -- (4.5,-0.8) node[below] {$3$};
      \draw[gray] (6,-0.6) -- (6,-0.8) node[below] {$4$};
      \draw[gray] (7.5,-0.6) -- (7.5,-0.8) node[below] {$5$};
      \draw[gray] (9,-0.6) -- (9,-0.8) node[below] {$6$};
      \draw[gray] (10.5,-0.6) -- (10.5,-0.8) node[below] {$7$};
    \end{tikzpicture}
  \end{center}
  \caption{A line metric for the state space $S = \{1, \ldots, 6\}$}
  \label{fig:wdist_example_metric}
\end{figure}

Here, we provide an example to illustrate the concept of Wasserstein distance.
Consider the state space $S = \{1, \ldots, 6\}$ with the line metric given
in \autoref{fig:wdist_example_metric}~-- the distance between two states is simply
the distance of their two locations on the line.
Let us further consider the probability distributions
$p = (0.35, \; 0.25, \; 0.05, \; 0.25, \; 0.1, \; 0)^{\transp}$ and
$q = (0.2, \; 0.45, \; 0.05, \; 0, \; 0.05, \; 0.25)^{\transp}$.
As mentioned already briefly, the Wasserstein distance is the cost of the
optimal transport plan for moving $p$'s mass such that it matches $q$'s mass.
With line metrics, the calculation of the Wasserstein distance is
relatively simple as there is always only a single path to transfer
mass from one location to another~-- along the line. Indeed, the mass difference
$\sum_{i=1}^k p(i) - \sum_{i=1}^k q(i)$ must always be shifted along
the line connecting state $k$ to state $k+1$, the direction of the shift
depending on the sign.

\begin{table}[htb]
  \begin{center}
    \begin{tikzpicture}[>={Latex[length=1.5mm,width=1.5mm]}]
      \node[rectangle,minimum size=1cm,inner sep=0pt,draw=gray] at (5,5) {$0$};
      \node[rectangle,minimum size=1cm,inner sep=0pt,draw=gray,fill=linkred!40!white] at (5,4) {$0.1$};
      \node[rectangle,minimum size=1cm,inner sep=0pt,draw=gray,fill=linkred!60!white] at (5,3) {$0.15$};
      \node[rectangle,minimum size=1cm,inner sep=0pt,draw=gray] at (5,2) {$0$};
      \node[rectangle,minimum size=1cm,inner sep=0pt,draw=gray] at (5,1) {$0$};
      \node[rectangle,minimum size=1cm,inner sep=0pt,draw=gray] at (5,0) {$0$};
      \node[rectangle,minimum size=1cm,inner sep=0pt,draw=gray] at (4,5) {$0$};
      \node[rectangle,minimum size=1cm,inner sep=0pt,draw=gray] at (4,4) {$0$};
      \node[rectangle,minimum size=1cm,inner sep=0pt,draw=gray,fill=linkred!20!white] at (4,3) {$0.05$};
      \node[rectangle,minimum size=1cm,inner sep=0pt,draw=gray] at (4,2) {$0$};
      \node[rectangle,minimum size=1cm,inner sep=0pt,draw=gray] at (4,1) {$0$};
      \node[rectangle,minimum size=1cm,inner sep=0pt,draw=gray] at (4,0) {$0$};
      \node[rectangle,minimum size=1cm,inner sep=0pt,draw=gray] at (3,5) {$0$};
      \node[rectangle,minimum size=1cm,inner sep=0pt,draw=gray] at (3,4) {$0$};
      \node[rectangle,minimum size=1cm,inner sep=0pt,draw=gray] at (3,3) {$0$};
      \node[rectangle,minimum size=1cm,inner sep=0pt,draw=gray] at (3,2) {$0$};
      \node[rectangle,minimum size=1cm,inner sep=0pt,draw=gray] at (3,1) {$0$};
      \node[rectangle,minimum size=1cm,inner sep=0pt,draw=gray] at (3,0) {$0$};
      \node[rectangle,minimum size=1cm,inner sep=0pt,draw=gray] at (2,5) {$0$};
      \node[rectangle,minimum size=1cm,inner sep=0pt,draw=gray] at (2,4) {$0$};
      \node[rectangle,minimum size=1cm,inner sep=0pt,draw=gray,fill=linkred!20!white] at (2,3) {$0.05$};
      \node[rectangle,minimum size=1cm,inner sep=0pt,draw=gray] at (2,2) {$0$};
      \node[rectangle,minimum size=1cm,inner sep=0pt,draw=gray] at (2,1) {$0$};
      \node[rectangle,minimum size=1cm,inner sep=0pt,draw=gray] at (2,0) {$0$};
      \node[rectangle,minimum size=1cm,inner sep=0pt,draw=gray] at (1,5) {$0$};
      \node[rectangle,minimum size=1cm,inner sep=0pt,draw=gray] at (1,4) {$0$};
      \node[rectangle,minimum size=1cm,inner sep=0pt,draw=gray] at (1,3) {$0$};
      \node[rectangle,minimum size=1cm,inner sep=0pt,draw=gray,fill=linkred!20!white] at (1,2) {$0.05$};
      \node[rectangle,minimum size=1cm,inner sep=0pt,draw=gray,fill=linkred,text=white] at (1,1) {$0.25$};
      \node[rectangle,minimum size=1cm,inner sep=0pt,draw=gray,fill=linkred!60!white] at (1,0) {$0.15$};
      \node[rectangle,minimum size=1cm,inner sep=0pt,draw=gray] at (0,5) {$0$};
      \node[rectangle,minimum size=1cm,inner sep=0pt,draw=gray] at (0,4) {$0$};
      \node[rectangle,minimum size=1cm,inner sep=0pt,draw=gray] at (0,3) {$0$};
      \node[rectangle,minimum size=1cm,inner sep=0pt,draw=gray] at (0,2) {$0$};
      \node[rectangle,minimum size=1cm,inner sep=0pt,draw=gray] at (0,1) {$0$};
      \node[rectangle,minimum size=1cm,inner sep=0pt,draw=gray,fill=linkred!80!white] at (0,0) {$0.2$};

      \fill[sectionblue!78!white] (-0.6, -0.4) -- (-0.6, 0.4) -- (-1.65, 0.4) -- node[left,text=black] {$0.35$} (-1.65, -0.4) -- cycle;
      \fill[sectionblue!56!white] (-0.6, 0.6) -- (-0.6, 1.4) -- (-1.35, 1.4) -- node[left,text=black] {$0.25$} (-1.35, 0.6) -- cycle;
      \fill[sectionblue!11!white] (-0.6, 1.6) -- (-0.6, 2.4) -- (-0.75, 2.4) -- node[left,text=black] {$0.05$} (-0.75, 1.6) -- cycle;
      \fill[sectionblue!56!white] (-0.6, 2.6) -- (-0.6, 3.4) -- (-1.35, 3.4) -- node[left,text=black] {$0.25$} (-1.35, 2.6) -- cycle;
      \fill[sectionblue!22!white] (-0.6, 3.6) -- (-0.6, 4.4) -- (-0.9, 4.4) -- node[left,text=black] {$0.1$} (-0.9, 3.6) -- cycle;
      \fill[white] (-0.6, 4.6) -- node[left,text=black] {$0$} (-0.6, 5.4) -- (-1.65, 5.4) -- (-1.65, 4.6) -- cycle;
      
      \fill[sectionblue!44!white] (-0.4, -0.6) -- (0.4, -0.6) -- (0.4, -1.2) -- node[below,text=black] {$0.2$} (-0.4, -1.2) -- cycle;
      \fill[sectionblue] (0.6, -0.6) -- (1.4, -0.6) -- (1.4, -1.95) -- node[below,text=black] {$0.45$} (0.6, -1.95) -- cycle;
      \fill[sectionblue!11!white] (1.6, -0.6) -- (2.4, -0.6) -- (2.4, -0.75) -- node[below,text=black] {$0.05$} (1.6, -0.75) -- cycle;
      \fill[white] (2.6, -0.6) -- node[below,text=black] {$0$} (3.4, -0.6) -- (3.4, -1.35) -- (2.6, -1.35) -- cycle;
      \fill[sectionblue!11!white] (3.6, -0.6) -- (4.4, -0.6) -- (4.4, -0.75) -- node[below,text=black] {$0.05$} (3.6, -0.75) -- cycle;
      \fill[sectionblue!56!white] (4.6, -0.6) -- (5.4, -0.6) -- (5.4, -1.35) -- node[below,text=black] {$0.25$} (4.6, -1.35) -- cycle;

      \node[sectionblue,left] at (-2.1, 2.5) {$p$};
      \node[sectionblue,below] at (2.5, -1.8) {$q$};
      \node[linkred,above right] at (5.7, 5.7) {$\gamma^\ast$};

      \node[circle,draw=gray,text=gray] at (0,6) {$1$};
      \node[circle,draw=gray,text=gray] at (1,6) {$2$};
      \node[circle,draw=gray,text=gray] at (2,6) {$3$};
      \node[circle,draw=gray,text=gray] at (3,6) {$4$};
      \node[circle,draw=gray,text=gray] at (4,6) {$5$};
      \node[circle,draw=gray,text=gray] at (5,6) {$6$};
      \node[circle,draw=gray,text=gray] at (6,0) {$1$};
      \node[circle,draw=gray,text=gray] at (6,1) {$2$};
      \node[circle,draw=gray,text=gray] at (6,2) {$3$};
      \node[circle,draw=gray,text=gray] at (6,3) {$4$};
      \node[circle,draw=gray,text=gray] at (6,4) {$5$};
      \node[circle,draw=gray,text=gray] at (6,5) {$6$};

      \node[gray,above] at (2.5, 6.4) {states};
      \node[gray,right] at (6.4, 2.5) {states};
    \end{tikzpicture}
  \end{center}
  \caption{An optimal coupling $\gamma^\ast$ for $p = (0.35, \; 0.25, \; 0.05, \; 0.25, \; 0.1, \; 0)^{\transp}$ and
  $q = (0.2, \; 0.45, \; 0.05, \; 0, \; 0.05, \; 0.25)^{\transp}$}
  \label{tab:wdist_example_coupling}
\end{table}
For the given $p$ and $q$, an optimal transport plan, or an optimal
coupling $\gamma^\ast$ from \eqref{eq:wasserstein_def_finite}, is given
in \autoref{tab:wdist_example_coupling}. It follows that
\begin{align*}
  \WD{p}{q} = \sum_{r,s} \dist{r}{s} \cdot \gamma^\ast(r,s) = 0.975
\end{align*}
Note that $\gamma^\ast$ from \autoref{tab:wdist_example_coupling} does not satisfy
\autoref{prop:wasserstein_equiv} \ref{prop:optgammaonlyincorout}. Indeed,
for $r = 3$, $\gamma^\ast(3, 2) = 0.05 > 0$ and $\gamma^\ast(4, 3) = 0.05 > 0$
(and the condition is also violated for $r = 5$). However, we can apply the
method given in the proof of \autoref{prop:wasserstein_equiv} to turn
$\gamma^\ast$ into a coupling which satisfies \autoref{prop:wasserstein_equiv} \ref{prop:optgammaonlyincorout}.
With $r = 3$, $u = 2$ and $s = 4$, the proof tells us to remove $0.05$ mass
from both pairs $(r,u)$ and $(s,r)$, and to then add $0.05$ mass
to the pair $(s,u)$, i.e., to the pair $(4, 2)$. Hidden in the proof is
that we should also add $0.05$ mass to the diagonal if we want to keep
$\gamma^\ast$ a coupling. Repeating the procedure for $r = 5$, we arrive
at the coupling $\gamma^\ast$ in \autoref{tab:wdist_example_coupling_modified},
now satisfying \autoref{prop:wasserstein_equiv} \ref{prop:optgammaonlyincorout}.
\begin{table}[htb]
  \begin{center}
    \begin{tikzpicture}[>={Latex[length=1.5mm,width=1.5mm]},scale=0.85]
      \node[rectangle,minimum size=0.85cm,inner sep=0pt,draw=gray] at (5,5) {$0$};
      \node[rectangle,minimum size=0.85cm,inner sep=0pt,draw=gray] at (5,2) {$0$};
      \node[rectangle,minimum size=0.85cm,inner sep=0pt,draw=gray] at (5,1) {$0$};
      \node[rectangle,minimum size=0.85cm,inner sep=0pt,draw=gray] at (5,0) {$0$};
      \node[rectangle,minimum size=0.85cm,inner sep=0pt,draw=gray] at (4,5) {$0$};
      \node[rectangle,minimum size=0.85cm,inner sep=0pt,draw=gray] at (4,2) {$0$};
      \node[rectangle,minimum size=0.85cm,inner sep=0pt,draw=gray] at (4,1) {$0$};
      \node[rectangle,minimum size=0.85cm,inner sep=0pt,draw=gray] at (4,0) {$0$};
      \node[rectangle,minimum size=0.85cm,inner sep=0pt,draw=gray] at (3,5) {$0$};
      \node[rectangle,minimum size=0.85cm,inner sep=0pt,draw=gray] at (3,4) {$0$};
      \node[rectangle,minimum size=0.85cm,inner sep=0pt,draw=gray] at (3,3) {$0$};
      \node[rectangle,minimum size=0.85cm,inner sep=0pt,draw=gray] at (3,2) {$0$};
      \node[rectangle,minimum size=0.85cm,inner sep=0pt,draw=gray] at (3,1) {$0$};
      \node[rectangle,minimum size=0.85cm,inner sep=0pt,draw=gray] at (3,0) {$0$};
      \node[rectangle,minimum size=0.85cm,inner sep=0pt,draw=gray] at (2,5) {$0$};
      \node[rectangle,minimum size=0.85cm,inner sep=0pt,draw=gray] at (2,4) {$0$};
      \node[rectangle,minimum size=0.85cm,inner sep=0pt,draw=gray] at (2,1) {$0$};
      \node[rectangle,minimum size=0.85cm,inner sep=0pt,draw=gray] at (2,0) {$0$};
      \node[rectangle,minimum size=0.85cm,inner sep=0pt,draw=gray] at (1,5) {$0$};
      \node[rectangle,minimum size=0.85cm,inner sep=0pt,draw=gray] at (1,4) {$0$};
      \node[rectangle,minimum size=0.85cm,inner sep=0pt,draw=gray,fill=linkred,text=white] at (1,1) {$0.25$};
      \node[rectangle,minimum size=0.85cm,inner sep=0pt,draw=gray,fill=linkred!60!white] at (1,0) {$0.15$};
      \node[rectangle,minimum size=0.85cm,inner sep=0pt,draw=gray] at (0,5) {$0$};
      \node[rectangle,minimum size=0.85cm,inner sep=0pt,draw=gray] at (0,4) {$0$};
      \node[rectangle,minimum size=0.85cm,inner sep=0pt,draw=gray] at (0,3) {$0$};
      \node[rectangle,minimum size=0.85cm,inner sep=0pt,draw=gray] at (0,2) {$0$};
      \node[rectangle,minimum size=0.85cm,inner sep=0pt,draw=gray] at (0,1) {$0$};
      \node[rectangle,minimum size=0.85cm,inner sep=0pt,draw=gray,fill=linkred!80!white] at (0,0) {$0.2$};
      
      \node[rectangle,minimum size=0.85cm,inner sep=0pt,draw=black,line width=1pt,fill=linkred!20!white] at (5,4) {$0.05$};
      \node[rectangle,minimum size=0.85cm,inner sep=0pt,draw=black,line width=1pt,fill=linkred!80!white] at (5,3) {$0.2$};
      \node[rectangle,minimum size=0.85cm,inner sep=0pt,draw=black,line width=1pt,fill=linkred!20!white] at (4,4) {$0.05$};
      \node[rectangle,minimum size=0.85cm,inner sep=0pt,draw=black,line width=1pt] at (4,3) {$0$};
      \node[rectangle,minimum size=0.85cm,inner sep=0pt,draw=black,line width=1pt] at (2,3) {$0$};
      \node[rectangle,minimum size=0.85cm,inner sep=0pt,draw=black,line width=1pt,fill=linkred!20!white] at (2,2) {$0.05$};
      \node[rectangle,minimum size=0.85cm,inner sep=0pt,draw=black,line width=1pt,fill=linkred!20!white] at (1,3) {$0.05$};
      \node[rectangle,minimum size=0.85cm,inner sep=0pt,draw=black,line width=1pt] at (1,2) {$0$};

      \draw[line width=1pt,tumOrange] (0.4,-0.6) -- (1.6,-0.6) -- (1.6,0.6) -- (0.4,0.6) -- cycle;
      \draw[line width=1pt,tumOrange] (0.4,1.4) -- (1.6,1.4) -- (1.6,2.4) -- (2.6,2.4) -- (2.6,3.6) -- (0.4,3.6) -- cycle;
      \draw[line width=1pt,tumOrange] (3.4,2.4) -- (5.6,2.4) -- (5.6,4.6) -- (4.4,4.6) -- (4.4,3.6) -- (3.4,3.6) -- cycle;

      \node[left,sectionblue] at (-0.6,0) {$0.35$};
      \node[left,sectionblue] at (-0.6,1) {$0.25$};
      \node[left,sectionblue] at (-0.6,2) {$0.05$};
      \node[left,sectionblue] at (-0.6,3) {$0.25$};
      \node[left,sectionblue] at (-0.6,4) {$0.1$};
      \node[left,sectionblue] at (-0.6,5) {$0$};
      
      \node[below,sectionblue] at (0,-0.6) {$0.2$};
      \node[below,sectionblue] at (1,-0.6) {$0.45$};
      \node[below,sectionblue] at (2,-0.6) {$0.05$};
      \node[below,sectionblue] at (3,-0.6) {$0$};
      \node[below,sectionblue] at (4,-0.6) {$0.05$};
      \node[below,sectionblue] at (5,-0.6) {$0.25$};

      \node[sectionblue,left] at (-1.5, 2.5) {$p$};
      \node[sectionblue,below] at (2.5, -1.2) {$q$};
      \node[linkred,above] at (2.5, 5.7) {$\gamma^\ast$};

      \fill[tumGreen] (6.6,1.5) -- (7.4,1.5) -- (7.4,3.5) -- node[above,text=black] {$2$} (6.6,3.5) -- cycle;
      \node[above] at (8,1.5) {$0$};
      \fill[tumGreen] (8.6,1.5) -- (9.4,1.5) -- (9.4,2.5) -- node[above,text=black] {$1$} (8.6,2.5) -- cycle;
      \fill[tumGreen] (9.6,1.5) -- (10.4,1.5) -- (10.4,4) -- node[above,text=black] {$2.5$} (9.6,4) -- cycle;
      \fill[tumGreen] (10.6,1.5) -- (11.4,1.5) -- (11.4,2.5) -- node[above,text=black] {$1$} (10.6,2.5) -- cycle;
      \node[above] at (12,1.5) {$0$};
      \node[above,tumGreen] at (9.5, 5.7) {$f^\ast$};
      \draw[->,dashed,line width=1pt,tumOrange] (6.1, 5) -- (6.8, 5);
      \draw[dashed,line width=1pt,tumOrange] (6.8, 5) -- (7.1, 5) node[right] {maximally descending slopes};

      \draw[->,dashed,line width=1pt,tumOrange] (7,3.5) -- (7.5,2.5);
      \draw[dashed,line width=1pt,tumOrange] (7.5,2.5) -- (8,1.5);
      \draw[->,dashed,line width=1pt,tumOrange] (9,2.5) -- (8.5,2);
      \draw[dashed,line width=1pt,tumOrange] (8.5,2) -- (8,1.5);
      \draw[->,dashed,line width=1pt,tumOrange] (10,4) -- (9.5,3.25);
      \draw[dashed,line width=1pt,tumOrange] (9.5,3.25) -- (9,2.5);
      \draw[->,dashed,line width=1pt,tumOrange] (10,4) -- (10.5,3.25);
      \draw[dashed,line width=1pt,tumOrange] (10.5,3.25) -- (11,2.5);
      \draw[->,dashed,line width=1pt,tumOrange] (11,2.5) -- (11.5,2);
      \draw[dashed,line width=1pt,tumOrange] (11.5,2) -- (12,1.5);

      \node[circle,draw=gray,gray] at (7,1) {$1$};
      \node[circle,draw=gray,gray] at (8,1) {$2$};
      \node[circle,draw=gray,gray] at (9,1) {$3$};
      \node[circle,draw=gray,gray] at (10,1) {$4$};
      \node[circle,draw=gray,gray] at (11,1) {$5$};
      \node[circle,draw=gray,gray] at (12,1) {$6$};
      \node[below,gray] at (9.5,0.5) {states};
    \end{tikzpicture}
  \end{center}
  \caption{An optimal coupling $\gamma^\ast$ for $p = (0.35, \; 0.25, \; 0.05, \; 0.25, \; 0.1, \; 0)^{\transp}$ and
  $q = (0.2, \; 0.45, \; 0.05, \; 0, \; 0.05, \; 0.25)^{\transp}$ which
  also satisfies \autoref{prop:wasserstein_equiv} \ref{prop:optgammaonlyincorout}. Squares with black borders
  show changes to \autoref{tab:wdist_example_coupling}. The cost of the coupling / transport plan
  remains unchanged. On the right: the
  $f^\ast$ corresponding to $\gamma^\ast$ given by \autoref{prop:wasserstein_equiv}.
  In orange: the areas of maximally descending slope of $f^\ast$.}
  \label{tab:wdist_example_coupling_modified}
\end{table}

An optimal $f^\ast$ for \eqref{eq:wasserstein_dual_finite} which satisfies
(together with $\gamma^\ast$ from \autoref{tab:wdist_example_coupling_modified})
\autoref{prop:wasserstein_equiv} \ref{prop:optgammacoupling}-\ref{prop:optfgammarelation}
is given by
\begin{align*}
  f = (2, \; 0, \; 1, \; 2.5, \; 1, \; 0)^{\transp}
  \;\; \implies \;\;
  (p^{\transp} - q^{\transp}) f
  = (0.15, \; -0.2, \; 0, \; 0.25, \; 0.05, \; -0.25) \cdot f =
  0.975
\end{align*}
If $f^\ast$ is pictured as a height map, then the mass travels along descending
slopes of $f^\ast$ in the optimal transport plan $\gamma^\ast$ from
$p$ to $q$, and even only along slopes which are as steep as allowed
by the Lipschitz condition on $f^\ast$. $f^\ast$ is also shown in
\autoref{tab:wdist_example_coupling_modified}, together with the slopes along
which mass may travel in $\gamma^\ast$. Mass on the diagonal of $\gamma^\ast$
does not travel at all (which is allowed by the ``travel along steep slopes of $f^\ast$'' restriction),
and does not give rise to any cost.

\subsection{Ricci curvature}

The so-called Ricci curvature, which was originally defined from a geometric
point of view for a metric space \cite{riccicurvatureorig}, has been extended by \cite{riccimarkovmetricspaces}
to the setting of DTMCs.
\begin{definition}
  \label{def:ricci_curvature_dtmc}
  Given a DTMC with transition matrix $P \in \bbR^{n \times n}$, two states $r, s \in S = \{1, \ldots, n\}$,
  and a metric $\DIST$ on $S$, we define the coarse Ricci curvature of the DTMC along the states $r$ and $s$, with $r \neq s$, as
  (cf.\ \cite[Definition 3]{riccimarkovmetricspaces})
  \begin{align*}
    \kappa(r, s) := 1 - \frac{\WD{P_r}{P_s}}{\dist{r}{s}}
  \end{align*}
  where $P_r$ and $P_s$ are the $r$-th and $s$-th row of $P$, respectively.
  Furthermore, we define
  $\underline{\kappa}(P) := \min_{r \neq s} \kappa(r, s)$.
\end{definition}

In \cite[after Example 4 on page 814]{riccimarkovmetricspaces},
the extension to CTMCs is also briefly touched upon: Let $P^{(t)}(r, s) = \Prb{X_t = s \midd| X_0 = r}$
where $X_t$ is the state of the CTMC at time $t$. In particular, we have $\lim_{t \to 0} \frac{1}{t} P^{(t)}(r, s) = Q(r, s)$ for $r \neq s$.
Then, the Ricci curvature of the CTMC along the states $r$ and $s$, for $r \neq s$, is defined as
\begin{align*}
  \kappa(r, s) := - \left.\dufrac{t}\right|_{t = 0} \frac{\WD{P^{(t)}(r, \cdot)}{P^{(t)}(s, \cdot)}}{\dist{r}{s}}
\end{align*}
if this derivative exists. Of course, we have
$P^{(t)}(r, \cdot) = \delta_r^{\transp} e^{tQ}$ with $\delta_r \in \bbR^n$,
$\delta_r(r) = 1$ and $\delta_r(s) = 0$ for $s \neq r$ in our notation.
Hence, we define Ricci curvature as follows:
\begin{definition}
  \label{def:ricci_curvature}
  Given a CTMC with generator $Q \in \bbR^{n \times n}$, two states $r, s \in S = \{1, \ldots, n\}$,
  and a metric $\DIST$ on $S$, we define the coarse Ricci curvature of the CTMC along the states $r$ and $s$, with $r \neq s$, as
  \begin{align*}
    \kappa(r, s) := -\frac{1}{\dist{r}{s}} \cdot \left.\dufrac{t}\right|_{t = 0^+} \WD{\delta_r^{\transp} e^{tQ}}{\delta_s^{\transp} e^{tQ}}
  \end{align*}
  with $\delta_r, \delta_s \in \bbR^n$ being the
  Dirac measures concentrated on $r$ and $s$, respectively.
  The derivative exists by \autoref{lem:wdderivative} and \autoref{cor:wdderiv_error}.
  We also set $\underline{\kappa}(Q) := \min_{r \neq s} \kappa(r, s)$.
\end{definition}
The concept of Ricci curvature will help us bound the error caused by
our aggregation scheme.

\section{Wasserstein error bounds}
\label{sec:w_error_bounds}

\subsection{The CTMC case}
\label{ssec:w_error_bounds_ctmc}

Recall that $p_0 \in \bbR^n$ is the initial distribution,
that the transient distribution of the CTMC is given by $p_t^{\transp} = p_0^{\transp} e^{t Q}$,
and that we approximate $p_t$ by $\widetilde{p}_t$, defined as
\begin{align*}
  \widetilde{p}_t^{\transp} = \pi_0^{\transp} e^{t \Theta} A \;\;\;
  \textrm{ with } \pi_0 \in \bbR^m, \Theta \in \bbR^{m \times m}, A \in \bbR^{m \times n}
\end{align*}

Our goal is to prove \autoref{thm:wasserstein_error_growth_bound} below, which
bounds the rate at which the error $\WD{\widetilde{p}_t}{p_t}$ can grow, and thereby
lets us bound the error at any point in time.
\autoref{thm:wasserstein_error_growth_bound} is a generalization
of \cite[Theorem 5]{formalbndsstatespaceredmc} (at least in the setting
where our reduced model on the lower-dimensional state space is also
a Markov chain). For the proof of \autoref{thm:wasserstein_error_growth_bound}, we
will split the error growth into two classes, which are treated in
\autoref{sssec:accum_error_growth} and \autoref{sssec:approx_error_growth}:
first, we consider how the error accumulated up to a given time point
will propagate, and then, we look at the error caused by the approximation
of the dynamics on a lower-dimensional state space. For a bound on the accumulated
error propagation, we will rely on the Ricci curvature from \autoref{def:ricci_curvature}.

We start with a general result which gives us a way to calculate the
derivative of the Wasserstein distance between two probability distributions
which depend on a time parameter.

\begin{lemma}
  \label{lem:wdderivative}
  Let $p_u, q_u \in \bbR^n$ be probability measures depending on a parameter $u \geq 0$.
  Further assume that $p_u$ and $q_u$ are continuous for $u \geq 0$, that $p_u$ and $q_u$ have one-sided right derivatives for $u \geq 0$
  which are locally bounded near $0^+$, and denote the (one-sided)
  derivatives in $u$ by $\dot{p}_u$ and $\dot{q}_u$.
  Then, the one-sided derivative of $\WD{p_u}{q_u}$ at $u = 0^+$ exists and
  \begin{align*}
    \left.\dufrac{u}\right|_{u = 0^+} \WD{p_u}{q_u} &= \max_{f \in M} \big(\dot{p}_0^{\transp} -  \dot{q}_0^{\transp}\big) f \\
    \textrm{where } M &:= \argmax_{\substack{f \in \bbR^n \textrm{ is }1\textrm{-Lipschitz w.r.t.\ }\DIST\\\forall s \in S: 0 \leq f(s) \leq d_{\max}}} \; \big(p_0^{\transp} - q_0^{\transp}\big) f \qedhere
  \end{align*}
\end{lemma}

\begin{tproof}
  By \eqref{eq:wasserstein_dual_finite}, we have
  \begin{align}
    \WD{p_u}{q_u}
    &= \max_{\substack{f \in \bbR^n \textrm{ is }1\textrm{-Lipschitz w.r.t.\ }\DIST\\\forall s \in S: 0 \leq f(s) \leq d_{\max}}} \; (p_u^{\transp} f - q_u^{\transp} f)
    \label{eq:danskin_max}
  \end{align}
  We will use Danskin's Theorem \cite[Theorem I on page 22]{theorymaxmin}. In particular,
  we use the version proven in \cite{danskinappendix}, which requires that
  the maximization in \eqref{eq:danskin_max} is over a compact subset of a Banach space. Indeed, the set $\calV$
  of all $f \in \bbR^n$ which are $1$-Lipschitz w.r.t.\ $\DIST$ and which
  satisfy $0 \leq f(s) \leq d_{\max}$ for all $s$ is clearly a compact subset
  of the vector space $\bbR^n$ with the Euclidean norm, which is a Banach space.
  We further have to verify the three hypotheses from \cite{danskinappendix}:
  \begin{itemize}
    \item \textbf{H1}. The map $\bbR^n \ni f \mapsto (p_u^{\transp} f - q_u^{\transp} f) \in \bbR$ is clearly continuous
    with respect to the Euclidean topology.
    \item \textbf{H2}. For all $f \in \calV$ and for all $u \in [0, \varepsilon)$ (for some $\varepsilon > 0$),
    the one-sided derivative
    \begin{align*}
      \dufrac{u^+} (p_u^{\transp} f - q_u^{\transp} f)
      = \dufrac{u^+} (p_u^{\transp} - q_u^{\transp}) f
    \end{align*}
    clearly exists: it is equal to the linear combination with weights $f(s)$ of
    $\dot{p}_u - \dot{q}_u$, which we assumed to exist. The derivatives are locally bounded
    by assumption.
    \item \textbf{H3}. The map $(u, f) \mapsto (p_u^{\transp} f - q_u^{\transp} f)$
    is clearly continuous: it is linear in $f$, continuous in $u$ because $p_u$ and
    $q_u$ are continuous in $u$ by assumption, and $(a, b) \mapsto a^{\transp}b$ is
    a continuous map as well.
  \end{itemize}
  Therefore, by \cite[Theorem 10.1]{danskinappendix}, the right derivative (or one-sided derivative at $u = 0^+$)
  of $\WD{p_u}{q_u}$ exists and is given by
  \begin{align*}
    \left.\dufrac{u}\right|_{u = 0^+} \WD{p_u}{q_u}
    &= \max_{f \in M} \left.\dufrac{u}\right|_{u = 0^+} (p_u^{\transp} f - q_u^{\transp} f)
    = \max_{f \in M} (\dot{p}_0^{\transp} -  \dot{q}_0^{\transp}) f \\
    \textrm{where } M &:= \argmax_{\substack{f \in \bbR^n \textrm{ is }1\textrm{-Lipschitz w.r.t.\ }\DIST\\\forall s \in S: 0 \leq f(s) \leq d_{\max}}} \; (p_0^{\transp} f - q_0^{\transp} f)
    \qedhere
  \end{align*}
\end{tproof}

\begin{remark}
  In the subsequent applications of \autoref{lem:wdderivative}, we will not
  mention the assumption of the locally bounded right derivatives anymore. However,
  this assumption does indeed hold when we apply \autoref{lem:wdderivative} later.
  As an example, the right derivative (in $t$) of all components of $p^{\transp} e^{tQ}$
  (for $p \in \bbR^n$ a probability measure and $Q \in \bbR^{n \times n}$ a generator
  matrix) is bounded for all $t \geq 0$ by the maximal exit rate of $Q$. This
  argumentation can be extended to the cases where we apply \autoref{lem:wdderivative} below.
  We also do not explicitly mention continuity and existence of the right
  derivatives whenever these assumptions are straightforward to show.
\end{remark}

\subsubsection{Bounding the growth of the accumulated errror}
\label{sssec:accum_error_growth}

We now consider how to bound the growth of the error which has already accumulated up
to the current time point of the calculation in the aggregation scheme.
Given the actual transient distribution $p_t$ and its approximation $\widetilde{p}_t$,
we will look at how the Wasserstein distance between the two would develop
if both distributions were to evolve according to the original generator $Q$ from time
$t$ onwards. That is, we ignore the approximation of the dynamics and just look at
the way in which the accumulated error propagates.

A first step in that direction is the following direct corollary of \autoref{lem:wdderivative}:
\begin{corollary}
  \label{cor:wdderiv_error}
  Let $p, q \in \bbR^n$ be probability measures, and let $Q \in \bbR^{n \times n}$ be a generator matrix.
  Then,
  \begin{align*}
    \left.\dufrac{u}\right|_{u = 0^+} \WD{p^{\transp} e^{uQ}}{q^{\transp} e^{uQ}}
    &= \max_{f \in M} \; (p^{\transp} -  q^{\transp}) Q f
    \leq \;\; \max_{\substack{f \in \bbR^n \textrm{ is }1\textrm{-Lipschitz w.r.t.\ }\DIST\\\forall s \in S: 0 \leq f(s) \leq d_{\max}}} \; (p^{\transp} -  q^{\transp}) Q f \\
    \textrm{with }  M &= \argmax_{\substack{f \in \bbR^n \textrm{ is }1\textrm{-Lipschitz w.r.t.\ }\DIST\\\forall s \in S: 0 \leq f(s) \leq d_{\max}}} \; (p^{\transp} - q^{\transp})f
    \qedhere
  \end{align*}
\end{corollary}

\autoref{cor:wdderiv_error} gives us a way to calculate the coarse
Ricci curvature from \autoref{def:ricci_curvature} with a linear program,
which is helpful for applications, but also for the subsequent theory.

\begin{lemma}
  \label{lem:ricci_linearprogram}
  For $r \neq s$, we have that
  \begin{align*}
    &\kappa(r, s) = -\frac{1}{\dist{r}{s}} \cdot V \quad \textrm{ where } V \textrm{ is the solution of } \\
    &\max_{f \in \bbR^n, f \geq 0} (Q_r - Q_s)f \quad \textrm{s.t.} \quad
    f(r) - f(s) = \dist{r}{s}
    \textrm{ and }
    \forall \widetilde{r}, \widetilde{s}: f(\widetilde{r}) - f(\widetilde{s}) \leq \dist{\widetilde{r}}{\widetilde{s}}
  \end{align*}
  where $Q_r$ and $Q_s$ are the $r$-th and $s$-th row of $Q$.
\end{lemma}

\begin{tproof}
  Let $\delta_r, \delta_s \in \bbR^n$ be the Dirac probability measures concentrated
  on $r$ and $s$, respectively, and let $Q$ be a generator matrix.
  By \autoref{cor:wdderiv_error}, we have
  \begin{align*}
    \left.\dufrac{u}\right|_{u = 0^+} \WD{\delta_r^{\transp} e^{uQ}}{\delta_s^{\transp} e^{uQ}}
    &= \max_{f \in M} (\delta_r^{\transp} - \delta_s^{\transp})Qf
    = \max_{f \in M} (Q_r - Q_s)f 
  \end{align*}
  where $Q_r$ and $Q_s$ are again the $r$-th and $s$-th row of $Q$, and where
  \begin{align*}
    M &= \argmax_{\substack{f \in \bbR^n \textrm{ is }1\textrm{-Lipschitz w.r.t.\ }\DIST\\\forall \widetilde{s} \in S: 0 \leq f(\widetilde{s}) \leq d_{\max}}} \; \underbrace{(\delta_r^{\transp} - \delta_s^{\transp})f}_{= f(r) - f(s)} \\
    &= \left\{ f \in \bbR^n : \; f \textrm{ is }1\textrm{-Lip.\ w.r.t.\ }\DIST, \; \forall \widetilde{s} \in S: 0 \leq f(\widetilde{s}) \leq d_{\max}, \; f(r) - f(s) = \dist{r}{s}\right\}
  \end{align*}
  Dropping the restriction to $f$ with $f(\widetilde{s}) \leq d_{\max}$ does not change anything
  by the remark after \eqref{eq:wasserstein_dual_finite}.
\end{tproof}

The next lemma shows how the coarse Ricci curvature from \autoref{def:ricci_curvature}
can be used to bound the rate at which the Wasserstein distance between two
transient distributions of a CTMC grows. This will later help us to bound
the rate at which the accumulated error continues to grow in our aggregation
scheme.

\begin{lemma}
  \label{lem:wdderiv_error_ricciinfbound}
  Let $p, q \in \bbR^n$ be probability measures, and let $Q$ be a generator matrix.
  Then,
  \begin{align*}
    \left.\dufrac{u}\right|_{u = 0^+} \WD{p^{\transp} e^{uQ}}{q^{\transp} e^{uQ}}
    &\leq -\underline{\kappa}(Q) \cdot \WD{p}{q}
  \end{align*}
  where $\underline{\kappa}(Q)$ was defined in \autoref{def:ricci_curvature}.
\end{lemma}

\begin{tproof}
  This is essentially a corollary of \cite[Theorem 1.9]{markovprocricci}. However, we include a
  proof specifically for our setting.

  By \autoref{cor:wdderiv_error}, we have
  \begin{align*}
    \left.\dufrac{u}\right|_{u = 0^+} \WD{p^{\transp} e^{uQ}}{q^{\transp} e^{uQ}}
    &= \max_{f \in M} (p^{\transp} -  q^{\transp}) Q f
    \textrm{ with }  M &= \argmax_{\substack{f \in \bbR^n \textrm{ is }1\textrm{-Lipschitz w.r.t.\ }\DIST\\\forall s \in S: 0 \leq f(s) \leq d_{\max}}} (p^{\transp} - q^{\transp})f
  \end{align*}
  Assume that $\gamma$ is a coupling
  achieving the Wasserstein distance between $p$ and $q$,
  that is, $\WD{p}{q} = \sum_{r,s} \dist{r}{s}\gamma(r, s)$.
  Indeed, we choose a $\gamma$ of the form given in
  \autoref{prop:wasserstein_equiv} \ref{prop:optgammacoupling}--\ref{prop:optfgammarelation}.
  We have
  \begin{align*}
    (p^{\transp} -  q^{\transp}) Q f &= \sum_r \big(p(r) - q(r)\big) \cdot (Qf)(r) = \sum_{r, s} \gamma(r,s) \cdot \big((Qf)(r) - (Qf)(s)\big) \\
    &= \sum_{r \neq s} \gamma(r,s) \cdot \big((Qf)(r) - (Qf)(s)\big)
  \end{align*}
  Hence,
  \begin{align*}
    \max_{f \in M} \; (p^{\transp} -  q^{\transp}) Q f
    &\leq
    \sum_{r \neq s} \gamma(r,s) \cdot \max_{f \in M} \big((Qf)(r) - (Qf)(s)\big)
    = \sum_{r \neq s} \gamma(r,s) \cdot \max_{f \in M} (Q_r - Q_s) f
  \end{align*}
  Now, consider the set $M$, which is the set of optimal solutions
  of \eqref{eq:wasserstein_finite_primallp} (actually a subset due to the
  restriction $\leq d_{\max}$). As we did already in the proof of
  \autoref{prop:wasserstein_equiv} \ref{prop:optfgammarelation}, we can
  invoke complementary slackness \cite[Theorem 5.3]{linearprogramming}
  to see that:
  \begin{align*}
    \forall f \in M: \gamma(r,s) > 0 \implies f(r) - f(s) = \dist{r}{s}
  \end{align*}
  This implies
  \begin{align*}
    \forall r \neq s: \gamma(r,s) \cdot \max_{f \in M} (Q_r - Q_s) f
    \leq \gamma(r,s) \cdot \underbrace{\left(\max_{\substack{f \in \bbR^n \textrm{ is }1\textrm{-Lipschitz w.r.t.\ }\DIST\\\forall \widetilde{s} \in S: 0 \leq f(\widetilde{s}) \leq d_{\max},\; f(r) - f(s) = \dist{r}{s}}} (Q_r - Q_s) f\right)}_{\displaystyle =: V(r,s)}
  \end{align*}
  By \autoref{lem:ricci_linearprogram} and its proof, we have
  $V(r, s) = -\dist{r}{s} \cdot \kappa(r, s)$. We can therefore conclude that
  \begin{align}
    \label{eq:ricci_bound_sharp}
    \max_{f \in M} \; (p^{\transp} -  q^{\transp}) Q f
    &\leq \sum_{r \neq s} \gamma(r,s) \dist{r}{s} \cdot \big(-\kappa(r, s)\big) \\
    &\leq \sum_{r \neq s} \gamma(r,s) \dist{r}{s} \cdot \left(\max_{r \neq s} -\kappa(r, s)\right)
    = \WD{p}{q} \cdot (-\underline{\kappa}(Q)) \notag
  \end{align}
  Note that \eqref{eq:ricci_bound_sharp} gives a sharper bound, which we will
  use from time to time instead of the final bound relying on $\underline{\kappa}(Q)$.

  \textbf{\textit{Alternative proof.}} Let $\gamma$ be a coupling
  achieving the Wasserstein distance between $p$ and $q$,
  that is, $\WD{p}{q} = \sum_{r,s} \dist{r}{s}\gamma(r, s)$.
  For every $u \geq 0$ and all state pairs $r,s$, let $\eta_u^{(r,s)}$ be the
  coupling of $\delta_r^{\transp} e^{uQ}$ and $\delta_s^{\transp} e^{uQ}$ achieving
  the Wasserstein distance between the two distributions. Then,
  $\beta_u := \sum_{r, s} \gamma(r, s) \cdot \eta_u^{(r,s)}$ is a coupling
  between $p^{\transp} e^{uQ}$ and $q^{\transp} e^{uQ}$. Thus,
  \begin{align*}
    \WD{p^{\transp} e^{uQ}}{q^{\transp} e^{uQ}} &\leq \sum_{i, j} \beta_u(i, j) \cdot \dist{i}{j}
    = \sum_{i, j} \sum_{r, s} \gamma(r, s) \cdot \eta_u^{(r,s)}(i, j) \cdot \dist{i}{j} \\
    &= \sum_{r, s} \gamma(r, s) \cdot \sum_{i, j} \eta_u^{(r,s)}(i, j) \cdot \dist{i}{j}
    = \sum_{r, s} \gamma(r, s) \cdot \WD{\delta_r^{\transp} e^{uQ}}{\delta_s^{\transp} e^{uQ}}
  \end{align*}
  Differentiating, we obtain (note that the inequality above is an equality for $u = 0$)
  \begin{align*}
    \left.\dufrac{u}\right|_{u = 0^+} \WD{p^{\transp} e^{uQ}}{q^{\transp} e^{uQ}}
    &\leq \left.\dufrac{u}\right|_{u = 0^+} \sum_{r, s} \gamma(r, s) \cdot \WD{\delta_r^{\transp} e^{uQ}}{\delta_s^{\transp} e^{uQ}} \\
    &= \sum_{r \neq s} \gamma(r, s) \cdot \underbrace{\left.\dufrac{u}\right|_{u = 0^+}\WD{\delta_r^{\transp} e^{uQ}}{\delta_s^{\transp} e^{uQ}}}_{\displaystyle \dist{r}{s} \cdot \big(-\kappa(r, s)\big)}
  \end{align*}
  The existence of the derivatives follows from \autoref{lem:wdderivative}.
  The proof can then be finished as shown in the line after \eqref{eq:ricci_bound_sharp}.
\end{tproof}

The following lemma provides a lower bound for $\kappa$ (i.e., an upper bound
for the derivative of the Wasserstein distance between two transient distributions).
The bound is straightforward to compute with simple matrix-vector multiplications,
and therefore computationally less expensive than solving the linear program
from \autoref{lem:ricci_linearprogram}.

\begin{lemma}
  \label{lem:ricci_lower_bound}
  For $r \neq s$, it holds that
  \begin{align*}
    \kappa(r, s) \geq -\frac{\min\left\{Q_{r} \dist{r}{\cdot}, \;\; Q_{r} \dist{s}{\cdot}\right\} + \min\left\{Q_{s} \dist{s}{\cdot}, \;\; Q_{s} \dist{r}{\cdot}\right\}}{\dist{r}{s}} =: k(r, s)
  \end{align*}
  where $Q_r$ is the $r$-th row of $Q$ and $\dist{r}{\cdot} := \big(\dist{r}{1}, \ldots, \dist{r}{n}\big)^{\transp}$.
\end{lemma}

\begin{tproof}
  By \autoref{lem:ricci_linearprogram}, we have
  \begin{align*}
    \kappa(r,s) &= -\frac{\max_{f \in M} (Q_r - Q_s)f}{\dist{r}{s}} \\
    \textrm{with } M &= \left\{ f \in \bbR^n : \; f \textrm{ is }1\textrm{-Lip.\ w.r.t.\ }\DIST, \; \forall \widetilde{s} \in S: 0 \leq f(\widetilde{s}), \; f(r) - f(s) = \dist{r}{s}\right\}
  \end{align*}
  Hence, it suffices to show that
  \begin{align*}
    \max_{f \in M} (Q_r - Q_s)f \leq \min\left\{Q_{r} \dist{r}{\cdot}, \;\; Q_{r} \dist{s}{\cdot}\right\} + \min\left\{Q_{s} \dist{s}{\cdot}, \;\; Q_{s} \dist{r}{\cdot}\right\}
  \end{align*}
  Indeed, we have, for arbitrary $\widetilde{r},\widetilde{s} \in S$,
  \begin{align}
    (Q_r - Q_s)f &= (Qf)(r) - (Qf)(s)
    = \sum_k Q(r,k) f(k) - \sum_k Q(s,k) f(k) \notag \\
    &\overset{\circledast}{=} \sum_k Q(r,k) \big(f(k) - f(\widetilde{r})\big) + \sum_k Q(s,k) \big(f(\widetilde{s}) - f(k)\big) \notag \\
    &= Q(r,r) \big(f(r) - f(\widetilde{r})\big) + Q(s,s) \big(f(\widetilde{s}) - f(s)\big) \notag \\
    &\hphantom{\;=\;} +\sum_{k \neq r} \underbrace{Q(r,k)}_{\geq 0} \big(f(k) - f(\widetilde{r})\big) + \sum_{k \neq s} \underbrace{Q(s,k)}_{\geq 0} \big(f(\widetilde{s}) - f(k)\big) \notag \\
    &\hskip-0.53cm \overset{\textrm{for }f \textrm{ 1-Lip.}}{\leq} Q(r,r) \big(f(r) - f(\widetilde{r})\big) + Q(s,s) \big(f(\widetilde{s}) - f(s)\big) \notag \\
    &\hskip-0.53cm \hphantom{\;\overset{\textrm{for }f \textrm{ 1-Lip.}}{\leq}\;} + \sum_{k \neq r} Q(r,k) \dist{k}{\widetilde{r}} + \sum_{k \neq s} Q(s,k) \dist{\widetilde{s}}{k} \notag \\
    &= Q_r \dist{\widetilde{r}}{\cdot} + Q_s \dist{\widetilde{s}}{\cdot} \notag \\
    &\hphantom{\;=\;}
    + Q(r,r) \big(f(r) - f(\widetilde{r}) - \dist{\widetilde{r}}{r}\big)
    + Q(s,s) \big(f(\widetilde{s}) - f(s) - \dist{\widetilde{s}}{s}\big) \label{eq:extra_term}
  \end{align}
  where $\circledast$ holds because each row of $Q$ sums to $0$. We can now
  finish the proof by considering all four possible combinations of $\widetilde{r} \in \{r, s\}$ and
  $\widetilde{s} \in \{r, s\}$ and by noting that the extra term in \eqref{eq:extra_term}
  always disappears if $f \in M$: for $\widetilde{r} = r$, we have
  $f(r) - f(\widetilde{r}) - \dist{\widetilde{r}}{r} = f(r) - f(r) - \dist{r}{r} = 0$;
  for $\widetilde{r} = s$, we have
  $f(r) - f(\widetilde{r}) - \dist{\widetilde{r}}{r} = f(r) - f(s) - \dist{r}{s} = 0$
  because $f(r) - f(s) = \dist{r}{s}$ for all $f \in M$; and the same argumentation
  applies to $Q(s,s) \big(f(\widetilde{s}) - f(s) - \dist{\widetilde{s}}{s}\big)$.
\end{tproof}

\begin{corollary}
  \label{cor:wdderiv_error_bound}
  Let $p, q \in \bbR^n$ be probability measures, and let $Q \in \bbR^{n \times n}$ be a generator matrix.
  Then,
  \begin{align*}
    \left.\dufrac{u}\right|_{u = 0^+} \WD{p^{\transp} e^{uQ}}{q^{\transp} e^{uQ}}
    &\begin{cases}
      \displaystyle \leq -\underline{k}(Q) \cdot \WD{p}{q} \qquad \textrm{ with } \underline{k}(Q) := \min_{r \neq s} k(r, s) \\
      \displaystyle \leq K(Q) := \max\left\{ 0, \; -\min_{r \neq s} \dist{r}{s} \cdot k(r, s) \right\}
    \end{cases}
  \end{align*}
  with $k(r,s)$ defined in \autoref{lem:ricci_lower_bound}.
\end{corollary}

\begin{tproof}
  By \autoref{lem:ricci_lower_bound}, $\kappa(r, s) \geq k(r, s)$ for $r \neq s$.
  This implies $\underline{\kappa}(Q) = \min_{r \neq s} \kappa(r, s) \geq \min_{r \neq s} k(r, s) = \underline{k}(Q)$ and thus,
  by \autoref{lem:wdderiv_error_ricciinfbound}, we have
  \begin{align*}
    \left.\dufrac{u}\right|_{u = 0^+} \WD{p^{\transp} e^{uQ}}{q^{\transp} e^{uQ}}
    &\leq -\underline{\kappa}(Q) \cdot \WD{p}{q} \leq -\underline{k}(Q) \cdot \WD{p}{q}
  \end{align*}
  For the second bound, we use \eqref{eq:ricci_bound_sharp} which implies,
  with $\gamma$ being an optimal coupling achieving the Wasserstein distance
  between $p$ and $q$,
  \begin{align*}
    \left.\dufrac{u}\right|_{u = 0^+} \WD{p^{\transp} e^{uQ}}{q^{\transp} e^{uQ}}
    &\leq \sum_{r \neq s} \gamma(r,s) \dist{r}{s} \cdot \big(-\kappa(r, s)\big) \\
    &\overset{\textrm{\autoref{lem:ricci_lower_bound}}}{\leq} \sum_{r \neq s} \gamma(r,s) \cdot \big(-\dist{r}{s} \cdot k(r, s)\big) \\
    &\leq \left(\sum_{r \neq s} \gamma(r,s)\right) \cdot \left(-\min_{r \neq s} \dist{r}{s} \cdot k(r, s)\right) \\
    &\leq \max\left\{ 0, \; -\min_{r \neq s} \dist{r}{s} \cdot k(r, s) \right\}
  \end{align*}
  In the last inequality, we have to insert the maximum of $0$ and the term
  from the previous line because it could be that $-\min_{r \neq s} \ldots < 0$
  (and it typically holds that $\sum_{r \neq s} \gamma(r, s) < 1$).
\end{tproof}

For the discrete metric, we can simplify the expression for $k(r, s)$:
\begin{lemma}
  \label{lem:krs_discrmet}
  Let $Q$ be a generator matrix.
  If $\DIST$ is the discrete metric, then, for $r \neq s$,
  \begin{align*}
    k(r, s) = Q(r, s) + Q(s, r)
  \end{align*}
  where $k(r, s)$ was defined in \autoref{lem:ricci_lower_bound}.
\end{lemma}

\begin{tproof}
  Note that, for the discrete metric, we have, for $r \neq s$,
  \begin{align*}
    Q_r \dist{r}{\cdot} = -Q(r, r) \geq 0, \quad
    Q_r \dist{s}{\cdot} &= -Q(r, s) \leq 0, \\
    Q_s \dist{s}{\cdot} = -Q(s, s) \geq 0, \quad
    Q_s \dist{r}{\cdot} &= -Q(s, r) \leq 0 \\
    \implies \min\left\{Q_{r} \dist{r}{\cdot}, \;\; Q_{r} \dist{s}{\cdot}\right\} &= -Q(r, s), \\
    \min\left\{Q_{s} \dist{s}{\cdot}, \;\; Q_{s} \dist{r}{\cdot}\right\} &= -Q(s, r) \\
    \implies k(r, s) &= -\frac{-Q(r,s)-Q(s,r)}{\dist{r}{s}} = Q(r, s) + Q(s, r) \qedhere
  \end{align*}
\end{tproof}

We can now show that the Wasserstein distance between
transient distributions is necessarily non-increasing for the discrete metric.

\begin{corollary}
  \label{cor:wdderiv_error_discrmet}
  Let $p, q \in \bbR^n$ be probability measures, and let $Q$ be a generator matrix.
  If $\DIST$ is the discrete metric, then
  \begin{align*}
    \left.\dufrac{u}\right|_{u = 0^+} \WD{p^{\transp} e^{uQ}}{q^{\transp} e^{uQ}}
    &\leq 0 \qedhere
  \end{align*}
\end{corollary}

\begin{tproof}
  For the discrete metric and for $r \neq s$,
  it holds by \autoref{lem:krs_discrmet} that
  $k(r, s) = Q(r, s) + Q(s, r) \geq 0$
  which implies that $-\underline{k}(Q) = \left(-\min_{r \neq s} k(r, s)\right) \leq 0$.
  By \autoref{cor:wdderiv_error_bound}, we therefore have
  \begin{align*}
    \left.\dufrac{u}\right|_{u = 0^+} \WD{p^{\transp} e^{uQ}}{p^{\transp} e^{uQ}}
    &\leq -\underline{k}(Q) \cdot \WD{p}{q} \leq 0 \qedhere
  \end{align*}
\end{tproof}

\subsubsection{Bounding the approximation error}
\label{sssec:approx_error_growth}

Up to now, we have seen how to bound the rate of growth of the Wasserstein distance
between two transient distributions of the same CTMC (i.e., with the same generator).
This helps us to analyze how the accumulated error in our aggregation scheme might
blow up (or even decrease over time).
To provide complete error bounds, we also need to consider the error
caused by approximating the generator of the original CTMC on a lower-dimensional state space.
The following corollary of \autoref{lem:wdderivative} fills that gap:
\begin{corollary}
  \label{cor:wdderiv_approxdyn}
  Let $\pi \in \bbR^m$ be a probability measure, and let $\Theta \in \bbR^{m \times m}, Q \in \bbR^{n \times n}$ be generator matrices.
  Further, let $A \in \bbR^{m \times n}$ be a matrix with non-negative entries and with each row summing to $1$.
  Then,
  \begin{align*}
    \left.\dufrac{u}\right|_{u = 0^+} \WD{\pi^{\transp} e^{u\Theta} A}{\pi^{\transp} A e^{uQ}}
    &= \max_{\substack{f \in \bbR^n \textrm{ is }1\textrm{-Lipschitz w.r.t.\ }\DIST\\\forall s \in S: 0 \leq f(s) \leq d_{\max}}} \; \pi^{\transp} (\Theta A - A Q) f \qedhere
  \end{align*}
\end{corollary}

\begin{tproof}
  Note that the set $M$ from \autoref{lem:wdderivative} is the set of all 1-Lipschitz functions in this
  case because the Wasserstein distance between the two probability measures $\pi^{\transp} e^{u\Theta} A$
  and $\pi^{\transp} A e^{uQ}$ is $0$ for $u = 0$.
\end{tproof}

We can immediately derive an upper bound for the derivative in \autoref{cor:wdderiv_approxdyn}
which is easier to compute:
\begin{corollary}
  \label{cor:wdderiv_approxdyn_bound}
  Assume that $\Theta \in \bbR^{m \times m}, Q \in \bbR^{n \times n}, A \in \bbR^{m \times n}$ are fixed. Then,
  \begin{align*}
    \forall \pi \in \bbR^m \textrm{ prob.\ measure}: \;\;
    \left.\dufrac{u}\right|_{u = 0^+} \WD{\pi^{\transp} e^{u\Theta} A}{\pi^{\transp} A e^{uQ}}
    &\leq \pi^{\transp} \dabs{\Theta A - A Q}_{\mathrm{W}}
  \end{align*}
  where $\dabs{\cdot}_{\mathrm{W}}$ was defined in \autoref{def:row_wd_norm}.
\end{corollary}

\subsubsection{Overall error bound}

\begin{theorem}
  \label{thm:wasserstein_error_growth_bound}
  Consider an initial distribution $p_0 \in \bbR^n$ of a CTMC with generator $Q \in \bbR^{n \times n}$
  and transient distribution $p_t^{\transp} = p_0^{\transp} e^{t Q}$,
  and consider the approximation $\widetilde{p}_t$, defined as
  \begin{align*}
    \widetilde{p}_t^{\transp} = \pi_0^{\transp} e^{t \Theta} A \;\;\;
    \textrm{ with } \pi_0 \in \bbR^m, \Theta \in \bbR^{m \times m}, A \in \bbR^{m \times n}
  \end{align*}
  where $\pi_0$ is our approximated initial distribution on a lower-dimensional state space,
  $\Theta$ is the generator matrix for the CTMC on the lower dimensional state space, and
  $A$ is the disaggregation matrix (non-negative entries, each row sums to $1$).

  Then, $\WD{\widetilde{p}_{t}}{p_{t}}$ is continuous in $t$, and
  \begin{align*}
    \dufrac{t^+} \WD{\widetilde{p}_{t}}{p_{t}}
    &\leq \pi_t^{\transp} \dabs{\Theta A - A Q}_{\mathrm{W}}
    + K(Q) \\
    \textrm{and } \;\; \dufrac{t^+} \WD{\widetilde{p}_{t}}{p_{t}}
    &\leq \pi_t^{\transp} \dabs{\Theta A - A Q}_{\mathrm{W}}
    + \WD{\widetilde{p}_{t}}{p_{t}} \cdot \big(-\underline{k}(Q)\big)
  \end{align*}
  where $\pi_t^{\transp} = \pi_0^{\transp} e^{t \Theta}$, and where $t^+$ indicates
  that we consider the one-sided derivative into the positive $t$-direction
  (the right derivative). This derivative exists for all $t \geq 0$. $K(Q)$ and $\underline{k}(Q)$ are defined
  in \autoref{cor:wdderiv_error_bound}. We can also replace $\underline{k}(Q)$ with
  $\underline{\kappa}(Q)$ from \autoref{def:ricci_curvature} without affecting the validity of the bound above.
\end{theorem}

\begin{tproof}
  $\WD{\widetilde{p}_{t}}{p_{t}}$ is continuous in $t$ because $\widetilde{p}_{t}$ and $p_t$ are continuous in $t$.
  Indeed, by the triangle inequality for the Wasserstein distance,
  \begin{align*}
    \WD{\widetilde{p}_{t + u}}{p_{t + u}} &\leq \WD{\widetilde{p}_{t+u}}{\widetilde{p}_{t}} + \WD{\widetilde{p}_{t}}{p_{t}} + \WD{p_{t}}{p_{t+u}} \\
    \WD{\widetilde{p}_{t}}{p_{t}} &\leq \WD{\widetilde{p}_{t}}{\widetilde{p}_{t+u}} + \WD{\widetilde{p}_{t+u}}{p_{t+u}} + \WD{p_{t+u}}{p_{t}}
  \end{align*}
  Subtracting $\WD{\widetilde{p}_{t}}{p_{t}}$ from both sides in the first equation,
  and subtracting $\WD{\widetilde{p}_{t+u}}{p_{t+u}}$ in the second equation, we get
  \begin{align*}
    \abs{\WD{\widetilde{p}_{t + u}}{p_{t + u}} - \WD{\widetilde{p}_{t}}{p_{t}}}
    &\leq \WD{\widetilde{p}_{t}}{\widetilde{p}_{t+u}} + \WD{p_{t}}{p_{t+u}}
  \end{align*}
  Furthermore, by continuity of $p_t$,
  \begin{align*}
    \WD{p_{t}}{p_{t+u}} \;\overset{\textrm{\eqref{eq:wd_tovar_bound}}}{\leq}\; \frac{d_{\max}}{2} \cdot \norm{p_{t} - p_{t+u}}_1 \; \to 0 \quad \textrm{ for } \quad u \to 0
  \end{align*}
  The same holds for $\WD{\widetilde{p}_{t}}{\widetilde{p}_{t+u}}$, showing
  that $\WD{\widetilde{p}_{t}}{p_{t}}$ is indeed continuous in $t$.

  For the proof of the main statement, note that
  \begin{align*}
    \WD{\widetilde{p}_{t + u}}{p_{t + u}}
    &\overset{\triangle\textrm{-inequ.}}{\leq} \WD{\widetilde{p}_{t + u}^{\transp}}{\widetilde{p}_t^{\transp} e^{uQ}} + \WD{\widetilde{p}_t^{\transp} e^{uQ}}{p_{t + u}^{\transp}} \\
    &\;\;\;\;=\;\;\;\; \WD{\pi_t^{\transp} e^{u\Theta} A}{\pi_t^{\transp} A e^{uQ}} + \WD{\widetilde{p}_t^{\transp} e^{uQ}}{p_t^{\transp} e^{uQ}}
  \end{align*}
  and apply \autoref{cor:wdderiv_approxdyn_bound} as well as \autoref{cor:wdderiv_error_bound}
  (or \autoref{lem:wdderiv_error_ricciinfbound} when $\underline{k}(Q)$ is replaced by $\underline{\kappa}(Q)$).
  Also note that for $u = 0$, the inequality in the equation above is actually an
  equality, which is why we can bound the one-sided derivative of the left-hand side
  at $u = 0^+$ by the one-sided derivative of the right-hand side.
  The existence of the derivative was shown in \autoref{lem:wdderivative}.
\end{tproof}

Similarly to \cite[Theorem 5]{formalbndsstatespaceredmc} we can deduce that
\begin{itemize}
  \item $\displaystyle \WD{\widetilde{p}_{t}}{p_{t}} \leq \WD{\widetilde{p}_{0}}{p_{0}} + \int_0^t \pi_s^{\transp} \dabs{\Theta A - A Q}_{\mathrm{W}} \dx{s} + t \cdot K(Q)$
  \item $\displaystyle \WD{\widetilde{p}_{t}}{p_{t}} \leq \WD{\widetilde{p}_{0}}{p_{0}} + t \cdot \big(\norm{\Theta A - A Q}_{\mathrm{W}} + K(Q)\big)$
\end{itemize}

\begin{remark}
  If $\DIST$ is the discrete metric, then $\WD{\widetilde{p}_{t}}{p_{t}} = \frac{1}{2}\norm{\widetilde{p}_{t} - p_t}_1$,
  $\norm{\Theta A - A Q}_{\mathrm{W}} = \frac{1}{2} \norm{\Theta A - A Q}_{\infty}$
  and $K(Q) = 0$ (by the proof of \autoref{cor:wdderiv_error_discrmet} and by the
  definition of $K(Q)$ in \autoref{cor:wdderiv_error_bound}), hence
  \begin{align*}
    \frac{1}{2}\norm{\widetilde{p}_{t} - p_t}_1 \leq \frac{1}{2}\norm{\widetilde{p}_{0} - p_0}_1 + t \cdot \frac{1}{2} \norm{\Theta A - A Q}_{\infty}
  \end{align*}
  i.e., we recover the total variation result from \cite[Theorem 5 (iii)]{formalbndsstatespaceredmc}.
  In contrast, \autoref{thm:wasserstein_error_growth_bound} is not applicable
  if we approximate the original process with vectors $\widetilde{p}_t$
  which are no longer probability distributions, or with a matrix $\Theta$ which
  is not necessarily a generator matrix.
\end{remark}

From \autoref{thm:wasserstein_error_growth_bound}, we can also deduce that
\begin{align}
  \WD{\widetilde{p}_{t}}{p_{t}} &\leq \begin{cases}
    \left(\WD{\widetilde{p}_{0}}{p_{0}} - \frac{\norm{\Theta A - A Q}_{\mathrm{W}}}{\underline{k}(Q)}\right) \cdot e^{-\underline{k}(Q) \cdot t} + \frac{\norm{\Theta A - A Q}_{\mathrm{W}}}{\underline{k}(Q)} & \textrm{ if } \underline{k}(Q) \neq 0 \\
    \WD{\widetilde{p}_{0}}{p_{0}} + t \cdot \norm{\Theta A - A Q}_{\mathrm{W}} & \textrm{ otherwise}
  \end{cases}
  \label{eq:werror_growth_integrated}
\end{align}
Again, we could also replace $\underline{k}(Q)$ with $\underline{\kappa}(Q)$ in \eqref{eq:werror_growth_integrated}.
Note the following subtle point: \eqref{eq:werror_growth_integrated} does provide an upper
bound for $\WD{\widetilde{p}_{t}}{p_{t}}$, but the derivative of the bound
in \eqref{eq:werror_growth_integrated} need not be an upper bound for
$\dufrac{t^+} \WD{\widetilde{p}_{t}}{p_{t}}$ when $\underline{k}(Q) > 0$ because
$\WD{\widetilde{p}_{t}}{p_{t}}$ might be strictly smaller than its bound on the right-hand
side of \eqref{eq:werror_growth_integrated}. Indeed, to derive \eqref{eq:werror_growth_integrated}
from \autoref{thm:wasserstein_error_growth_bound}, it is crucial that we have
shown continuity of $\WD{\widetilde{p}_{t}}{p_{t}}$. If we only knew that the
right derivatives of $\WD{\widetilde{p}_{t}}{p_{t}}$ existed at every $t$, $\WD{\widetilde{p}_{t}}{p_{t}}$
could have upwards jumps (while staying right-continuous and still having right derivatives)
crossing the bound from \eqref{eq:werror_growth_integrated}.

\subsubsection{A class of CTMCs with non-negative Ricci curvature}
\label{sssec:class_ctmc_nonneg_curv}

In \eqref{eq:werror_growth_integrated}, we can see that the bound from \autoref{thm:wasserstein_error_growth_bound}
using $\underline{k}(Q)$ or $\underline{\kappa}(Q)$ will grow exponentially if
$\underline{k}(Q) < 0$ (respectively, $\underline{\kappa}(Q) < 0$). We already saw
in \autoref{cor:wdderiv_error_discrmet} that the discrete metric results in $\underline{\kappa}(Q) \geq 0$,
that is, \eqref{eq:werror_growth_integrated} does not grow exponentially. In this section,
we will see another example of a class of CTMCs and a class of metrics with the same property.

We call this class translation-invariant CTMCs, which we define as follows:
consider a CTMC on the (infinite) state space $\bbZ^d$ with the same jump rate (exit rate) $\lambda$
in every state. Furthermore, if $X_{t^-} \in \bbZ^d$ denotes the state of
the CTMC before a jump occurring at time $t$, then $X_t = X_{t^-} + J$ where $J \in \bbZ^d$
has the same distribution for every state $X_{t^-}$. That is,
the jump offsets are identically distributed everywhere in the state space. We truncate
the state space to a finite box $S = ([l_1, u_1] \times \ldots \times [l_d, u_d]) \cap \bbZ^d$
as follows: at every jump, we set $X_t = \rho_S(X_{t^-} + J)$ where $\rho_S$ is the closest-point
projection onto $S$ (according to the usual Euclidean distance). Finally, we assume
that the metric $\DIST$ on $S$ is the usual Euclidean distance, i.e.,
$\dist{r}{s} = \norm{r - s}_2$.

\begin{proposition}
  \label{prop:transl_inv_ctmc_nonneg_ricci}
  Consider a translation-invariant CTMC with jump rate $\lambda$, with jumps
  distributed according to the random variable $J \in \bbZ^d$, and truncated
  to the state space $S \subseteq \bbZ^d$, an axis-aligned box (for details,
  see the paragraph above). Further assume that
  the metric $\DIST$ on $S$ is the usual Euclidean distance.
  Let $Q$ denote the generator of the CTMC.

  Then, we have $\underline{\kappa}(Q) \geq 0$.
\end{proposition}

\begin{tproof}
  Let $r,s \in S$ with $r \neq s$, and let $X^{(r)}_t$ denote the state of the CTMC at time $t$ when started in
  $X^{(r)}_0 = r$. We now define a coupling $\gamma$ between the processes
  $(X^{(r)}_t)_{t \geq 0}$ and $(X^{(s)}_t)_{t \geq 0}$. By assumption, the
  jump times $t_1, t_2, \ldots$ satisfy the following: the inter-arrival times
  $t_1, t_2 - t_1, t_3 - t_2, \ldots$ are iid with distribution $\mathrm{Exp}(\lambda)$.
  Furthermore, the jump offsets $J_1, J_2, \ldots$ (before projection back onto the state space $S$
  in case a jump would leave $S$) are also iid for both processes, with the
  same distribution as $J$. We can therefore define the coupling $\gamma$ such
  that the jump times and offsets agree for the two processes $(X^{(r)}_t)_{t \geq 0}$ and $(X^{(s)}_t)_{t \geq 0}$.

  Let us now consider $\dist{X^{(r)}_t}{X^{(s)}_t}$ under the given coupling. The Wasserstein distance
  stays constant whenever the processes do not jump. When a synchronous jump
  occurs at time $t_i$, then
  \begin{align*}
    X^{(r)}_{t_i} = \rho_S\left(X^{(r)}_{t_i^-} + J_i\right), \qquad
    X^{(s)}_{t_i} = \rho_S\left(X^{(s)}_{t_i^-} + J_i\right)
  \end{align*}
  Now, note that $\rho_S$ maps each point in $\bbZ^d$ to a unique point
  in $S$ because $S$ is an axis-aligned box. Hence, $\rho_S$ simply projects
  each coordinate independently of the others onto the unique closest value
  in the respective coordinate range of the box. In addition,
  $\norm{\rho_S(z_1) - \rho_S(z_2)}_2 = \dist{\rho_S(z_1)}{\rho_S(z_2)} \leq \dist{z_1}{z_2} = \norm{z_1 - z_2}_2$
  for all points $z_1, z_2$,
  which can easily be verified by considering the squared distance and then again
  each coordinate separately. Hence,
  \begin{align*}
    \dist{X^{(r)}_{t_i}}{X^{(s)}_{t_i}}
    &= \dist{\rho_S\left(X^{(r)}_{t_i^-} + J_i\right)}{\rho_S\left(X^{(s)}_{t_i^-} + J_i\right)}
    \leq \dist{X^{(r)}_{t_i^-} + J_i}{X^{(s)}_{t_i^-} + J_i} \\
    &= \norm{X^{(r)}_{t_i^-} + J_i - X^{(s)}_{t_i^-} - J_i}_2
    = \norm{X^{(r)}_{t_i^-} - X^{(s)}_{t_i^-}}_2
    = \dist{X^{(r)}_{t_i^-}}{X^{(s)}_{t_i^-}}
  \end{align*}
  Thus, the distance between $X^{(r)}_t$ and $X^{(s)}_t$ is non-increasing
  at the jump times. As it is constant when no jump occurs, it follows that
  the distance is non-increasing in $t$ in general under the coupling $\gamma$.
  In particular, $\dist{X^{(r)}_t}{X^{(s)}_t} \leq \dist{X^{(r)}_0}{X^{(s)}_0} = \dist{r}{s}$.

  By definition of the Wasserstein distance (\autoref{def:wasserstein_polish}
  and \eqref{eq:wasserstein_def_finite}) and of $X^{(r)}_t, X^{(s)}_t$, we conclude
  \begin{align*}
    \WD{\delta_r^{\transp} e^{tQ}}{\delta_s^{\transp} e^{tQ}}
    &\leq \Exc{\gamma}{\dist{X^{(r)}_t}{X^{(s)}_t}}
    \leq \dist{X^{(r)}_0}{X^{(s)}_0}
    \textrm{ with equalities at } t = 0
  \end{align*}
  It follows that
  \begin{align*}
    \left.\dufrac{t}\right|_{t = 0^+} \WD{\delta_r^{\transp} e^{tQ}}{\delta_s^{\transp} e^{tQ}}
    \leq \left.\dufrac{t}\right|_{t = 0^+} \dist{X^{(r)}_0}{X^{(s)}_0} = 0
  \end{align*}
  and hence, according to \autoref{def:ricci_curvature},
  \begin{align*}
    \kappa(r, s) = -\frac{1}{\dist{r}{s}} \cdot \left.\dufrac{t}\right|_{t = 0^+} \WD{\delta_r^{\transp} e^{tQ}}{\delta_s^{\transp} e^{tQ}} \geq 0
  \end{align*}
  As $r$ and $s$ were arbitrary, we conclude that $\underline{\kappa}(Q) \geq 0$.
\end{tproof}

\begin{remark}
  In \autoref{prop:transl_inv_ctmc_nonneg_ricci}, we don't actually need to
  truncate the CTMC's state space to an axis-aligned box. Instead, it would suffice
  to truncate $\bbZ^d$ to an $S \subseteq \bbZ^d$ such that $\rho_S^{\DIST}(z)$ (the closest-point projection
  onto $S$ according to the metric $\DIST$) is unique for
  every $z \in \bbZ^d$ and, at the same time, $\dist{\rho_S(z_1)}{\rho_S(z_2)} \leq \dist{z_1}{z_2}$
  for all $z_1, z_2 \in \bbZ^d$. $\DIST$ need not be the Euclidean metric,
  but it needs to be translation-invariant, i.e., there should be some $\theta : \bbZ^d \to [0, \infty)$ such
  that $\dist{z_1}{z_2} = \theta(z_1 - z_2)$.
\end{remark}

\subsubsection{An improved linear error bound}
\label{sssec:curv_improvement}

In \autoref{thm:wasserstein_error_growth_bound},
it is possible to improve upon the bound
\begin{align*}
  \dufrac{t^+} \WD{\widetilde{p}_{t}}{p_{t}}
    &\leq \pi_t^{\transp} \dabs{\Theta A - A Q}_{\mathrm{W}}
    + K(Q)
\end{align*}
This bound arises from the inequality
\begin{align}
  \label{eq:linear_k_bound}
  \left.\dufrac{u}\right|_{u = 0^+} \WD{\widetilde{p}_t^{\transp} e^{uQ}}{p_t^{\transp} e^{uQ}}
  &\leq K(Q)
\end{align}
as shown in \autoref{cor:wdderiv_error_bound}, which relies on
\autoref{lem:wdderiv_error_ricciinfbound}. Both \autoref{cor:wdderiv_error_bound}
and \autoref{lem:wdderiv_error_ricciinfbound} apply to the derivative
of the Wasserstein distance between two arbitrary initial distributions.
While we indeed want \eqref{eq:linear_k_bound} to hold for all possible
distributions $p_t$ (we don't want to compute $p_t$ explicitly, so we
cannot make any assumptions about it), we do \emph{not} need
\eqref{eq:linear_k_bound} to hold for any probability distribution
$\widetilde{p}_t$: we actually know $\widetilde{p}_t$ because this is
the approximate transient distribution which we compute.

Now, recall \eqref{eq:ricci_bound_sharp}:
\begin{align*}
  \left.\dufrac{u}\right|_{u = 0^+} \WD{p^{\transp} e^{uQ}}{q^{\transp} e^{uQ}}
  &\leq \sum_{r \neq s} \gamma(r,s) \dist{r}{s} \cdot \big(-\kappa(r, s)\big)
\end{align*}
where $\gamma$ was a coupling achieving the Wasserstein distance between
$p$ and $q$. We can use \eqref{eq:ricci_bound_sharp} with $p = \widetilde{p}_t$ to improve the bound
\eqref{eq:linear_k_bound}: for all probability distributions $q$, we have
\begin{align*}
  \left.\dufrac{u}\right|_{u = 0^+} \WD{\widetilde{p}_t^{\transp} e^{uQ}}{q^{\transp} e^{uQ}}
  &\leq \sum_{r \neq s} \gamma(r,s) \dist{r}{s} \cdot \big(-\kappa(r, s)\big) \\
  &\leq \sum_{r \in S} \underbrace{\left(\sum_{\substack{s \in S\\s \neq r}}\gamma(r,s)\right)}_{\displaystyle \overset{\circledast}{\leq} \widetilde{p}_t(r)} \left(-\min_{\substack{s \in S\\s \neq r}} \dist{r}{s} \kappa(r,s)\right) \\
  &\leq \sum_{r \in S} \widetilde{p}_t(r) \cdot \max\left\{0, -\min_{\substack{s \in S\\s \neq r}} \dist{r}{s} \kappa(r,s)\right\} \\ 
  &\leq \sum_{r \in S} \widetilde{p}_t(r) \cdot \underbrace{\max\left\{0, -\min_{\substack{s \in S\\s \neq r}} \dist{r}{s} k(r,s)\right\}}_{\displaystyle =: K_{\mathrm{loc}}(r)}
\end{align*}
$\circledast$ holds because $\gamma$ is a coupling between $p = \widetilde{p}_t$
and $q$ in \eqref{eq:ricci_bound_sharp}. $K_{\mathrm{loc}}(r)$ is essentially a local version
of $K(Q)$ at the state $r$. Depending on the model, replacing $K(Q)$ with
the $\widetilde{p}_t$-weighted average of the $K_{\mathrm{loc}}(r)$ improved the resulting
error bounds by a factor of around 2-10 if $K(Q)$ was bigger than $\norm{\Theta A - A Q}_{\mathrm{W}}$
in our experiments. However, in these situations, both the bound using $K(Q)$ as
well as the improved version using $K_{\mathrm{loc}}(r)$ were too large to be useful
in practice, which is why we will not mention the improved bound in the
following examples anymore.

\subsubsection{A toy example for illustration}
\label{sssec:toy_example}

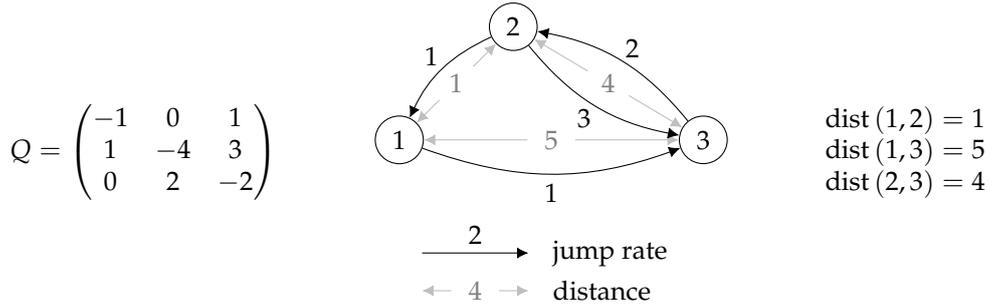
\begin{figure}[htb]
  \begin{center}
    \begin{tabular}{ccccc}
      $\displaystyle Q = \begin{pmatrix}
        -1 & 0 & 1 \\
        1 & -4 & 3 \\
        0 & 2 & -2
      \end{pmatrix}$ & ~~~~~ &
      \begin{tikzpicture}[baseline=(current bounding box.center),>={Latex[length=1.5mm,width=1.5mm]}]
        \node[circle,draw=black] (A) at (0, -1.5) {$1$};
        \node[circle,draw=black] (B) at (1.5, 0) {$2$};
        \node[circle,draw=black] (C) at (4, -1.5) {$3$};
        \draw[gray!50!white,<->] (A) -- node[circle,gray,fill=white] {$1$} (B);
        \draw[gray!50!white,<->] (B) -- node[circle,gray,fill=white] {$4$} (C);
        \draw[gray!50!white,<->] (A) -- node[circle,gray,fill=white] {$5$} (C);
        \draw[->] (A) edge[bend right=20] node[below] {$1$} (C);
        \draw[->] (B) edge[bend right=20] node[above left=-1mm] {$1$} (A);
        \draw[->] (B) edge[bend right=20] node[below left=-1mm] {$3$} (C);
        \draw[->] (C) edge[bend right=20] node[above right=-1mm] {$2$} (B);
        \draw[->] (0.3,-3) -- node[above] {$2$} (1.7,-3) node[right=2mm] {jump rate};
        \draw[gray!50!white,<->] (0.3,-3.5) -- node[circle,gray,fill=white] {$4$} (1.7,-3.5) node[right=2mm,black,fill=white] {distance};
      \end{tikzpicture} & ~~~~~ &
      \begin{minipage}{2.3cm}
        $\dist{1}{2} = 1$ \\
        $\dist{1}{3} = 5$ \\
        $\dist{2}{3} = 4$
      \end{minipage}
    \end{tabular}
  \end{center}
  \caption{A toy CTMC used for illustrating some of the concepts}
  \label{fig:toy_ctmc}
\end{figure}

Let us consider the CTMC on the state space $S = \{1, 2, 3\}$ with generator $Q$ and metric $\DIST$ given in \autoref{fig:toy_ctmc}.

\noindent\textbf{\textit{Ricci curvature.}} First, we calculate $\kappa(1, 2)$ using \autoref{lem:ricci_linearprogram},
which tells us that
\begin{align*}
  &\kappa(1, 2) = -\frac{1}{\dist{1}{2}} \cdot V = -\frac{V}{1} \quad \textrm{ where } V \textrm{ is the solution of } \\
  &\max_{f \in \bbR^n, f \geq 0} (Q_1 - Q_2)f \quad \textrm{s.t.} \quad
  f(1) - f(2) = \dist{1}{2} = 1
  \textrm{ and } f \textrm{ is } 1\textrm{-Lip.}
\end{align*}
We have $Q_1 - Q_2 = \big( -2, \; 4, \; -2 \big)$, so we maximize
$(Q_1 - Q_2)f =  -2 f(1) + 4 f(2) -2 f(3)$. As $f(1) - f(2) = 1$, it follows
that an optimal $f$ is $\big( 5, \; 4, \; 0 \big)^{\transp}$ achieving the objective value $6$, hence $\kappa(1, 2) = -\frac{6}{1} = -6$.
As the distance between states $1$ and $2$ is $\dist{1}{2} = 1$, this implies
(see \autoref{def:ricci_curvature}) that
\begin{align*}
  \left.\dufrac{t}\right|_{t = 0^+} \WD{\delta_1^{\transp} e^{tQ}}{\delta_2^{\transp} e^{tQ}}
  = -\kappa(1,2) \cdot \dist{1}{2} = 6
\end{align*}
That is, the Wasserstein distance between the two transient distributions obtained when starting
in states $1$ and $2$, respectively, initially increases at a rate of $6$.

\begin{figure}[htb]
  \begin{center}
    \begin{tikzpicture}[>={Latex[length=1.5mm,width=1.5mm]}]
      \draw[gray,dashed] (2.50, 0) -- (2.50, 6.70);
      \draw[gray,dashed] (5.00, 0) -- (5.00, 6.70);
      \draw[gray,dashed] (7.50, 0) -- (7.50, 6.70);
      \draw[gray,dashed] (10.00, 0) -- (10.00, 6.70);
      \draw[gray,dashed] (0, 1.50) -- (10.20, 1.50);
      \draw[gray,dashed] (0, 3.00) -- (10.20, 3.00);
      \draw[gray,dashed] (0, 4.50) -- (10.20, 4.50);
      \draw[gray,dashed] (0, 6.00) -- (10.20, 6.00);
      \draw[->] (-0.20, 0) -- (10.30, 0) node[right] {$t$};
      \draw[->] (0, -0.20) -- (0, 6.80) node[above] {Wasserstein distance};
      \draw (2.50, 0.1) -- (2.50, -0.1) node[below] {0.25};
      \draw (5.00, 0.1) -- (5.00, -0.1) node[below] {0.5};
      \draw (7.50, 0.1) -- (7.50, -0.1) node[below] {0.75};
      \draw (10.00, 0.1) -- (10.00, -0.1) node[below] {1};
      \draw (0.1, 1.50) -- (-0.1, 1.50) node[left] {0.5};
      \draw (0.1, 3.00) -- (-0.1, 3.00) node[left] {1};
      \draw (0.1, 4.50) -- (-0.1, 4.50) node[left] {1.5};
      \draw (0.1, 6.00) -- (-0.1, 6.00) node[left] {2};
      \draw[tumOrange,line width=2pt] (0,3) -- (2,6.6) node[below right,fill=white] {line with slope $-\kappa(1, 2) \cdot \dist{1}{2} = 6$};
      \draw[linkred,line width=2pt] (0,3) -- node[above right=4mm,fill=white] {line with slope $-k(1, 2) \cdot \dist{1}{2} = 14$} (0.86,6.6);
      \draw[sectionblue,line width=2pt] (0.00, 3.000) -- (0.02, 3.036) -- (0.04, 3.071) -- (0.07, 3.114) -- (0.09, 3.148) -- (0.11, 3.181) -- (0.13, 3.222) -- (0.15, 3.262) -- (0.18, 3.293) -- (0.20, 3.331) -- (0.22, 3.369) -- (0.24, 3.398) -- (0.27, 3.427) -- (0.29, 3.455) -- (0.31, 3.490) -- (0.33, 3.518) -- (0.35, 3.544) -- (0.38, 3.577) -- (0.40, 3.603) -- (0.42, 3.628) -- (0.44, 3.653) -- (0.46, 3.683) -- (0.48, 3.707) -- (0.50, 3.730) -- (0.52, 3.753) -- (0.54, 3.775) -- (0.56, 3.797) -- (0.58, 3.819) -- (0.60, 3.845) -- (0.62, 3.866) -- (0.65, 3.886) -- (0.67, 3.906) -- (0.69, 3.925) -- (0.71, 3.944) -- (0.73, 3.967) -- (0.75, 3.986) -- (0.77, 4.003) -- (0.79, 4.021) -- (0.81, 4.038) -- (0.83, 4.054) -- (0.85, 4.071) -- (0.87, 4.087) -- (0.90, 4.106) -- (0.92, 4.121) -- (0.94, 4.136) -- (0.96, 4.150) -- (0.98, 4.165) -- (1.00, 4.178) -- (1.02, 4.192) -- (1.04, 4.205) -- (1.06, 4.221) -- (1.08, 4.233) -- (1.10, 4.246) -- (1.12, 4.257) -- (1.14, 4.269) -- (1.16, 4.280) -- (1.18, 4.291) -- (1.21, 4.304) -- (1.23, 4.315) -- (1.25, 4.325) -- (1.27, 4.335) -- (1.29, 4.344) -- (1.31, 4.353) -- (1.33, 4.362) -- (1.35, 4.371) -- (1.37, 4.379) -- (1.39, 4.387) -- (1.41, 4.395) -- (1.43, 4.403) -- (1.45, 4.412) -- (1.47, 4.419) -- (1.49, 4.426) -- (1.51, 4.433) -- (1.53, 4.439) -- (1.55, 4.445) -- (1.57, 4.451) -- (1.59, 4.457) -- (1.61, 4.463) -- (1.63, 4.468) -- (1.65, 4.473) -- (1.67, 4.478) -- (1.69, 4.482) -- (1.71, 4.487) -- (1.73, 4.491) -- (1.75, 4.495) -- (1.77, 4.500) -- (1.79, 4.504) -- (1.81, 4.507) -- (1.83, 4.510) -- (1.85, 4.513) -- (1.88, 4.516) -- (1.90, 4.519) -- (1.92, 4.521) -- (1.94, 4.523) -- (1.96, 4.525) -- (1.98, 4.527) -- (2.00, 4.529) -- (2.02, 4.531) -- (2.04, 4.532) -- (2.06, 4.534) -- (2.08, 4.535) -- (2.10, 4.536) -- (2.12, 4.537) -- (2.14, 4.537) -- (2.16, 4.538) -- (2.26, 4.538) -- (2.28, 4.537) -- (2.30, 4.537) -- (2.32, 4.536) -- (2.34, 4.535) -- (2.36, 4.534) -- (2.38, 4.533) -- (2.40, 4.531) -- (2.42, 4.530) -- (2.44, 4.529) -- (2.46, 4.527) -- (2.48, 4.525) -- (2.50, 4.523) -- (2.52, 4.521) -- (2.54, 4.519) -- (2.56, 4.517) -- (2.58, 4.515) -- (2.60, 4.512) -- (2.62, 4.510) -- (2.65, 4.507) -- (2.67, 4.504) -- (2.69, 4.501) -- (2.71, 4.498) -- (2.73, 4.495) -- (2.75, 4.492) -- (2.77, 4.489) -- (2.79, 4.485) -- (2.81, 4.482) -- (2.83, 4.478) -- (2.85, 4.475) -- (2.87, 4.471) -- (2.89, 4.466) -- (2.91, 4.462) -- (2.93, 4.458) -- (2.95, 4.454) -- (2.97, 4.450) -- (2.99, 4.446) -- (3.01, 4.441) -- (3.03, 4.437) -- (3.05, 4.432) -- (3.07, 4.428) -- (3.09, 4.423) -- (3.11, 4.418) -- (3.13, 4.413) -- (3.15, 4.409) -- (3.17, 4.404) -- (3.19, 4.398) -- (3.21, 4.393) -- (3.23, 4.388) -- (3.25, 4.383) -- (3.27, 4.378) -- (3.29, 4.372) -- (3.31, 4.367) -- (3.33, 4.361) -- (3.35, 4.356) -- (3.37, 4.350) -- (3.39, 4.344) -- (3.41, 4.339) -- (3.43, 4.333) -- (3.45, 4.327) -- (3.47, 4.321) -- (3.49, 4.315) -- (3.51, 4.309) -- (3.53, 4.301) -- (3.55, 4.295) -- (3.57, 4.289) -- (3.59, 4.283) -- (3.61, 4.276) -- (3.63, 4.270) -- (3.65, 4.264) -- (3.67, 4.257) -- (3.69, 4.251) -- (3.71, 4.244) -- (3.73, 4.238) -- (3.75, 4.231) -- (3.77, 4.224) -- (3.79, 4.218) -- (3.81, 4.211) -- (3.83, 4.204) -- (3.85, 4.197) -- (3.88, 4.190) -- (3.90, 4.184) -- (3.92, 4.177) -- (3.94, 4.170) -- (3.96, 4.163) -- (3.98, 4.156) -- (4.00, 4.148) -- (4.02, 4.141) -- (4.04, 4.132) -- (4.06, 4.125) -- (4.08, 4.116) -- (4.11, 4.109) -- (4.13, 4.100) -- (4.16, 4.091) -- (4.18, 4.082) -- (4.20, 4.074) -- (4.22, 4.065) -- (4.25, 4.057) -- (4.27, 4.048) -- (4.29, 4.041) -- (4.31, 4.031) -- (4.33, 4.024) -- (4.36, 4.016) -- (4.38, 4.006) -- (4.40, 3.999) -- (4.42, 3.989) -- (4.45, 3.982) -- (4.47, 3.972) -- (4.49, 3.964) -- (4.52, 3.955) -- (4.54, 3.945) -- (4.56, 3.937) -- (4.58, 3.927) -- (4.61, 3.920) -- (4.63, 3.910) -- (4.65, 3.902) -- (4.68, 3.892) -- (4.70, 3.882) -- (4.72, 3.874) -- (4.75, 3.864) -- (4.77, 3.856) -- (4.79, 3.846) -- (4.81, 3.837) -- (4.83, 3.829) -- (4.86, 3.819) -- (4.88, 3.811) -- (4.90, 3.803) -- (4.92, 3.793) -- (4.95, 3.784) -- (4.97, 3.774) -- (4.99, 3.766) -- (5.02, 3.756) -- (5.04, 3.746) -- (5.06, 3.738) -- (5.09, 3.728) -- (5.11, 3.718) -- (5.13, 3.708) -- (5.16, 3.698) -- (5.18, 3.688) -- (5.21, 3.679) -- (5.23, 3.669) -- (5.25, 3.659) -- (5.28, 3.649) -- (5.30, 3.641) -- (5.32, 3.631) -- (5.34, 3.623) -- (5.37, 3.614) -- (5.39, 3.604) -- (5.42, 3.594) -- (5.44, 3.584) -- (5.46, 3.576) -- (5.48, 3.566) -- (5.50, 3.557) -- (5.53, 3.549) -- (5.55, 3.539) -- (5.58, 3.529) -- (5.60, 3.519) -- (5.62, 3.511) -- (5.64, 3.500) -- (5.67, 3.492) -- (5.69, 3.484) -- (5.71, 3.474) -- (5.74, 3.464) -- (5.76, 3.454) -- (5.79, 3.443) -- (5.81, 3.433) -- (5.83, 3.425) -- (5.86, 3.415) -- (5.88, 3.405) -- (5.91, 3.395) -- (5.93, 3.387) -- (5.95, 3.377) -- (5.97, 3.369) -- (5.99, 3.360) -- (6.02, 3.350) -- (6.04, 3.340) -- (6.07, 3.330) -- (6.09, 3.322) -- (6.11, 3.312) -- (6.13, 3.304) -- (6.15, 3.296) -- (6.18, 3.286) -- (6.20, 3.276) -- (6.23, 3.266) -- (6.25, 3.258) -- (6.27, 3.248) -- (6.29, 3.240) -- (6.31, 3.232) -- (6.34, 3.222) -- (6.36, 3.212) -- (6.39, 3.202) -- (6.41, 3.192) -- (6.43, 3.183) -- (6.46, 3.175) -- (6.48, 3.165) -- (6.50, 3.155) -- (6.53, 3.145) -- (6.55, 3.137) -- (6.57, 3.127) -- (6.59, 3.120) -- (6.62, 3.112) -- (6.64, 3.102) -- (6.67, 3.092) -- (6.69, 3.082) -- (6.71, 3.075) -- (6.73, 3.065) -- (6.75, 3.057) -- (6.78, 3.049) -- (6.80, 3.040) -- (6.83, 3.030) -- (6.85, 3.021) -- (6.87, 3.013) -- (6.89, 3.003) -- (6.92, 2.996) -- (6.94, 2.988) -- (6.96, 2.978) -- (6.99, 2.969) -- (7.01, 2.959) -- (7.03, 2.952) -- (7.05, 2.942) -- (7.08, 2.935) -- (7.10, 2.927) -- (7.12, 2.918) -- (7.14, 2.910) -- (7.16, 2.903) -- (7.19, 2.893) -- (7.21, 2.884) -- (7.24, 2.874) -- (7.26, 2.867) -- (7.28, 2.858) -- (7.30, 2.850) -- (7.32, 2.843) -- (7.35, 2.834) -- (7.37, 2.824) -- (7.40, 2.815) -- (7.42, 2.808) -- (7.44, 2.799) -- (7.46, 2.791) -- (7.48, 2.784) -- (7.51, 2.775) -- (7.53, 2.766) -- (7.56, 2.757) -- (7.58, 2.749) -- (7.60, 2.740) -- (7.62, 2.733) -- (7.64, 2.726) -- (7.67, 2.717) -- (7.69, 2.708) -- (7.71, 2.699) -- (7.74, 2.690) -- (7.76, 2.681) -- (7.79, 2.674) -- (7.81, 2.665) -- (7.83, 2.656) -- (7.86, 2.647) -- (7.88, 2.640) -- (7.90, 2.632) -- (7.92, 2.625) -- (7.95, 2.618) -- (7.97, 2.609) -- (8.00, 2.600) -- (8.02, 2.593) -- (8.04, 2.584) -- (8.06, 2.576) -- (8.09, 2.569) -- (8.11, 2.560) -- (8.13, 2.552) -- (8.16, 2.543) -- (8.19, 2.535) -- (8.21, 2.526) -- (8.23, 2.518) -- (8.26, 2.511) -- (8.28, 2.502) -- (8.30, 2.494) -- (8.33, 2.486) -- (8.36, 2.477) -- (8.38, 2.469) -- (8.40, 2.462) -- (8.43, 2.454) -- (8.45, 2.446) -- (8.47, 2.437) -- (8.50, 2.431) -- (8.52, 2.422) -- (8.54, 2.414) -- (8.57, 2.408) -- (8.59, 2.399) -- (8.62, 2.391) -- (8.64, 2.383) -- (8.67, 2.375) -- (8.69, 2.367) -- (8.71, 2.361) -- (8.74, 2.353) -- (8.76, 2.345) -- (8.79, 2.337) -- (8.81, 2.329) -- (8.84, 2.321) -- (8.86, 2.313) -- (8.88, 2.306) -- (8.91, 2.298) -- (8.93, 2.291) -- (8.96, 2.283) -- (8.98, 2.275) -- (9.00, 2.267) -- (9.03, 2.261) -- (9.05, 2.253) -- (9.07, 2.245) -- (9.10, 2.238) -- (9.12, 2.232) -- (9.14, 2.224) -- (9.17, 2.216) -- (9.19, 2.210) -- (9.21, 2.203) -- (9.24, 2.195) -- (9.27, 2.188) -- (9.29, 2.180) -- (9.31, 2.172) -- (9.34, 2.166) -- (9.36, 2.159) -- (9.38, 2.152) -- (9.41, 2.144) -- (9.44, 2.137) -- (9.46, 2.129) -- (9.48, 2.122) -- (9.51, 2.116) -- (9.53, 2.109) -- (9.55, 2.102) -- (9.58, 2.094) -- (9.61, 2.087) -- (9.63, 2.080) -- (9.65, 2.074) -- (9.68, 2.067) -- (9.70, 2.060) -- (9.72, 2.053) -- (9.75, 2.047) -- (9.77, 2.040) -- (9.79, 2.033) -- (9.82, 2.027) -- (9.84, 2.020) -- (9.87, 2.013) -- (9.89, 2.006) -- (9.92, 1.999) -- (9.94, 1.992) -- (9.96, 1.986) -- (9.99, 1.980) -- (10.00, 1.975) node[below left,fill=white] {$\WD{\delta_1^{\transp} e^{tQ}}{\delta_2^{\transp} e^{tQ}}$};
    \end{tikzpicture}
  \end{center}
  \caption{Evolution of the Wasserstein distance between the two transient distributions obtained when starting
  in states $1$ and $2$ of the toy CTMC}
  \label{fig:ricci_example_curvalongstates}
\end{figure}
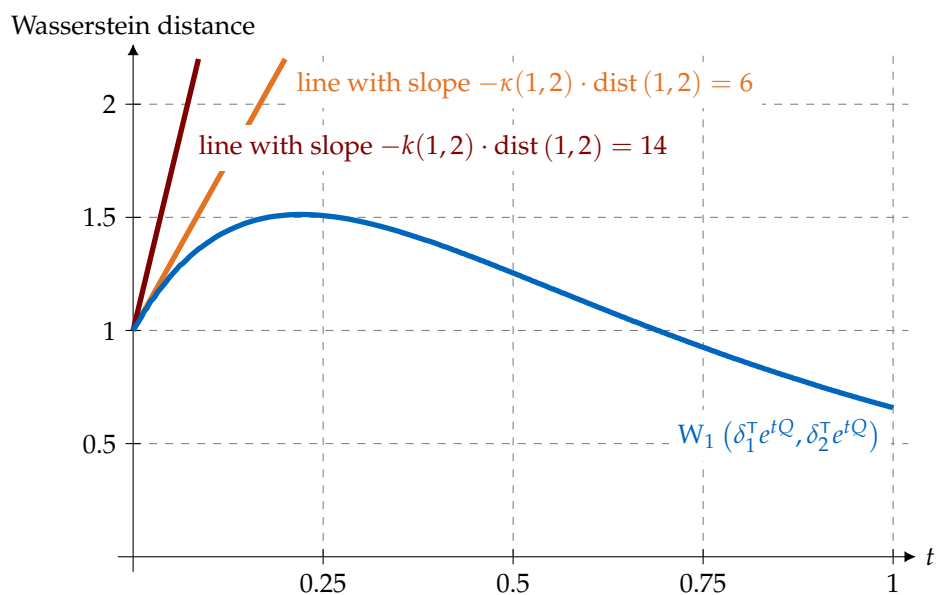

The lower bound $k(1,2)$ from \autoref{lem:ricci_lower_bound} yields:
\begin{align*}
  k(1,2) &= -\frac{\min\left\{Q_{1} \dist{1}{\cdot}, \;\; Q_{1} \dist{2}{\cdot}\right\} + \min\left\{Q_{2} \dist{2}{\cdot}, \;\; Q_{2} \dist{1}{\cdot}\right\}}{\dist{1}{2}} \\
  &= -\frac{\min\left\{5, \;\; 3\right\} + \min\left\{13, \;\; 11\right\}}{1} = -\frac{14}{1} = -14
\end{align*}
\autoref{fig:ricci_example_curvalongstates} shows how the Wasserstein distance
between $\delta_1^{\transp} e^{tQ}$ and $\delta_2^{\transp} e^{tQ}$ evolves.
The initial slope of this curve matches $-\kappa(1, 2) \cdot \dist{1}{2}$. Since $\kappa(1, 2) \geq k(1, 2)$,
$-k(1, 2) \cdot \dist{1}{2}$ is an upper bound on the initial derivative. This example
shows that the distance between transient distributions can grow (if the underlying metric
is not the discrete metric).

\begin{figure}[htb]
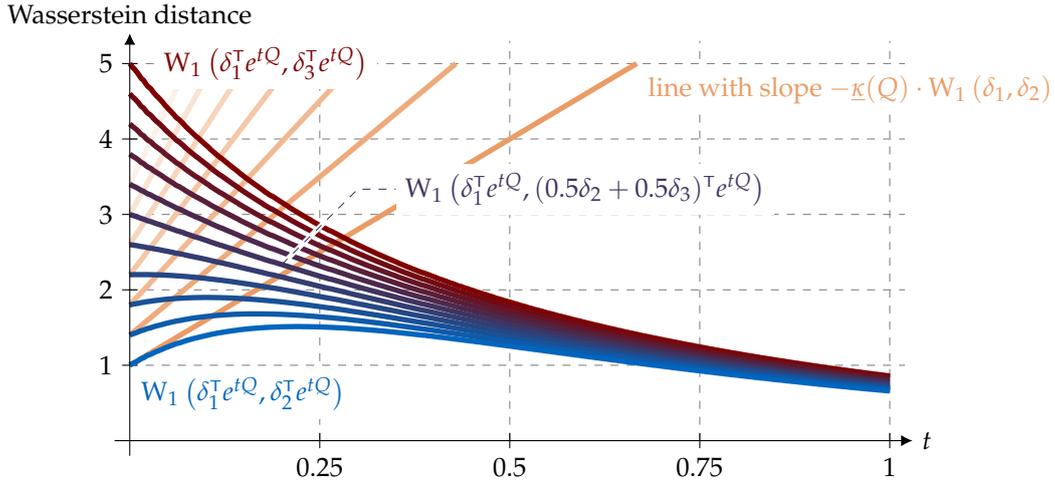

  \begin{center}

  \end{center}
  \caption{Evolution of the Wasserstein distance between two transient distributions of the toy CTMC:
  on the one hand the transient distribution obtained when starting in state $1$;
  on the other hand the initial distribution $\alpha \delta_2 + (1 - \alpha) \delta_3$ for $\alpha \in [0, 1]$
  is used. The orange lines show the bound on the initial derivative of the
  Wasserstein distance given by \autoref{lem:wdderiv_error_ricciinfbound}.}
  \label{fig:ricci_example_mixedcurvature}
\end{figure}

We can calculate the coarse Ricci curvature and the lower bound
for all pairs of states; the result is shown in \autoref{tab:ricci_example} (on the left).
As we can see, the Ricci curvature is positive (and exactly matches the lower bound)
for the two other state pairs. 
\begin{table}[htb]
  \begin{center}
    \begin{minipage}{6cm}
      \centering
      with metric $\DIST$ from \autoref{fig:toy_ctmc}:\\
      \def\arraystretch{1.2}
      \begin{tabular}{|c|c|c|c|}
        \hline
        $r$ & $s$ & $\kappa(r, s)$ & $k(r, s)$ \\ \hline
        $1$ & $2$ & $-6$ & $-14$ \\ \hline
        $1$ & $3$ & $2.6$ & $2.6$ \\ \hline
        $2$ & $3$ & $4.75$ & $4.75$ \\ \hline
      \end{tabular}
    \end{minipage}
    ~~
    \begin{minipage}{6cm}
      \centering
      with $\DIST$ being the discrete metric:\\
      \def\arraystretch{1.2}
      \begin{tabular}{|c|c|c|c|}
        \hline
        $r$ & $s$ & $\kappa(r, s)$ & $k(r, s)$ \\ \hline
        $1$ & $2$ & $2$ & $1$ \\ \hline
        $1$ & $3$ & $1$ & $1$ \\ \hline
        $2$ & $3$ & $5$ & $5$ \\ \hline
      \end{tabular}
    \end{minipage}
  \end{center}
  \caption{Ricci curvature and lower bounds $k(r,s)$ for the toy CTMC}
  \label{tab:ricci_example}
\end{table}
Using \autoref{tab:ricci_example} (left side), we can also calculate $\underline{\kappa}(Q) = -6$, $\underline{k}(Q) = -14$ and $K(Q) = 14$.
By \autoref{lem:wdderiv_error_ricciinfbound}, the initial derivative of 
$\WD{p^{\transp}e^{tQ}}{q^{\transp}e^{tQ}}$ for any two initial distributions $p \in \bbR^3$ and $q \in \bbR^3$ is upper bounded by
$-\underline{\kappa}(Q) \cdot \WD{p}{q}$. This bound is attained
when $p = \delta_1$ and $q = \delta_2$, i.e., when we choose Dirac distributions
on the two states $r, s$ for which $\underline{\kappa}(Q) = \kappa(r, s)$. However,
for other initial distributions, the initial derivative can be much lower.
\autoref{tab:ricci_example} (left side) gives us an indication for which initial
distributions this might be the case: the more the initial distributions resemble
Dirac measures on $r$ and $s$, the closer the initial derivative will be to
$-\kappa(r, s) \cdot \dist{r}{s}$. For example, if we choose $p = \delta_1$ and $q = \delta_3$
as initial distributions, then we have
\begin{align*}
  \left.\dufrac{t}\right|_{t = 0^+} \WD{p^{\transp}e^{tQ}}{q^{\transp}e^{tQ}}
  = -\kappa(1, 3) \cdot \underbrace{\WD{p}{q}}_{\dist{1}{3}} = -2.6 \cdot 5 \ll 6 \cdot 5 = -\underline{\kappa}(Q) \cdot \WD{p}{q}
\end{align*}
\autoref{fig:ricci_example_mixedcurvature} visualizes $\WD{p^{\transp}e^{tQ}}{q^{\transp}e^{tQ}}$
for $p = \delta_1$ and $q = \alpha \delta_2 + (1 - \alpha) \delta_3$ with $\alpha \in [0, 1]$.
The orange lines show the upper bound on $\dufrac{t}\big|_{t = 0^+} \WD{p^{\transp}e^{tQ}}{q^{\transp}e^{tQ}}$
given by \autoref{lem:wdderiv_error_ricciinfbound}, i.e., these lines have slope
$-\underline{\kappa}(Q) \cdot \WD{p}{q}$. When $p = \delta_1$ and $q = \delta_2$,
the bound matches the initial slope exactly, but the bound then gradually becomes worse
as $q$ has less resemblance with $\delta_2$ and approaches $\delta_3$.

Single state pairs like $1$ and $2$ in the toy CTMC can thus cause $\underline{\kappa}(Q)$ to
become negative (which is bad in the sense that this entails a bound which predicts an
increasing Wasserstein distance), even though the other state pairs are better behaved
(have positive coarse Ricci curvature $\kappa(r, s)$). This can be problematic for bounding
the growth of the accumulated error in the aggregation setting: the bound on the accumulated
error will grow exponentially if $\underline{\kappa}(Q) < 0$ and $\underline{k}(Q) < 0$, even though
the accumulated error will actually decrease in the long run as the transient distributions
approach stationarity. In such a scenario, it makes sense to use the bound with $K(Q)$ from \autoref{thm:wasserstein_error_growth_bound} instead,
as this will result in an accumulated error bound which grows at most linearly.

If the metric from \autoref{fig:toy_ctmc} is replaced by the discrete metric,
then by the proof of \autoref{cor:wdderiv_error_discrmet}, $\underline{\kappa}(Q) \geq \underline{k}(Q) \geq 0$. Indeed,
as we can see in \autoref{tab:ricci_example} on the right,
we even have $\underline{\kappa}(Q) = \underline{k}(Q) = 1 > 0$ in this example. Hence,
the bound from \autoref{cor:wdderiv_error_bound} implies that, regardless of the initial
distributions, the Wasserstein distance between two transient distributions is always
decreasing (unless it is already equal to zero). As $K(Q) = 0$ in this case, the second
bound from \autoref{cor:wdderiv_error_bound} would not reflect that (we actually always
have $K(Q) \geq 0$, so using $K(Q)$ will never result in a bound showing that the
Wasserstein distance will initially decrease). In contrast to what
we have seen for the metric from \autoref{fig:toy_ctmc}, the bound
on the accumulated error would now decrease exponentially if we use
the bound with $\underline{k}(Q)$ from \autoref{thm:wasserstein_error_growth_bound}.

\begin{figure}[htb]
  \begin{center}
    \begin{tikzpicture}[>={Latex[length=1.5mm,width=1.5mm]}]
      \draw[gray,dashed] (2.50, 0) -- (2.50, 5.30);
      \draw[gray,dashed] (5.00, 0) -- (5.00, 5.30);
      \draw[gray,dashed] (7.50, 0) -- (7.50, 5.30);
      \draw[gray,dashed] (10.00, 0) -- (10.00, 5.30);
      \draw[gray,dashed] (0, 1.00) -- (10.20, 1.00);
      \draw[gray,dashed] (0, 2.00) -- (10.20, 2.00);
      \draw[gray,dashed] (0, 3.00) -- (10.20, 3.00);
      \draw[gray,dashed] (0, 4.00) -- (10.20, 4.00);
      \draw[gray,dashed] (0, 5.00) -- (10.20, 5.00);
      \draw (2.50, 0.1) -- (2.50, -0.1) node[below] {0.25};
      \draw (5.00, 0.1) -- (5.00, -0.1) node[below] {0.5};
      \draw (7.50, 0.1) -- (7.50, -0.1) node[below] {0.75};
      \draw (10.00, 0.1) -- (10.00, -0.1) node[below] {1};
      \draw (0.1, 1.00) -- (-0.1, 1.00) node[left] {1};
      \draw (0.1, 2.00) -- (-0.1, 2.00) node[left] {2};
      \draw (0.1, 3.00) -- (-0.1, 3.00) node[left] {3};
      \draw (0.1, 4.00) -- (-0.1, 4.00) node[left] {4};
      \draw (0.1, 5.00) -- (-0.1, 5.00) node[left] {5};
      \draw[tumOrange!70!white,dashed,line width=2pt,domain=0:5.8,samples=30,smooth] plot (\x,{(exp(6*\x/10)-1)/6}) node[below right,fill=white] {bound using $\underline{\kappa}(Q)$};
      \draw[linkred,line width=2pt] (0, 0) -- node[above left=0.8cm and -0.4cm,align=center,fill=white] {bound\\using\\$K(Q)$} (3.47, 5.2);
      \draw[tumOrange,line width=2pt,domain=0:3.07,samples=30,smooth] plot (\x,{(exp(14*\x/10)-1)/14}) node[above right=0mm and -6mm,fill=white] {bound using $\underline{k}(Q)$};
      \draw[sectionblue,line width=2pt] (0.00, 0) -- (0.10, 0.010) -- (0.20, 0.021) -- (0.30, 0.031) -- (0.40, 0.041) -- (0.51, 0.052) -- (0.61, 0.062) -- (0.71, 0.073) -- (0.81, 0.083) -- (0.91, 0.093) -- (1.01, 0.103) -- (1.11, 0.113) -- (1.21, 0.123) -- (1.31, 0.132) -- (1.41, 0.141) -- (1.52, 0.150) -- (1.62, 0.159) -- (1.72, 0.168) -- (1.82, 0.176) -- (1.92, 0.184) -- (2.02, 0.191) -- (2.12, 0.199) -- (2.22, 0.206) -- (2.32, 0.213) -- (2.42, 0.219) -- (2.53, 0.225) -- (2.63, 0.231) -- (2.73, 0.237) -- (2.83, 0.243) -- (2.93, 0.248) -- (3.03, 0.253) -- (3.13, 0.257) -- (3.23, 0.262) -- (3.33, 0.266) -- (3.43, 0.269) -- (3.54, 0.273) -- (3.64, 0.276) -- (3.74, 0.279) -- (3.84, 0.282) -- (3.94, 0.285) -- (4.04, 0.287) -- (4.14, 0.290) -- (4.24, 0.292) -- (4.34, 0.293) -- (4.44, 0.295) -- (4.55, 0.296) -- (4.65, 0.297) -- (4.75, 0.299) -- (4.85, 0.299) -- (4.95, 0.300) -- (5.05, 0.301) -- (5.66, 0.301) -- (5.76, 0.300) -- (5.86, 0.300) -- (5.96, 0.299) -- (6.06, 0.298) -- (6.16, 0.297) -- (6.26, 0.296) -- (6.36, 0.295) -- (6.46, 0.294) -- (6.57, 0.293) -- (6.67, 0.292) -- (6.77, 0.290) -- (6.87, 0.289) -- (6.97, 0.287) -- (7.07, 0.286) -- (7.17, 0.284) -- (7.27, 0.282) -- (7.37, 0.280) -- (7.47, 0.278) -- (7.58, 0.277) -- (7.68, 0.275) -- (7.78, 0.273) -- (7.88, 0.271) -- (7.98, 0.269) -- (8.08, 0.266) -- (8.18, 0.264) -- (8.28, 0.262) -- (8.38, 0.260) -- (8.48, 0.258) -- (8.59, 0.256) -- (8.69, 0.253) -- (8.79, 0.251) -- (8.89, 0.249) -- (8.99, 0.246) -- (9.09, 0.244) -- (9.19, 0.242) -- (9.29, 0.239) -- (9.39, 0.237) -- (9.49, 0.235) -- (9.60, 0.232) -- (9.70, 0.230) -- (9.80, 0.228) -- (9.90, 0.225) -- (10.00, 0.223) node[above left=3mm,fill=white] {actual error $\WD{\widetilde{p}_t}{p_t}$};
      \draw[->] (-0.20, 0) -- (10.30, 0) node[right] {$t$};
      \draw[->] (0, -0.20) -- (0, 5.40) node[above] {Wasserstein distance};
    \end{tikzpicture}
  \end{center}
  \caption{Error evolution for the toy CTMC using the given aggregation and
  the initial distribution $p_0 = (\frac{1}{2}, \; \frac{1}{2}, \; 0)^{\transp}$. The red
  and orange lines show the error bounds obtained from \autoref{thm:wasserstein_error_growth_bound}.}
  \label{fig:toy_aggregation_error}
\end{figure}
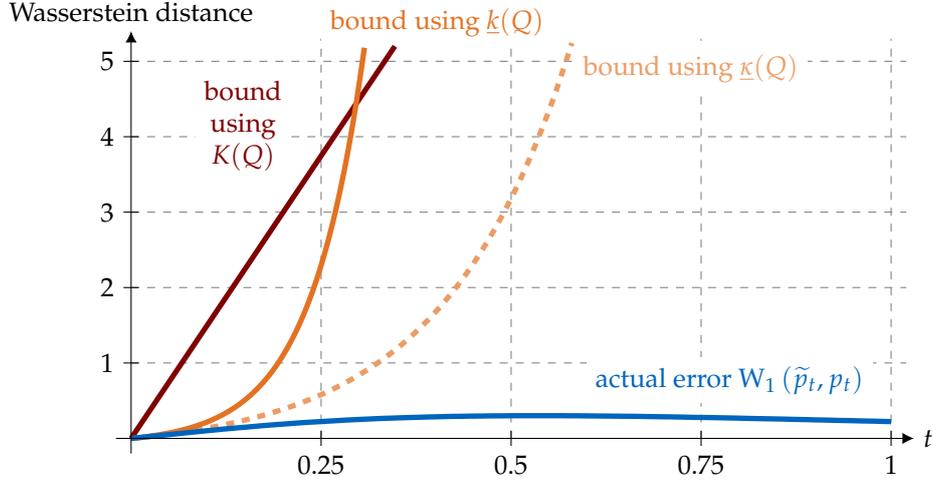

\noindent\textbf{\textit{Aggregation.}} In this example, an aggregation to approximate transient distributions doesn't make much sense
due to the low dimension of the state space, and a more realistic scenario is shown in the next section. Still,
for illustration, it is interesting to consider an aggregation of the toy CTMC.
For the following, we consider again
the metric defined in \autoref{fig:toy_ctmc} for the toy example.
We simply aggregate the two states which are closest
according to $\DIST$, that is, states $1$ and $2$, and both are assigned
equal weight within the aggregate. Hence, we put (for the definition
of $\Lambda$, see the end of \autoref{ssec:prelim_mc})
\begin{align*}
  A = \begin{pmatrix}
    \frac{1}{2} & \frac{1}{2} & 0 \\
    0 & 0 & 1
  \end{pmatrix}, \;
  \Lambda = \begin{pmatrix}
    1 & 0 \\
    1 & 0 \\
    0 & 1
  \end{pmatrix}, \;
  \Theta = A Q \Lambda = \begin{pmatrix}
    -2 & 2 \\
    2 & -2
  \end{pmatrix}, \;
  \Theta A - A Q = \begin{pmatrix}
    -1 & 1 & 0 \\
    1 & -1 & 0
  \end{pmatrix}
\end{align*}
Here, $\dabs{\Theta A - A Q}_{\mathrm{W}}$ is easy to calculate. We can express
\begin{align*}
  \Theta A - A Q &= \begin{pmatrix}
    0 & 1 & 0 \\
    1 & 0 & 0
  \end{pmatrix} - \begin{pmatrix}
    1 & 0 & 0 \\
    0 & 1 & 0
  \end{pmatrix} = \begin{pmatrix}
    \hpipe \, \delta_2^{\transp} \, \hpipe \\
    \hpipe \, \delta_1^{\transp} \, \hpipe
  \end{pmatrix} - \begin{pmatrix}
    \hpipe \, \delta_1^{\transp} \, \hpipe \\
    \hpipe \, \delta_2^{\transp} \, \hpipe
  \end{pmatrix} \\
  \implies \dabs{\Theta A - A Q}_{\mathrm{W}} &= \begin{pmatrix}
    \WD{\delta_2}{\delta_1} \\
    \WD{\delta_1}{\delta_2}
  \end{pmatrix} = \begin{pmatrix}
    \dist{2}{1} \\
    \dist{1}{2}
  \end{pmatrix} = \begin{pmatrix}
    1 \\
    1
  \end{pmatrix}
\end{align*}
(Compare with the remark after \autoref{def:row_wd_norm}.) Therefore,
by \autoref{thm:wasserstein_error_growth_bound}, we have
\begin{align*}
  \dufrac{t^+} \WD{\widetilde{p}_{t}}{p_{t}}
  &\leq \pi_t^{\transp} \dabs{\Theta A - A Q}_{\mathrm{W}}
  + K(Q)
  = \pi_t^{\transp} \begin{pmatrix} 1 \\ 1 \end{pmatrix} + 14 = 15 \\
  \textrm{and } \;\; \dufrac{t^+} \WD{\widetilde{p}_{t}}{p_{t}}
  &\leq \pi_t^{\transp} \dabs{\Theta A - A Q}_{\mathrm{W}}
  + \WD{\widetilde{p}_{t}}{p_{t}} \cdot \big(-\underline{k}(Q)\big)
  = 1 + 14 \cdot \WD{\widetilde{p}_{t}}{p_{t}}
\end{align*}
Integrating, we get the bounds
\begin{align}
  \WD{\widetilde{p}_{t}}{p_{t}} \leq \WD{\widetilde{p}_{0}}{p_{0}} + 15t
  \quad &\textrm{ and } \quad
  \WD{\widetilde{p}_{t}}{p_{t}} \leq \left(\WD{\widetilde{p}_{0}}{p_{0}} + \frac{1}{14}\right) e^{14t} - \frac{1}{14}
  \label{eq:example_error_bounds}
\end{align}
A slight improvement would be possible by using the second bound as long as its
derivative is smaller than $15$, and switch to the first bound otherwise.
As an example, we consider $p_0 = (\frac{1}{2}, \; \frac{1}{2}, \; 0)^{\transp}$ so that
$\pi_0 = (1, \; 0)^{\transp}$ and $\widetilde{p}_0^{\transp} = \pi_0^{\transp} A = p_0^{\transp}$ (there is thus
no initial error). We can actually calculate $\pi_t$ explicitly in this case:
\begin{align*}
  \pi_t^{\transp} &= \pi_0^{\transp} e^{t \Theta} = \left(\frac{1}{2}\big(1 + e^{-4t}\big), \; \frac{1}{2}\big(1 - e^{-4t}\big)\right) \\
  \implies \widetilde{p}_t^{\transp} &= \pi_t^{\transp}A = \left(\frac{1}{4}\big(1 + e^{-4t}\big), \; \frac{1}{4}\big(1 + e^{-4t}\big), \; \frac{1}{2}\big(1 - e^{-4t}\big)\right)
\end{align*}
The analytical expression for the actual $p_t$ is already quite complicated, so we omit
it here. \autoref{fig:toy_aggregation_error} compares the actual error
$\WD{\widetilde{p}_t}{p_t}$ to the error bounds from \eqref{eq:example_error_bounds}
(and a third error bound obtained when using $\underline{\kappa}(Q)$ instead
of $\underline{k}(Q)$). We can see that the bounds using $\underline{k}(Q)$
and $\underline{\kappa}(Q)$ are close to the actual error for $t$ near $0$, but then,
the bounds grow exponentially while the actual error plateaus and even decreases near $t = 1$.
The bound using $K(Q)$ does not grow exponentially, but is already far off near $t = 0$.

\begin{figure}[htb]
  \begin{center}
    \begin{tikzpicture}[>={Latex[length=1.5mm,width=1.5mm]}]
      \draw[gray,dashed] (2.50, 0) -- (2.50, 5.30);
      \draw[gray,dashed] (5.00, 0) -- (5.00, 5.30);
      \draw[gray,dashed] (7.50, 0) -- (7.50, 5.30);
      \draw[gray,dashed] (10.00, 0) -- (10.00, 5.30);
      \draw[gray,dashed] (0, 1.00) -- (10.20, 1.00);
      \draw[gray,dashed] (0, 2.00) -- (10.20, 2.00);
      \draw[gray,dashed] (0, 3.00) -- (10.20, 3.00);
      \draw[gray,dashed] (0, 4.00) -- (10.20, 4.00);
      \draw[gray,dashed] (0, 5.00) -- (10.20, 5.00);
      \draw (2.50, 0.1) -- (2.50, -0.1) node[below] {0.25};
      \draw (5.00, 0.1) -- (5.00, -0.1) node[below] {0.5};
      \draw (7.50, 0.1) -- (7.50, -0.1) node[below] {0.75};
      \draw (10.00, 0.1) -- (10.00, -0.1) node[below] {1};
      \draw (0.1, 1.00) -- (-0.1, 1.00) node[left] {0.2};
      \draw (0.1, 2.00) -- (-0.1, 2.00) node[left] {0.4};
      \draw (0.1, 3.00) -- (-0.1, 3.00) node[left] {0.6};
      \draw (0.1, 4.00) -- (-0.1, 4.00) node[left] {0.8};
      \draw (0.1, 5.00) -- (-0.1, 5.00) node[left] {1};
      \draw[linkred,line width=2pt] (0, 0) -- node[above left,fill=white] {bound using $K(Q)$} (10, 5);
      \draw[tumOrange,line width=2pt,domain=0:10,samples=50,smooth] plot (\x,{(1 - exp(-\x/10))*5}) node[right,align=center,fill=white] {bound\\using\\$\underline{k}(Q)$};
      \draw[sectionblue,line width=2pt] (0.00, 0) -- (0.10, 0.049) -- (0.20, 0.095) -- (0.30, 0.138) -- (0.40, 0.179) -- (0.51, 0.217) -- (0.61, 0.253) -- (0.71, 0.287) -- (0.81, 0.318) -- (0.91, 0.347) -- (1.01, 0.375) -- (1.11, 0.400) -- (1.21, 0.424) -- (1.31, 0.447) -- (1.41, 0.467) -- (1.52, 0.486) -- (1.62, 0.504) -- (1.72, 0.520) -- (1.82, 0.535) -- (1.92, 0.549) -- (2.02, 0.562) -- (2.12, 0.574) -- (2.22, 0.584) -- (2.32, 0.594) -- (2.42, 0.602) -- (2.53, 0.610) -- (2.63, 0.617) -- (2.73, 0.623) -- (2.83, 0.629) -- (2.93, 0.634) -- (3.03, 0.638) -- (3.13, 0.641) -- (3.23, 0.644) -- (3.33, 0.646) -- (3.43, 0.648) -- (3.54, 0.649) -- (3.64, 0.650) -- (3.84, 0.650) -- (3.94, 0.649) -- (4.04, 0.648) -- (4.14, 0.647) -- (4.24, 0.646) -- (4.34, 0.644) -- (4.44, 0.642) -- (4.55, 0.639) -- (4.65, 0.636) -- (4.75, 0.633) -- (4.85, 0.630) -- (4.95, 0.627) -- (5.05, 0.623) -- (5.15, 0.620) -- (5.25, 0.616) -- (5.35, 0.612) -- (5.45, 0.607) -- (5.56, 0.603) -- (5.66, 0.598) -- (5.76, 0.594) -- (5.86, 0.589) -- (5.96, 0.584) -- (6.06, 0.579) -- (6.16, 0.574) -- (6.26, 0.569) -- (6.36, 0.564) -- (6.46, 0.559) -- (6.57, 0.554) -- (6.67, 0.549) -- (6.77, 0.543) -- (6.87, 0.538) -- (6.97, 0.533) -- (7.07, 0.527) -- (7.17, 0.522) -- (7.27, 0.517) -- (7.37, 0.511) -- (7.47, 0.506) -- (7.58, 0.500) -- (7.68, 0.495) -- (7.78, 0.489) -- (7.88, 0.484) -- (7.98, 0.479) -- (8.08, 0.473) -- (8.18, 0.468) -- (8.28, 0.462) -- (8.38, 0.457) -- (8.48, 0.452) -- (8.59, 0.446) -- (8.69, 0.441) -- (8.79, 0.436) -- (8.89, 0.431) -- (8.99, 0.425) -- (9.09, 0.420) -- (9.19, 0.415) -- (9.29, 0.410) -- (9.39, 0.405) -- (9.49, 0.400) -- (9.60, 0.395) -- (9.70, 0.390) -- (9.80, 0.385) -- (9.90, 0.380) -- (10.00, 0.376) node[above left=3mm and 2mm,fill=white] {actual error $\WD{\widetilde{p}_t}{p_t}$};
      \draw[->] (-0.20, 0) -- (10.30, 0) node[right] {$t$};
      \draw[->] (0, -0.20) -- (0, 5.40) node[above] {Wasserstein distance};
    \end{tikzpicture}
  \end{center}
  \caption{Error evolution for the toy CTMC with the discrete metric using the given aggregation and
  the initial distribution $p_0 = (\frac{1}{2}, \; \frac{1}{2}, \; 0)^{\transp}$. The red
  and orange lines show the error bounds obtained from \autoref{thm:wasserstein_error_growth_bound}.}
  \label{fig:toy_aggregation_error_discrmet}
\end{figure}
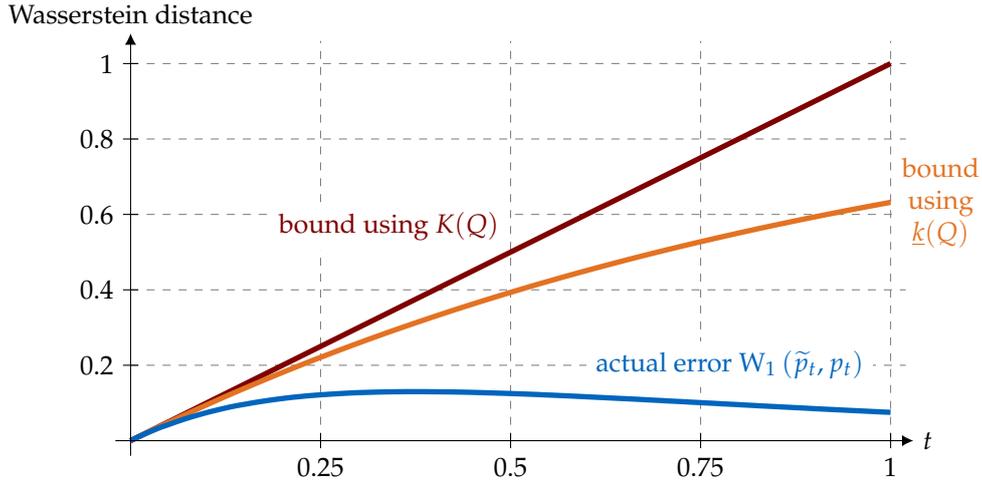

To conclude the section on the toy example, we consider again what happens if we use the example
with the discrete metric. This case is basically already covered in
\cite{formalbndsstatespaceredmc}, but we actually get a small improvement
in the error bounds for this particular example because $\underline{k}(Q) > 0$,
which tells us that the accumulated error will actually decrease over time.
More precisely, by \autoref{thm:wasserstein_error_growth_bound}, we
get (if we use the discrete metric; note that $\dabs{\Theta A - A Q}_{\mathrm{W}}$
remains unchanged here, which is not true in general if we use another metric)
\begin{align*}
  \dufrac{t^+} \WD{\widetilde{p}_{t}}{p_{t}}
  &\leq \pi_t^{\transp} \dabs{\Theta A - A Q}_{\mathrm{W}}
  + K(Q)
  = \pi_t^{\transp} \begin{pmatrix} 1 \\ 1 \end{pmatrix} + 0 = 1 \\
  \textrm{and } \;\; \dufrac{t^+} \WD{\widetilde{p}_{t}}{p_{t}}
  &\leq \pi_t^{\transp} \dabs{\Theta A - A Q}_{\mathrm{W}}
  + \WD{\widetilde{p}_{t}}{p_{t}} \cdot \big(-\underline{k}(Q)\big)
  = 1 - \WD{\widetilde{p}_{t}}{p_{t}}
\end{align*}
Integrating, we get the bounds (cf. \eqref{eq:werror_growth_integrated})
\begin{align*}
  \WD{\widetilde{p}_{t}}{p_{t}} \leq \WD{\widetilde{p}_{0}}{p_{0}} + t
  \quad \textrm{ and } \quad
  \WD{\widetilde{p}_{t}}{p_{t}} \leq \left(\WD{\widetilde{p}_{0}}{p_{0}} - 1\right) e^{-t} + 1
\end{align*}
The first bound is the one we would also have obtained from the technique in \cite{formalbndsstatespaceredmc}.
Using again $p_0 = (\frac{1}{2}, \; \frac{1}{2}, \; 0)^{\transp}$, we get
the picture shown in \autoref{fig:toy_aggregation_error_discrmet}. Both bounds
are close to the actual error near $t = 0$, but the bounds quickly become worse
as $t$ grows. Due to the positive Ricci curvature, the bound using $\underline{k}(Q)$
does not explode exponentially in this case, but grows more slowly with increasing $t$ instead.
This is because the accumulated error decreases over time due to the positive Ricci curvature~--
the growth in the error bound is caused solely by the bound from \autoref{cor:wdderiv_approxdyn_bound}.

\subsubsection{A more realistic example: the RSVP model}
\label{sssec:rsvp}

We next consider an example with a larger CTMC:
a compositional stochastic process algebra model,
the RSVP model from \cite{rsvp}. It comprises a lower network channel submodel
with capacity for $M$ calls, an upper network channel submodel with capacity for
$N$ calls, and a number of identical mobile nodes which request resources
for calls at a constant rate. Due to the mobile node symmetry in the model specification,
a lossless state space reduction is possible for this model. We use
$M = 7$, $N = 5$ and $3$ mobile nodes, resulting in a total of $842$ states.
If the considered initial distribution is compatible with the lossless
aggregation comprising $234$ aggregates (for details, see \cite{formalbndsstatespaceredmc}),
the aggregation scheme will calculate exact transient distributions,
i.e., $\widetilde{p}_t = p_t$.

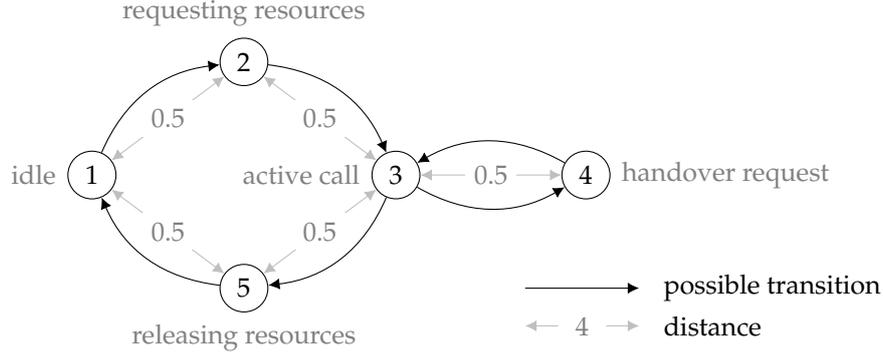
\begin{figure}[htb]
  \begin{center}
    \begin{tikzpicture}[>={Latex[length=1.5mm,width=1.5mm]}]
      \node[circle,draw=black] (A) at (-1,0) {$1$};
      \node[circle,draw=black] (B) at (1,1.5) {$2$};
      \node[circle,draw=black] (C) at (3,0) {$3$};
      \node[circle,draw=black] (D) at (5.5,0) {$4$};
      \node[circle,draw=black] (E) at (1,-1.5) {$5$};
      \node[gray,left=3.5mm] at (A) {idle};
      \node[gray,left=3.5mm] at (C) {active call};
      \node[gray,above=3.5mm] at (B) {requesting resources};
      \node[gray,below=3.5mm] at (E) {releasing resources};
      \node[gray,right=3.5mm] at (D) {handover request};
      \draw[gray!50!white,<->] (A) -- node[gray,fill=white] {$0.5$} (B);
      \draw[gray!50!white,<->] (B) -- node[gray,fill=white] {$0.5$} (C);
      \draw[gray!50!white,<->] (C) -- node[gray,fill=white] {$0.5$} (E);
      \draw[gray!50!white,<->] (E) -- node[gray,fill=white] {$0.5$} (A);
      \draw[gray!50!white,<->] (C) -- node[gray,fill=white] {$0.5$} (D);
      \draw[->] (A) edge[bend left=30] (B);
      \draw[->] (B) edge[bend left=30] (C);
      \draw[->] (C) edge[bend left=30] (E);
      \draw[->] (E) edge[bend left=30] (A);
      \draw[->] (C) edge[bend right=30] (D);
      \draw[->] (D) edge[bend right=30] (C);
      \draw[->] (4.7,-1.5) -- (6.2,-1.5) node[right=2mm] {possible transition};
      \draw[gray!50!white,<->] (4.7,-2) -- node[circle,gray,fill=white] {$4$} (6.2,-2) node[right=2mm,black,fill=white] {distance};
    \end{tikzpicture}
  \end{center}
  \caption{Suggested metric for measuring the distance between two mobile node
  states: the shortest path metric $d_{G}$ according to the gray distances in the graph given above}
  \label{fig:rsvp_mobile_metric}
\end{figure}

\textbf{\textit{The metric.}} The RSVP model was not defined in conjunction with a metric on its state space
originally. However, if we take a closer look at the model specification,
we can suggest a sensible choice of metric. The state of the CTMC consists of
$6$ components, $3$ for the states of the $3$ mobile nodes, $1$ component
for the lower network channel, $1$ component for the upper network channel,
and $1$ component for the channel monitor, which is responsible for handling
handover requests arising when a mobile node switches between network cells. That is, a state $s \in \bbN^6$
of the RSVP model looks as follows:
\begin{align*}
  s &= \big(s(1), \; s(2), \; s(3), \; s(4), \; s(5), \; s(6)\big)^{\transp} \in \bbN^6 \quad \textrm{ with } \\
  s(1), s(2), s(3) &\in \{1, \ldots, 5\} \textrm{ internal state of the mobile nodes} \\
  s(4) &\in \{0, \ldots, M\} \textrm{ number of resources currently used in the lower network channel} \\
  s(5) &\in \{0, \ldots, N\} \textrm{ number of resources currently used in the upper network channel} \\
  s(6) &\in \{0, \ldots, M\} \textrm{ number of handover requests handled by the channel monitor}
\end{align*}
Not all states are reachable, which is why the model with $M = 7$, $N = 5$ and $3$ mobile nodes
only contains $842$ states. For $s(4),s(5),s(6)$, it is natural to simply take
the absolute value of the component difference of two states $s$ and $\widetilde{s}$ of the
CTMC as a measure of how far apart the component states are. For the mobile nodes,
the five internal states are as follows:
\begin{align*}
  &1: \; \textrm{idle}, \;\;
  2: \; \textrm{requesting network resources}, \;\;
  3: \; \textrm{active call}, \;\; \\
  &4: \; \textrm{handover request}, \;\;
  5: \; \textrm{releasing resources}
\end{align*}
For each mobile node component, we suggest to set the distance between two states
as the shortest path distance $d_{G}$ in the graph given in \autoref{fig:rsvp_mobile_metric}.
Overall, we then suggest the following metric on the state space $S$ of the CTMC
arising from the RSVP model:
\begin{align*}
  \dist{r}{s} &:= \overbrace{d_G\big(r(1), s(1)\big) + d_G\big(r(2), s(2)\big) + d_G\big(r(3), s(3)\big)}^{\textrm{mobile nodes}} \\
  &\hphantom{\;:=\;} + \underbrace{\abs{r(4) - s(4)}}_{\textrm{lower n.c.}} + \underbrace{\abs{r(5) - s(5)}}_{\textrm{upper n.c.}} + \frac{1}{2}\underbrace{\abs{r(6) - s(6)}}_{\textrm{channel m.}} \qquad \textrm{for } r,s \in S \subseteq \bbR^6
\end{align*}
The state space then has a diameter of $d_{\max} = \max_{r, s} \dist{r}{s} = 18$.

\textbf{\textit{Aggregation and erros.}} We aggregated the CTMC using \cite[Algorithm 3]{formalbndsstatespaceredmc} with $\varepsilon = 0.1$,
resulting in $123$ aggregates. This aggregation is not exact, we have $\Theta A \neq AQ$.
The computation of $K(Q)$ and $\underline{k}(Q)$ takes less than $10$ seconds on our test machine
(single-threaded execution on an Intel Core i7-1260P CPU with a maximum frequency of 4.7 GHz)
and results in $K(Q) \approx 130$ and $\underline{k}(Q) \approx -254$. Comparing these
to the diameter $d_{\max} = 18$, we see that the error bounds will not be very useful in practice,
growing beyond the diameter of the state space (i.e., the maximal possible error)
very quickly. It does not help that $\norm{\Theta A - A Q}_{\mathrm{W}} \approx 0.165$
is of a more reasonable size. 

\begin{figure}[htb]
  \begin{center}
    \begin{tikzpicture}[>={Latex[length=1.5mm,width=1.5mm]}]
      \draw[gray!20!white,line width=3pt,decorate,decoration={snake,amplitude=0.4mm,segment length=3mm}] (-0.5,9.66) -- (10.6,9.66);
      \draw[gray,dashed] (2.50, 9.93) -- (2.50, 13.84);
      \draw[gray,dashed] (5.00, 9.93) -- (5.00, 13.84);
      \draw[gray,dashed] (7.50, 9.93) -- (7.50, 13.84);
      \draw[gray,dashed] (10.00, 9.93) -- (10.00, 13.84);
      \draw[gray,dashed] (2.50, 9.09) -- (2.50, 9.39);
      \draw[gray,dashed] (5.00, 9.09) -- (5.00, 9.39);
      \draw[gray,dashed] (7.50, 9.09) -- (7.50, 9.39);
      \draw[gray,dashed] (10.00, 9.09) -- (10.00, 9.39);
      \draw[gray,dashed] (0, 10.23) -- (10.20, 10.23);
      \draw[gray,dashed] (0, 11.36) -- (10.20, 11.36);
      \draw[gray,dashed] (0, 12.50) -- (10.20, 12.50);
      \draw[gray,dashed] (0, 13.64) -- (10.20, 13.64);
      \draw[->] (-0.20, 9.09) -- (10.30, 9.09) node[right] {$t$};
      \draw (0, 8.89) -- (0, 9.39);
      \draw[->] (0, 9.93) -- (0, 13.94) node[above] {$\WD{\widetilde{p}_t}{p_t}$};
      \draw (2.50, 9.19) -- (2.50, 8.99) node[below] {5};
      \draw (5.00, 9.19) -- (5.00, 8.99) node[below] {10};
      \draw (7.50, 9.19) -- (7.50, 8.99) node[below] {15};
      \draw (10.00, 9.19) -- (10.00, 8.99) node[below] {20};
      \draw (0.1, 10.23) -- (-0.1, 10.23) node[left] {2.25};
      \draw (0.1, 11.36) -- (-0.1, 11.36) node[left] {2.5};
      \draw (0.1, 12.50) -- (-0.1, 12.50) node[left] {2.75};
      \draw (0.1, 13.64) -- (-0.1, 13.64) node[left] {3};
      \draw[sectionblue,line width=2pt] (0.00, 13.636) -- (0.20, 13.636) -- (0.30, 13.631) -- (0.40, 13.619) -- (0.51, 13.601) -- (0.61, 13.576) -- (0.71, 13.544) -- (0.81, 13.504) -- (0.91, 13.457) -- (1.01, 13.403) -- (1.11, 13.340) -- (1.21, 13.271) -- (1.31, 13.194) -- (1.41, 13.110) -- (1.52, 13.023) -- (1.62, 12.935) -- (1.72, 12.841) -- (1.82, 12.751) -- (1.92, 12.659) -- (2.02, 12.583) -- (2.12, 12.510) -- (2.22, 12.446) -- (2.32, 12.394) -- (2.42, 12.346) -- (2.53, 12.310) -- (2.63, 12.276) -- (2.73, 12.244) -- (2.83, 12.214) -- (2.93, 12.184) -- (3.03, 12.156) -- (3.13, 12.128) -- (3.23, 12.101) -- (3.33, 12.075) -- (3.43, 12.049) -- (3.54, 12.023) -- (3.64, 11.998) -- (3.74, 11.973) -- (3.84, 11.948) -- (3.94, 11.923) -- (4.04, 11.899) -- (4.14, 11.874) -- (4.24, 11.850) -- (4.34, 11.825) -- (4.44, 11.801) -- (4.55, 11.777) -- (4.65, 11.753) -- (4.75, 11.729) -- (4.85, 11.705) -- (4.95, 11.681) -- (5.05, 11.657) -- (5.15, 11.633) -- (5.25, 11.609) -- (5.35, 11.585) -- (5.45, 11.561) -- (5.56, 11.537) -- (5.66, 11.513) -- (5.76, 11.490) -- (5.86, 11.466) -- (5.96, 11.442) -- (6.06, 11.419) -- (6.16, 11.395) -- (6.26, 11.371) -- (6.36, 11.348) -- (6.46, 11.324) -- (6.57, 11.301) -- (6.67, 11.277) -- (6.77, 11.254) -- (6.87, 11.230) -- (6.97, 11.207) -- (7.07, 11.183) -- (7.17, 11.160) -- (7.27, 11.137) -- (7.37, 11.113) -- (7.47, 11.090) -- (7.58, 11.067) -- (7.68, 11.044) -- (7.78, 11.021) -- (7.88, 10.998) -- (7.98, 10.974) -- (8.08, 10.951) -- (8.18, 10.928) -- (8.28, 10.905) -- (8.38, 10.882) -- (8.48, 10.860) -- (8.59, 10.837) -- (8.69, 10.814) -- (8.79, 10.791) -- (8.89, 10.768) -- (8.99, 10.745) -- (9.09, 10.723) -- (9.19, 10.700) -- (9.29, 10.677) -- (9.39, 10.655) -- (9.49, 10.632) -- (9.60, 10.610) -- (9.70, 10.590) -- (9.80, 10.583) -- (9.90, 10.575) -- (10.00, 10.574);
    \end{tikzpicture}
  \end{center}
  \caption{Evolution of the actual error $\WD{\widetilde{p}_t}{p_t}$ for the CTMC
  arising from the RSVP model with $M = 7$, $N = 5$ and $3$ mobile nodes (resulting in $842$ states),
  aggregated using \cite[Algorithm 3]{formalbndsstatespaceredmc} with $\varepsilon = 0.1$
  (resulting in $123$ aggregates). The initial distribution $p_0$ was chosen to be
  the Dirac measure on the initial state of the RSVP model (no active calls,
  all mobile nodes idle). As this initial state belongs to an aggregate with
  $4$ additional states, an error already occurs in the approximation of the
  initial distribution $p_0$ with $\widetilde{p}_0 = A^{\transp}\pi_0$.}
  \label{fig:rsvp_aggregation_error}
\end{figure}
The actual error $\WD{\widetilde{p}_t}{p_t}$ is plotted
in \autoref{fig:rsvp_aggregation_error} (when $p_0$ is the Dirac measure
on the initial state of the RSVP model).
We use the following, slightly modified versions of the error bounds in \autoref{thm:wasserstein_error_growth_bound}
(which are easier to integrate in order to obtain an error bound at a specific time point):
\begin{align}
  \begin{split}
    \dufrac{t^+} \WD{\widetilde{p}_{t}}{p_{t}}
    &\leq \norm{\Theta A - A Q}_{\mathrm{W}}
    + K(Q) \\
    \textrm{and } \;\; \dufrac{t^+} \WD{\widetilde{p}_{t}}{p_{t}}
    &\leq \norm{\Theta A - A Q}_{\mathrm{W}}
    + \WD{\widetilde{p}_{t}}{p_{t}} \cdot \big(-\underline{k}(Q)\big)
  \end{split}
  \label{eq:err_bnds_modified}
\end{align}
With these bounds, the first bound (using $K(Q)$) would already hit the state space
diameter around time $t \approx 0.12$, and the secound bound (using $\underline{k}(Q)$) would
hit the diameter around time $t \approx 0.007$. That is, the bounds are not
very useful (and not shown in \autoref{fig:rsvp_aggregation_error} because they would
grow out of the pictured range almost immediately).

\begin{figure}[htb]
  \begin{center}
    \begin{tikzpicture}[>={Latex[length=1.5mm,width=1.5mm]}]
      \draw[gray,dashed] (2.50, 0) -- (2.50, 4.75);
      \draw[gray,dashed] (5.00, 0) -- (5.00, 4.75);
      \draw[gray,dashed] (7.50, 0) -- (7.50, 4.75);
      \draw[gray,dashed] (10.00, 0) -- (10.00, 4.75);
      \draw[gray,dashed] (0, 0.73) -- (10.20, 0.73);
      \draw[gray,dashed] (0, 1.46) -- (10.20, 1.46);
      \draw[gray,dashed] (0, 2.18) -- (10.20, 2.18);
      \draw[gray,dashed] (0, 2.91) -- (10.20, 2.91);
      \draw[gray,dashed] (0, 3.64) -- (10.20, 3.64);
      \draw[gray,dashed] (0, 4.37) -- (10.20, 4.37);
      \draw[->] (-0.20, 0) -- (10.30, 0) node[right] {$t$};
      \draw[->] (0, -0.20) -- (0, 4.85) node[above] {$\WD{\widetilde{p}_t}{p_t}$};
      \draw (2.50, 0.1) -- (2.50, -0.1) node[below] {5};
      \draw (5.00, 0.1) -- (5.00, -0.1) node[below] {10};
      \draw (7.50, 0.1) -- (7.50, -0.1) node[below] {15};
      \draw (10.00, 0.1) -- (10.00, -0.1) node[below] {20};
      \draw (0.1, 0.73) -- (-0.1, 0.73) node[left] {0.25};
      \draw (0.1, 1.46) -- (-0.1, 1.46) node[left] {0.5};
      \draw (0.1, 2.18) -- (-0.1, 2.18) node[left] {0.75};
      \draw (0.1, 2.91) -- (-0.1, 2.91) node[left] {1};
      \draw (0.1, 3.64) -- (-0.1, 3.64) node[left] {1.25};
      \draw (0.1, 4.37) -- (-0.1, 4.37) node[left] {1.5};
      \draw[sectionblue,line width=2pt] (0.00, 0) -- (0.10, 0.011) -- (0.20, 0.031) -- (0.30, 0.058) -- (0.40, 0.089) -- (0.51, 0.126) -- (0.61, 0.167) -- (0.71, 0.211) -- (0.81, 0.257) -- (0.91, 0.307) -- (1.01, 0.359) -- (1.11, 0.414) -- (1.21, 0.470) -- (1.31, 0.527) -- (1.41, 0.587) -- (1.52, 0.647) -- (1.62, 0.708) -- (1.72, 0.769) -- (1.82, 0.831) -- (1.92, 0.893) -- (2.02, 0.956) -- (2.12, 1.018) -- (2.22, 1.080) -- (2.32, 1.142) -- (2.42, 1.204) -- (2.53, 1.266) -- (2.63, 1.327) -- (2.73, 1.387) -- (2.83, 1.447) -- (2.93, 1.507) -- (3.03, 1.566) -- (3.13, 1.625) -- (3.23, 1.683) -- (3.33, 1.741) -- (3.43, 1.799) -- (3.54, 1.855) -- (3.64, 1.912) -- (3.74, 1.967) -- (3.84, 2.023) -- (3.94, 2.077) -- (4.04, 2.131) -- (4.14, 2.185) -- (4.24, 2.238) -- (4.34, 2.291) -- (4.44, 2.343) -- (4.55, 2.395) -- (4.65, 2.446) -- (4.75, 2.497) -- (4.85, 2.547) -- (4.95, 2.597) -- (5.05, 2.646) -- (5.15, 2.695) -- (5.25, 2.744) -- (5.35, 2.792) -- (5.45, 2.839) -- (5.56, 2.886) -- (5.66, 2.933) -- (5.76, 2.979) -- (5.86, 3.024) -- (5.96, 3.069) -- (6.06, 3.114) -- (6.16, 3.159) -- (6.26, 3.202) -- (6.36, 3.246) -- (6.46, 3.289) -- (6.57, 3.331) -- (6.67, 3.374) -- (6.77, 3.415) -- (6.87, 3.457) -- (6.97, 3.498) -- (7.07, 3.538) -- (7.17, 3.578) -- (7.27, 3.618) -- (7.37, 3.657) -- (7.47, 3.696) -- (7.58, 3.734) -- (7.68, 3.772) -- (7.78, 3.810) -- (7.88, 3.847) -- (7.98, 3.884) -- (8.08, 3.921) -- (8.18, 3.957) -- (8.28, 3.993) -- (8.38, 4.028) -- (8.48, 4.063) -- (8.59, 4.098) -- (8.69, 4.132) -- (8.79, 4.166) -- (8.89, 4.200) -- (8.99, 4.233) -- (9.09, 4.266) -- (9.19, 4.298) -- (9.29, 4.330) -- (9.39, 4.362) -- (9.49, 4.393) -- (9.60, 4.424) -- (9.70, 4.455) -- (9.80, 4.486) -- (9.90, 4.516) -- (10.00, 4.545);
    \end{tikzpicture}
  \end{center}
  \caption{Evolution of the actual error $\WD{\widetilde{p}_t}{p_t}$ for the CTMC
  arising from the RSVP model with $M = 7$, $N = 5$ and $3$ mobile nodes (resulting in $842$ states),
  aggregated using \cite[Algorithm 3]{formalbndsstatespaceredmc} with $\varepsilon = 0.1$
  (resulting in $123$ aggregates). The initial distribution $p_0$ was chosen to be a uniform
  distribution over the aggregate containing the initial state of the RSVP model (no active calls,
  all mobile nodes idle).}
  \label{fig:rsvp_aggregation_error_zeroiniterror}
\end{figure}
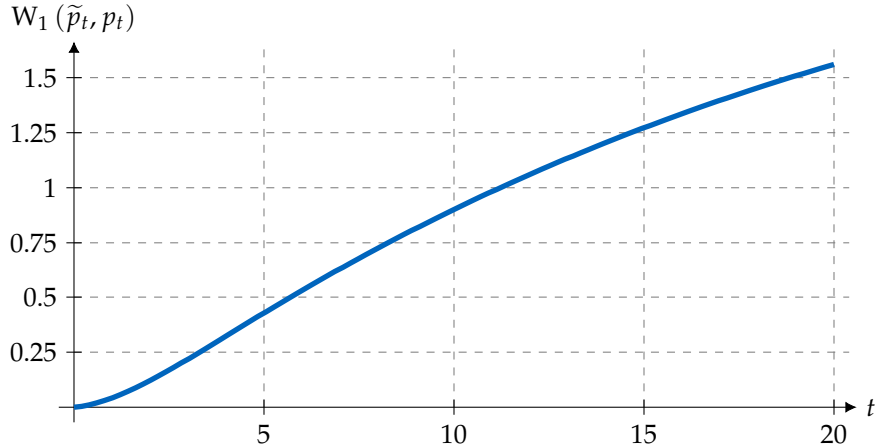

If we start with an initial distribution $p_0$ which does not cause any
error in the approximation of the initial distribution, then the actual
error behaves as in \autoref{fig:rsvp_aggregation_error_zeroiniterror}.
Again, the modified error bounds from \eqref{eq:err_bnds_modified}
are not useful, growing larger than the state space diameter near times
$t \approx 0.14$ (bound using $K(Q)$) and $t \approx 0.04$ (bound using $\underline{k}(Q)$).

\textbf{\textit{Ricci curvature.}} Calculating $\underline{\kappa}(Q)$ for the CTMC arising from the RSVP model with
$842$ states is computationally significantly more expensive than calculating
$\underline{k}(Q)$ or $K(Q)$ as a linear program has to be solved for every state pair
(with $842$ variables, $842^2$ inequality constraints and $1$ equality constraint,
see \autoref{lem:ricci_linearprogram}). Indeed, when using SciPy for solving these
linear programs, calculating $\kappa(r,s)$ for all state pairs would take too long on
our test machine. Instead, we can look at state pairs for which $k(r,s)$ is low. Looking
at one of the pairs for which $k(r,s)$ is minimal, we find that $\underline{\kappa}(Q) \lesssim -53.995$.
We can then calculate $\kappa(r,s)$ only for those pairs for which $k(r,s) < -53.99$,
which is sufficient to find $\underline{\kappa}(Q)$ by \autoref{lem:ricci_lower_bound}.
Using this strategy, we only have to calculate $\kappa(r,s)$ for around $0.7\%$ of the state pairs,
and we get that $\underline{\kappa}(Q) \approx -53.995$ (so the pair for which
$k(r,s)$ was minimal and which we chose was indeed a minimizer of $\kappa(r,s)$ as well).
The calculation still takes around 2.5 hours on our test machine (single-threaded).

Even when using $\underline{\kappa}(Q)$ instead of $\underline{k}(Q)$
in the modified error bounds from \eqref{eq:err_bnds_modified}, the error bound
still grows larger than the state space diameter near time $t \approx 0.16$
when the initial distribution $p_0$ is the uniform distribution over the aggregate
containing the initial state of the RSVP model. While this is an improvement
over the bound using $\underline{k}(Q)$, it is not sufficient to yield a
practically useful bound.

\begin{figure}[htb]
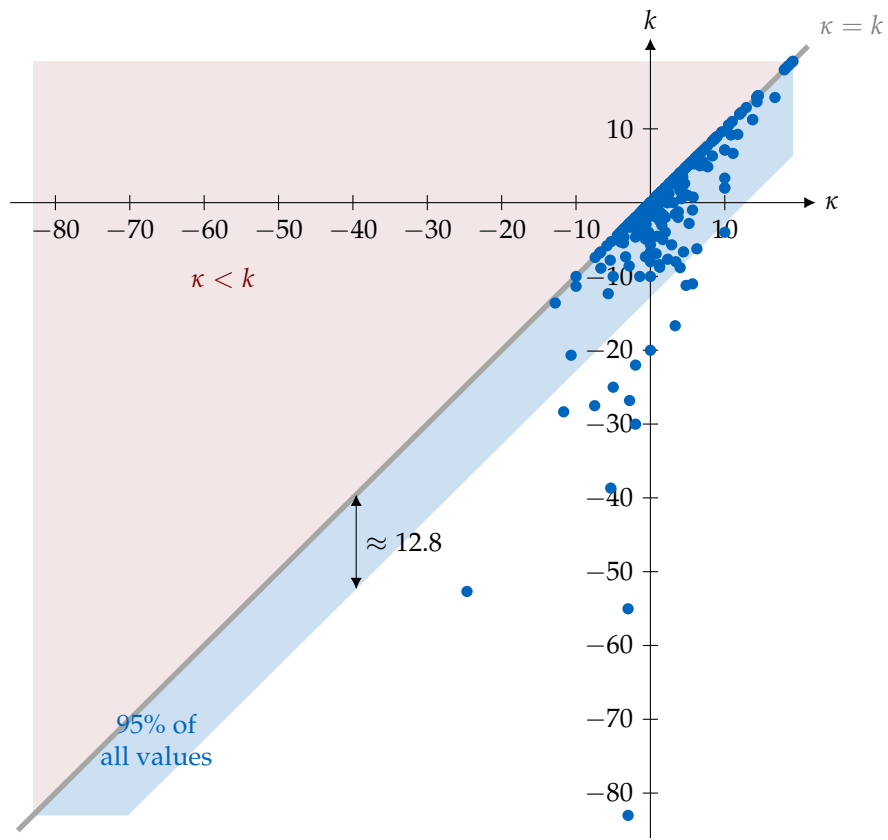

  \begin{center}

  \end{center}
  \caption{Comparison of $\kappa(r,s)$ and $k(r,s)$ for $300$ randomly selected state pairs
  of the CTMC arising from the RSVP model with $M = 7$, $N = 5$ and $3$ mobile nodes (resulting in $842$ states)}
  \label{fig:rsvp_kappa_k_comp}
\end{figure}

We now try to understand why the error bounds are so far from the actual measured errors.
First, we want to get an idea of how loose the bounds $k(r,s)$ are for the Ricci curvature.
In order to do that,
we randomly selected $300$ state pairs for which we calculated both $\kappa(r,s)$
and $k(r,s)$. The result is shown in \autoref{fig:rsvp_kappa_k_comp}.
We can see that for most of the sampled pairs, $k(r,s)$ is actually quite close to $\kappa(r,s)$. However,
there are some pairs where the bound $k(r,s)$ is significantly lower than $\kappa(r, s)$. If this happens
to be the case for the state pair where $\kappa(r,s)$ attains the minimum, $\underline{k}(Q)$ will be
much lower than $\underline{\kappa}(Q)$, which is the case in the RSVP model. Overall, $k(r,s)$ seems to be a reasonable
compromise between computation time and a tight bound, but finding better bounds would
still be a worthy research subject.

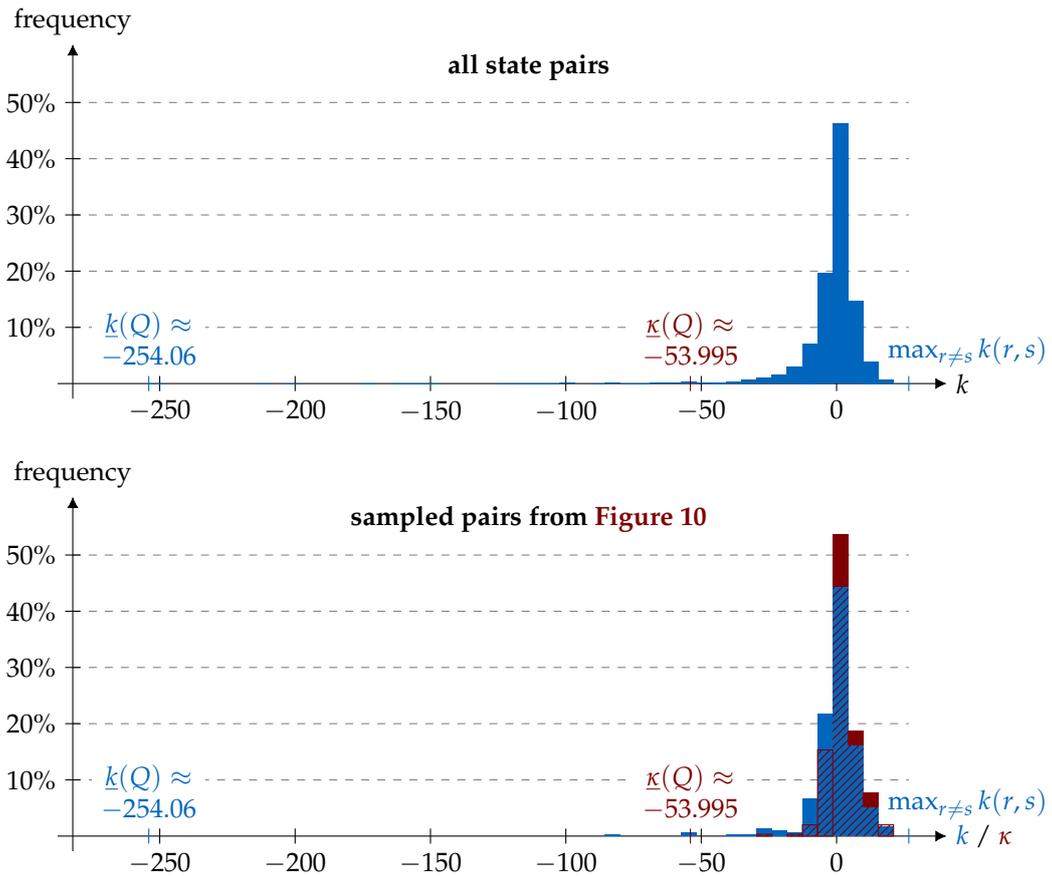
\begin{figure}[htb]
  \begin{center}
    \begin{tikzpicture}[>={Latex[length=1.5mm,width=1.5mm]}]
      \node at (5, 4.2) {\textbf{all state pairs}};
      \node at (5, -1.8) {\textbf{sampled pairs from \autoref{fig:rsvp_kappa_k_comp}}};
      \draw[gray,dashed] (-1,0.74534) -- (10,0.74534);
      \draw[gray,dashed] (-1,-5.2547) -- (10,-5.2547);
      \draw[gray,dashed] (-1,1.4907) -- (10,1.4907);
      \draw[gray,dashed] (-1,-4.5093) -- (10,-4.5093);
      \draw[gray,dashed] (-1,2.236) -- (10,2.236);
      \draw[gray,dashed] (-1,-3.764) -- (10,-3.764);
      \draw[gray,dashed] (-1,2.9814) -- (10,2.9814);
      \draw[gray,dashed] (-1,-3.0186) -- (10,-3.0186);
      \draw[gray,dashed] (-1,3.7267) -- (10,3.7267);
      \draw[gray,dashed] (-1,-2.2733) -- (10,-2.2733);
      \draw[->] (-1.2,0) -- (10.5,0) node[right] {$k$};
      \draw[->] (-1,-0.2) -- (-1,4.5) node[above] {frequency};
      \draw[sectionblue] (0,-0.1) -- (0,0.1) node[above,align=center,fill=white] {$\underline{k}(Q) \approx$\\$-254.06$};
      \draw[sectionblue] (10,-0.1) -- (10,0.1) node[above right=0mm and -4mm] {$\max_{r \neq s} k(r,s)$};
      \draw[->] (-1.2,-6) -- (10.5,-6) node[right] {$\textcolor{sectionblue}{k}$ / $\textcolor{linkred}{\kappa}$};
      \draw[->] (-1,-6.2) -- (-1,-1.5) node[above] {frequency};
      \draw[sectionblue] (0,-6.1) -- (0,-5.9) node[above,align=center,fill=white] {$\underline{k}(Q) \approx$\\$-254.06$};
      \draw[sectionblue] (10,-6.1) -- (10,-5.9) node[above right=0mm and -4mm] {$\max_{r \neq s} k(r,s)$};
      \draw (-0.9,0.74534) -- (-1.1,0.74534) node[left] {$10\%$};
      \draw (-0.9,-5.2547) -- (-1.1,-5.2547) node[left] {$10\%$};
      \draw (-0.9,1.4907) -- (-1.1,1.4907) node[left] {$20\%$};
      \draw (-0.9,-4.5093) -- (-1.1,-4.5093) node[left] {$20\%$};
      \draw (-0.9,2.236) -- (-1.1,2.236) node[left] {$30\%$};
      \draw (-0.9,-3.764) -- (-1.1,-3.764) node[left] {$30\%$};
      \draw (-0.9,2.9814) -- (-1.1,2.9814) node[left] {$40\%$};
      \draw (-0.9,-3.0186) -- (-1.1,-3.0186) node[left] {$40\%$};
      \draw (-0.9,3.7267) -- (-1.1,3.7267) node[left] {$50\%$};
      \draw (-0.9,-2.2733) -- (-1.1,-2.2733) node[left] {$50\%$};
      \draw (0.14445,0.1) -- (0.14445,-0.1) node[below] {$-250$};
      \draw (0.14445,-5.9) -- (0.14445,-6.1) node[below] {$-250$};
      \draw (1.9256,0.1) -- (1.9256,-0.1) node[below] {$-200$};
      \draw (1.9256,-5.9) -- (1.9256,-6.1) node[below] {$-200$};
      \draw (3.7067,0.1) -- (3.7067,-0.1) node[below] {$-150$};
      \draw (3.7067,-5.9) -- (3.7067,-6.1) node[below] {$-150$};
      \draw (5.4879,0.1) -- (5.4879,-0.1) node[below] {$-100$};
      \draw (5.4879,-5.9) -- (5.4879,-6.1) node[below] {$-100$};
      \draw (7.269,0.1) -- (7.269,-0.1) node[below] {$-50$};
      \draw (7.269,-5.9) -- (7.269,-6.1) node[below] {$-50$};
      \draw (9.0502,0.1) -- (9.0502,-0.1) node[below] {$0$};
      \draw (9.0502,-5.9) -- (9.0502,-6.1) node[below] {$0$};
      \fill[sectionblue] (0,0) -- (0.2,0) -- (0.2,0.0028419) -- (0,0.0028419) -- cycle;
      \fill[sectionblue] (1.4,0) -- (1.6,0) -- (1.6,0.0030314) -- (1.4,0.0030314) -- cycle;
      \fill[sectionblue] (1.8,0) -- (2,0) -- (2,6.3154e-05) -- (1.8,6.3154e-05) -- cycle;
      \fill[sectionblue] (2.8,0) -- (3,0) -- (3,0.000821) -- (2.8,0.000821) -- cycle;
      \fill[sectionblue] (3.2,0) -- (3.4,0) -- (3.4,0.0012631) -- (3.2,0.0012631) -- cycle;
      \fill[sectionblue] (3.4,0) -- (3.6,0) -- (3.6,0.0031577) -- (3.4,0.0031577) -- cycle;
      \fill[sectionblue] (3.6,0) -- (3.8,0) -- (3.8,0.001642) -- (3.6,0.001642) -- cycle;
      \fill[sectionblue] (4.6,0) -- (4.8,0) -- (4.8,0.00063154) -- (4.6,0.00063154) -- cycle;
      \fill[sectionblue] (4.8,0) -- (5,0) -- (5,0.0015157) -- (4.8,0.0015157) -- cycle;
      \fill[sectionblue] (5,0) -- (5.2,0) -- (5.2,0.00088415) -- (5,0.00088415) -- cycle;
      \fill[sectionblue] (5.2,0) -- (5.4,0) -- (5.4,0.0026525) -- (5.2,0.0026525) -- cycle;
      \fill[sectionblue] (5.4,0) -- (5.6,0) -- (5.6,0.0048628) -- (5.4,0.0048628) -- cycle;
      \fill[sectionblue] (5.8,0) -- (6,0) -- (6,0.0012631) -- (5.8,0.0012631) -- cycle;
      \fill[sectionblue] (6,0) -- (6.2,0) -- (6.2,0.0071995) -- (6,0.0071995) -- cycle;
      \fill[sectionblue] (6,-6) -- (6.2,-6) -- (6.2,-5.9752) -- (6,-5.9752) -- cycle;
      \fill[sectionblue] (6.2,0) -- (6.4,0) -- (6.4,0.00037892) -- (6.2,0.00037892) -- cycle;
      \fill[sectionblue] (6.4,0) -- (6.6,0) -- (6.6,0.0017051) -- (6.4,0.0017051) -- cycle;
      \fill[sectionblue] (6.6,0) -- (6.8,0) -- (6.8,0.0035998) -- (6.6,0.0035998) -- cycle;
      \fill[sectionblue] (6.8,0) -- (7,0) -- (7,0.005768) -- (6.8,0.005768) -- cycle;
      \fill[sectionblue] (7,0) -- (7.2,0) -- (7.2,0.025304) -- (7,0.025304) -- cycle;
      \fill[sectionblue] (7,-6) -- (7.2,-6) -- (7.2,-5.9503) -- (7,-5.9503) -- cycle;
      \fill[sectionblue] (7.2,0) -- (7.4,0) -- (7.4,0.0081047) -- (7.2,0.0081047) -- cycle;
      \fill[sectionblue] (7.4,0) -- (7.6,0) -- (7.6,0.01642) -- (7.4,0.01642) -- cycle;
      \fill[sectionblue] (7.6,0) -- (7.8,0) -- (7.8,0.01863) -- (7.6,0.01863) -- cycle;
      \fill[sectionblue] (7.6,-6) -- (7.8,-6) -- (7.8,-5.9752) -- (7.6,-5.9752) -- cycle;
      \fill[sectionblue] (7.8,0) -- (8,0) -- (8,0.045386) -- (7.8,0.045386) -- cycle;
      \fill[sectionblue] (7.8,-6) -- (8,-6) -- (8,-5.9752) -- (7.8,-5.9752) -- cycle;
      \fill[sectionblue] (8,0) -- (8.2,0) -- (8.2,0.073974) -- (8,0.073974) -- cycle;
      \fill[sectionblue] (8,-6) -- (8.2,-6) -- (8.2,-5.9006) -- (8,-5.9006) -- cycle;
      \fill[pattern=north east lines,pattern color=linkred] (8,-6) -- (8.2,-6) -- (8.2,-5.9752) -- (8,-5.9752) -- cycle;
      \fill[sectionblue] (8.2,0) -- (8.4,0) -- (8.4,0.11591) -- (8.2,0.11591) -- cycle;
      \fill[sectionblue] (8.2,-6) -- (8.4,-6) -- (8.4,-5.9255) -- (8.2,-5.9255) -- cycle;
      \fill[sectionblue] (8.4,0) -- (8.6,0) -- (8.6,0.21996) -- (8.4,0.21996) -- cycle;
      \fill[sectionblue] (8.4,-6) -- (8.6,-6) -- (8.6,-5.9503) -- (8.4,-5.9503) -- cycle;
      \fill[pattern=north east lines,pattern color=linkred] (8.4,-6) -- (8.6,-6) -- (8.6,-5.9752) -- (8.4,-5.9752) -- cycle;
      \fill[sectionblue] (8.6,0) -- (8.8,0) -- (8.8,0.5306) -- (8.6,0.5306) -- cycle;
      \fill[sectionblue] (8.6,-6) -- (8.8,-6) -- (8.8,-5.5031) -- (8.6,-5.5031) -- cycle;
      \fill[pattern=north east lines,pattern color=linkred] (8.6,-6) -- (8.8,-6) -- (8.8,-5.8509) -- (8.6,-5.8509) -- cycle;
      \fill[sectionblue] (8.8,0) -- (9,0) -- (9,1.4652) -- (8.8,1.4652) -- cycle;
      \fill[sectionblue] (8.8,-6) -- (9,-6) -- (9,-4.3851) -- (8.8,-4.3851) -- cycle;
      \fill[pattern=north east lines,pattern color=linkred] (8.8,-6) -- (9,-6) -- (9,-4.8571) -- (8.8,-4.8571) -- cycle;
      \fill[sectionblue] (9,0) -- (9.2,0) -- (9.2,3.4498) -- (9,3.4498) -- cycle;
      \fill[linkred] (9,-6) -- (9.2,-6) -- (9.2,-2) -- (9,-2) -- cycle;
      \fill[sectionblue] (9,-6) -- (9.2,-6) -- (9.2,-2.6957) -- (9,-2.6957) -- cycle;
      \fill[pattern=north east lines,pattern color=linkred] (9,-6) -- (9.2,-6) -- (9.2,-2.6957) -- (9,-2.6957) -- cycle;
      \fill[sectionblue] (9.2,0) -- (9.4,0) -- (9.4,1.0959) -- (9.2,1.0959) -- cycle;
      \fill[linkred] (9.2,-6) -- (9.4,-6) -- (9.4,-4.6087) -- (9.2,-4.6087) -- cycle;
      \fill[sectionblue] (9.2,-6) -- (9.4,-6) -- (9.4,-4.8075) -- (9.2,-4.8075) -- cycle;
      \fill[pattern=north east lines,pattern color=linkred] (9.2,-6) -- (9.4,-6) -- (9.4,-4.8075) -- (9.2,-4.8075) -- cycle;
      \fill[sectionblue] (9.4,0) -- (9.6,0) -- (9.6,0.29299) -- (9.4,0.29299) -- cycle;
      \fill[linkred] (9.4,-6) -- (9.6,-6) -- (9.6,-5.4286) -- (9.4,-5.4286) -- cycle;
      \fill[sectionblue] (9.4,-6) -- (9.6,-6) -- (9.6,-5.6273) -- (9.4,-5.6273) -- cycle;
      \fill[pattern=north east lines,pattern color=linkred] (9.4,-6) -- (9.6,-6) -- (9.6,-5.6273) -- (9.4,-5.6273) -- cycle;
      \fill[sectionblue] (9.6,0) -- (9.8,0) -- (9.8,0.051028) -- (9.6,0.051028) -- cycle;
      \fill[sectionblue] (9.6,-6) -- (9.8,-6) -- (9.8,-5.8758) -- (9.6,-5.8758) -- cycle;
      \fill[pattern=north east lines,pattern color=linkred] (9.6,-6) -- (9.8,-6) -- (9.8,-5.8509) -- (9.6,-5.8509) -- cycle;
      \fill[sectionblue] (9.8,0) -- (10,0) -- (10,0.00088415) -- (9.8,0.00088415) -- cycle;
      \draw[linkred] (8,-6) -- (8.2,-6) -- (8.2,-5.9752) -- (8,-5.9752) -- cycle;
      \draw[linkred] (8.4,-6) -- (8.6,-6) -- (8.6,-5.9752) -- (8.4,-5.9752) -- cycle;
      \draw[linkred] (8.6,-6) -- (8.8,-6) -- (8.8,-5.8509) -- (8.6,-5.8509) -- cycle;
      \draw[linkred] (8.8,-6) -- (9,-6) -- (9,-4.8571) -- (8.8,-4.8571) -- cycle;
      \draw[linkred] (9.6,-6) -- (9.8,-6) -- (9.8,-5.8509) -- (9.6,-5.8509) -- cycle;
      \draw[linkred] (7.1267,-6.1) -- (7.1267,-5.9) node[above,align=center,fill=white] {$\underline{\kappa}(Q) \approx$\\$-53.995$};
      \draw[linkred] (7.1267,-0.1) -- (7.1267,0.1) node[above,align=center,fill=white] {$\underline{\kappa}(Q) \approx$\\$-53.995$};
    \end{tikzpicture}
  \end{center}
  \caption{Histogram showing the frequency of $k$ and $\kappa$ values for the CTMC
  arising from the RSVP model with $M = 7$, $N = 5$ and $3$ mobile nodes (resulting in $842$ states).
  The upper histogram shows how the values $k(r,s)$ for all state pairs $r,s$ are
  distributed. The lower histogram shows the distribution of both $k(r,s)$ and $\kappa(r,s)$, but
  only for the $300$ randomly selected pairs from \autoref{fig:rsvp_kappa_k_comp}
  ($k$ values in blue, $\kappa$ values in red).}
  \label{fig:rsvp_k_histogram}
\end{figure}

Next, we look at how many state pairs have $k(r,s)$ or $\kappa(r,s)$ near the minimum
values $\underline{k}(Q)$ and $\underline{\kappa}(Q)$.
In \autoref{fig:rsvp_k_histogram}, a histogram of the distribution of $k(r,s)$ over
all state pairs is shown, together with a histogram of the distribution of $k(r,s)$ and
$\kappa(r,s)$ for the state pairs selected for \autoref{fig:rsvp_kappa_k_comp}. We can
see that most values cluster around $0$.
The fraction of state pairs attaining a $k(r,s)$ with $k(r,s) < \underline{\kappa}(Q)$ is
so small that they are not visible in the histogram. That is, an almost negligible part of the
state pairs is responsible for the very low $\underline{k}(Q)$. Even near $\underline{\kappa}(Q) \approx -53.995$,
there are no visible bars in the histogram. The bars only become visible around $k \approx -30$.
It might be possible to exploit this (only a negligible fraction of
state pairs actually having $k(r,s)$ very close to $\underline{k}(Q)$, and
the same for $\kappa$) to achieve better error bounds,
even though it is not evident at all how~-- we don't want to compute $p_t$ exactly and
it therefore seems that the bound from \autoref{lem:wdderiv_error_ricciinfbound} needs
to hold for all probability distributions, and it is actually tight in that case.

The RSVP model shows that improvements over the current bounds are necessary
to achieve useful error bounds for this particular example and the chosen metric.

\subsubsection{Further examples}
\label{sssec:further_examples}

We also considered the workstation cluster model from \cite{modcheckdependability}. It consists of two clusters
of workstations connected by switches, where each workstation and switch can break down
and a repair unit can repair failed components. With $4$ workstation in each of the two clusters,
the model has $820$ states. Again, the model was not originally defined in
conjunction with a metric. We simply used the sum of the absolute value of the
differences of the single state components (sometimes multiplied with a factor of $0.5$
because the state of the repair unit is encoded in more than one state component,
which leads to redundant information in the state encoding),
resulting in a state space diameter of $12$.

\begin{figure}[htb]
  \begin{center}
    \begin{tikzpicture}[>={Latex[length=1.5mm,width=1.5mm]}]
      \node at (5, 4.3) {\textbf{all state pairs}};
      \node at (5, -1.7) {\textbf{300 sampled pairs}};
      \draw[gray,dashed] (-1,0.56075) -- (10,0.56075);
      \draw[gray,dashed] (-1,-5.4393) -- (10,-5.4393);
      \draw[gray,dashed] (-1,1.1215) -- (10,1.1215);
      \draw[gray,dashed] (-1,-4.8785) -- (10,-4.8785);
      \draw[gray,dashed] (-1,1.6822) -- (10,1.6822);
      \draw[gray,dashed] (-1,-4.3178) -- (10,-4.3178);
      \draw[gray,dashed] (-1,2.243) -- (10,2.243);
      \draw[gray,dashed] (-1,-3.757) -- (10,-3.757);
      \draw[gray,dashed] (-1,2.8037) -- (10,2.8037);
      \draw[gray,dashed] (-1,-3.1963) -- (10,-3.1963);
      \draw[gray,dashed] (-1,3.3645) -- (10,3.3645);
      \draw[gray,dashed] (-1,-2.6355) -- (10,-2.6355);
      \draw[gray,dashed] (-1,3.9252) -- (10,3.9252);
      \draw[gray,dashed] (-1,-2.0748) -- (10,-2.0748);
      \draw[->] (-1.2,0) -- (10.5,0) node[right] {$k$};
      \draw[->] (-1,-0.2) -- (-1,4.5) node[above] {frequency};
      \draw[sectionblue] (0,-0.1) -- (0,0.1) node[above,align=center,fill=white] {$\underline{k}(Q) \approx$\\$-100.01$};
      \draw[sectionblue] (10,-0.1) -- (10,0.1) node[above right=0mm and -4mm] {$\max_{r \neq s} k(r,s)$};
      \draw[->] (-1.2,-6) -- (10.5,-6) node[right] {$\textcolor{sectionblue}{k}$ / $\textcolor{linkred}{\kappa}$};
      \draw[->] (-1,-6.2) -- (-1,-1.5) node[above] {frequency};
      \draw[sectionblue] (0,-6.1) -- (0,-5.9) node[above,align=center,fill=white] {$\underline{k}(Q) \approx$\\$-100.01$};
      \draw[sectionblue] (10,-6.1) -- (10,-5.9) node[above right=0mm and -4mm] {$\max_{r \neq s} k(r,s)$};
      \draw (-0.9,0.56075) -- (-1.1,0.56075) node[left] {$10\%$};
      \draw (-0.9,-5.4393) -- (-1.1,-5.4393) node[left] {$10\%$};
      \draw (-0.9,1.1215) -- (-1.1,1.1215) node[left] {$20\%$};
      \draw (-0.9,-4.8785) -- (-1.1,-4.8785) node[left] {$20\%$};
      \draw (-0.9,1.6822) -- (-1.1,1.6822) node[left] {$30\%$};
      \draw (-0.9,-4.3178) -- (-1.1,-4.3178) node[left] {$30\%$};
      \draw (-0.9,2.243) -- (-1.1,2.243) node[left] {$40\%$};
      \draw (-0.9,-3.757) -- (-1.1,-3.757) node[left] {$40\%$};
      \draw (-0.9,2.8037) -- (-1.1,2.8037) node[left] {$50\%$};
      \draw (-0.9,-3.1963) -- (-1.1,-3.1963) node[left] {$50\%$};
      \draw (-0.9,3.3645) -- (-1.1,3.3645) node[left] {$60\%$};
      \draw (-0.9,-2.6355) -- (-1.1,-2.6355) node[left] {$60\%$};
      \draw (-0.9,3.9252) -- (-1.1,3.9252) node[left] {$70\%$};
      \draw (-0.9,-2.0748) -- (-1.1,-2.0748) node[left] {$70\%$};
      \draw (0.00090901,0.1) -- (0.00090901,-0.1) node[below] {$-100$};
      \draw (0.00090901,-5.9) -- (0.00090901,-6.1) node[below] {$-100$};
      \draw (0.90992,0.1) -- (0.90992,-0.1) node[below] {$-90$};
      \draw (0.90992,-5.9) -- (0.90992,-6.1) node[below] {$-90$};
      \draw (1.8189,0.1) -- (1.8189,-0.1) node[below] {$-80$};
      \draw (1.8189,-5.9) -- (1.8189,-6.1) node[below] {$-80$};
      \draw (2.7279,0.1) -- (2.7279,-0.1) node[below] {$-70$};
      \draw (2.7279,-5.9) -- (2.7279,-6.1) node[below] {$-70$};
      \draw (3.6369,0.1) -- (3.6369,-0.1) node[below] {$-60$};
      \draw (3.6369,-5.9) -- (3.6369,-6.1) node[below] {$-60$};
      \draw (4.546,0.1) -- (4.546,-0.1) node[below] {$-50$};
      \draw (4.546,-5.9) -- (4.546,-6.1) node[below] {$-50$};
      \draw (5.455,0.1) -- (5.455,-0.1) node[below] {$-40$};
      \draw (5.455,-5.9) -- (5.455,-6.1) node[below] {$-40$};
      \draw (6.364,0.1) -- (6.364,-0.1) node[below] {$-30$};
      \draw (6.364,-5.9) -- (6.364,-6.1) node[below] {$-30$};
      \draw (7.273,0.1) -- (7.273,-0.1) node[below] {$-20$};
      \draw (7.273,-5.9) -- (7.273,-6.1) node[below] {$-20$};
      \draw (8.182,0.1) -- (8.182,-0.1) node[below] {$-10$};
      \draw (8.182,-5.9) -- (8.182,-6.1) node[below] {$-10$};
      \draw (9.091,0.1) -- (9.091,-0.1) node[below] {$0$};
      \draw (9.091,-5.9) -- (9.091,-6.1) node[below] {$0$};
      \draw (10,0.1) -- (10,-0.1) node[below] {$10$};
      \draw (10,-5.9) -- (10,-6.1) node[below] {$10$};
      \fill[sectionblue] (0,0) -- (0.2,0) -- (0.2,0.00040078) -- (0,0.00040078) -- cycle;
      \fill[sectionblue] (0.8,0) -- (1,0) -- (1,0.00093516) -- (0.8,0.00093516) -- cycle;
      \fill[sectionblue] (1.8,0) -- (2,0) -- (2,0.0013025) -- (1.8,0.0013025) -- cycle;
      \fill[sectionblue] (2.6,0) -- (2.8,0) -- (2.8,0.0024381) -- (2.6,0.0024381) -- cycle;
      \fill[sectionblue] (3.6,0) -- (3.8,0) -- (3.8,0.0015029) -- (3.6,0.0015029) -- cycle;
      \fill[sectionblue] (4.4,0) -- (4.6,0) -- (4.6,0.002722) -- (4.4,0.002722) -- cycle;
      \fill[sectionblue] (5,0) -- (5.2,0) -- (5.2,0.0027387) -- (5,0.0027387) -- cycle;
      \fill[sectionblue] (5.4,0) -- (5.6,0) -- (5.6,0.0049096) -- (5.4,0.0049096) -- cycle;
      \fill[sectionblue] (5.8,0) -- (6,0) -- (6,0.006446) -- (5.8,0.006446) -- cycle;
      \fill[sectionblue] (5.8,-6) -- (6,-6) -- (6,-5.9813) -- (5.8,-5.9813) -- cycle;
      \fill[sectionblue] (6,0) -- (6.2,0) -- (6.2,0.00053438) -- (6,0.00053438) -- cycle;
      \fill[sectionblue] (6.2,0) -- (6.4,0) -- (6.4,0.010821) -- (6.2,0.010821) -- cycle;
      \fill[sectionblue] (6.6,0) -- (6.8,0) -- (6.8,0.0085167) -- (6.6,0.0085167) -- cycle;
      \fill[sectionblue] (6.6,-6) -- (6.8,-6) -- (6.8,-5.9626) -- (6.6,-5.9626) -- cycle;
      \fill[sectionblue] (6.8,0) -- (7,0) -- (7,0.016566) -- (6.8,0.016566) -- cycle;
      \fill[sectionblue] (7,0) -- (7.2,0) -- (7.2,0.0037407) -- (7,0.0037407) -- cycle;
      \fill[sectionblue] (7.2,0) -- (7.4,0) -- (7.4,0.026051) -- (7.2,0.026051) -- cycle;
      \fill[sectionblue] (7.4,0) -- (7.6,0) -- (7.6,0.026953) -- (7.4,0.026953) -- cycle;
      \fill[sectionblue] (7.4,-6) -- (7.6,-6) -- (7.6,-5.9439) -- (7.4,-5.9439) -- cycle;
      \fill[sectionblue] (7.6,0) -- (7.8,0) -- (7.8,0.028389) -- (7.6,0.028389) -- cycle;
      \fill[sectionblue] (7.6,-6) -- (7.8,-6) -- (7.8,-5.9626) -- (7.6,-5.9626) -- cycle;
      \fill[sectionblue] (7.8,0) -- (8,0) -- (8,0.038609) -- (7.8,0.038609) -- cycle;
      \fill[sectionblue] (7.8,-6) -- (8,-6) -- (8,-5.9626) -- (7.8,-5.9626) -- cycle;
      \fill[sectionblue] (8,0) -- (8.2,0) -- (8.2,0.0766) -- (8,0.0766) -- cycle;
      \fill[sectionblue] (8,-6) -- (8.2,-6) -- (8.2,-5.9626) -- (8,-5.9626) -- cycle;
      \fill[sectionblue] (8.2,0) -- (8.4,0) -- (8.4,0.043752) -- (8.2,0.043752) -- cycle;
      \fill[sectionblue] (8.2,-6) -- (8.4,-6) -- (8.4,-5.9439) -- (8.2,-5.9439) -- cycle;
      \fill[sectionblue] (8.4,0) -- (8.6,0) -- (8.6,0.033365) -- (8.4,0.033365) -- cycle;
      \fill[sectionblue] (8.4,-6) -- (8.6,-6) -- (8.6,-5.9252) -- (8.4,-5.9252) -- cycle;
      \fill[pattern=north east lines,pattern color=linkred] (8.4,-6) -- (8.6,-6) -- (8.6,-5.9439) -- (8.4,-5.9439) -- cycle;
      \fill[sectionblue] (8.6,0) -- (8.8,0) -- (8.8,0.037373) -- (8.6,0.037373) -- cycle;
      \fill[sectionblue] (8.6,-6) -- (8.8,-6) -- (8.8,-5.9813) -- (8.6,-5.9813) -- cycle;
      \fill[pattern=north east lines,pattern color=linkred] (8.6,-6) -- (8.8,-6) -- (8.8,-5.8879) -- (8.6,-5.8879) -- cycle;
      \fill[sectionblue] (8.8,0) -- (9,0) -- (9,0.23553) -- (8.8,0.23553) -- cycle;
      \fill[sectionblue] (8.8,-6) -- (9,-6) -- (9,-5.757) -- (8.8,-5.757) -- cycle;
      \fill[pattern=north east lines,pattern color=linkred] (8.8,-6) -- (9,-6) -- (9,-5.8131) -- (8.8,-5.8131) -- cycle;
      \fill[sectionblue] (9,0) -- (9.2,0) -- (9.2,3.6688) -- (9,3.6688) -- cycle;
      \fill[linkred] (9,-6) -- (9.2,-6) -- (9.2,-2) -- (9,-2) -- cycle;
      \fill[sectionblue] (9,-6) -- (9.2,-6) -- (9.2,-2.2243) -- (9,-2.2243) -- cycle;
      \fill[pattern=north east lines,pattern color=linkred] (9,-6) -- (9.2,-6) -- (9.2,-2.2243) -- (9,-2.2243) -- cycle;
      \fill[sectionblue] (9.2,0) -- (9.4,0) -- (9.4,1.1887) -- (9.2,1.1887) -- cycle;
      \fill[sectionblue] (9.2,-6) -- (9.4,-6) -- (9.4,-4.8785) -- (9.2,-4.8785) -- cycle;
      \fill[pattern=north east lines,pattern color=linkred] (9.2,-6) -- (9.4,-6) -- (9.4,-4.8411) -- (9.2,-4.8411) -- cycle;
      \fill[sectionblue] (9.4,0) -- (9.6,0) -- (9.6,0.12545) -- (9.4,0.12545) -- cycle;
      \fill[sectionblue] (9.4,-6) -- (9.6,-6) -- (9.6,-5.9065) -- (9.4,-5.9065) -- cycle;
      \fill[pattern=north east lines,pattern color=linkred] (9.4,-6) -- (9.6,-6) -- (9.6,-5.9065) -- (9.4,-5.9065) -- cycle;
      \fill[sectionblue] (9.6,0) -- (9.8,0) -- (9.8,0.0040078) -- (9.6,0.0040078) -- cycle;
      \fill[sectionblue] (9.8,0) -- (10,0) -- (10,0.010354) -- (9.8,0.010354) -- cycle;
      \draw[linkred] (8.4,-6) -- (8.6,-6) -- (8.6,-5.9439) -- (8.4,-5.9439) -- cycle;
      \draw[linkred] (8.6,-6) -- (8.8,-6) -- (8.8,-5.8879) -- (8.6,-5.8879) -- cycle;
      \draw[linkred] (8.8,-6) -- (9,-6) -- (9,-5.8131) -- (8.8,-5.8131) -- cycle;
      \draw[linkred] (9.2,-6) -- (9.4,-6) -- (9.4,-4.8411) -- (9.2,-4.8411) -- cycle;
      \draw[linkred] (9.4,-6) -- (9.6,-6) -- (9.6,-5.9065) -- (9.4,-5.9065) -- cycle;
      \draw[linkred] (8.182,-6.1) -- (8.182,-5.9) node[above left=0mm and -3mm,align=center,fill=white] {$\underline{\kappa}(Q) \approx$\\$-9.9998$};
      \draw[linkred] (8.182,-0.1) -- (8.182,0.1) node[above left=0mm and -3mm,align=center,fill=white] {$\underline{\kappa}(Q) \approx$\\$-9.9998$};
    \end{tikzpicture}
  \end{center}
  \caption{Histogram showing the frequency of $k$ and $\kappa$ values for the CTMC
  arising from the workstation cluster model with $4$ workstations per cluster (resulting in $820$ states).
  The upper histogram shows how the values $k(r,s)$ for all state pairs $r,s$ are
  distributed. The lower histogram shows the distribution of both $k(r,s)$ and $\kappa(r,s)$, but
  only for $300$ randomly selected pairs
  ($k$ values in blue, $\kappa$ values in red).}
  \label{fig:wrkstcls_k_histogram}
\end{figure}
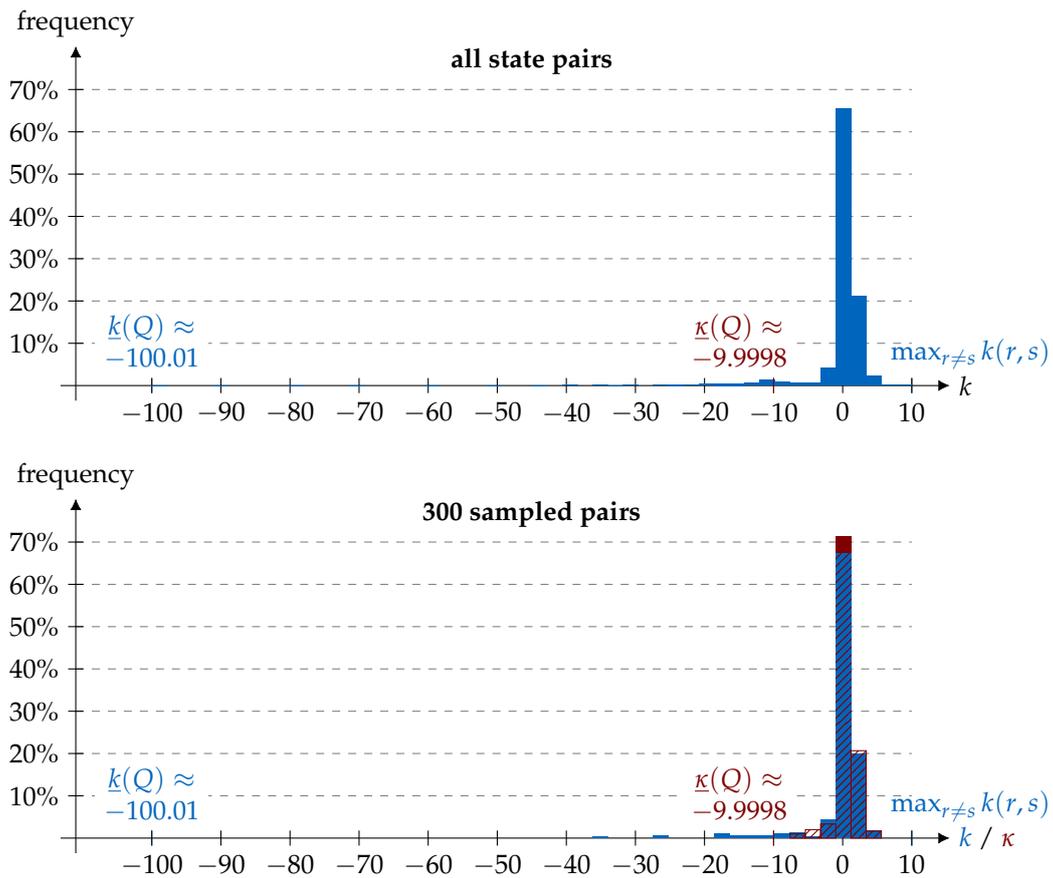

For the workstation cluster model, we get $\underline{k}(Q) \approx -100$, $K(Q) \approx 100$
and $\underline{\kappa}(Q) \approx -10$ (the computation of the latter taking around $15$ hours
on our machine, using the same strategy as with the RSVP model, requiring the calculation
of $\kappa(r,s)$ for $4.6\%$ of all pairs).
\autoref{fig:wrkstcls_k_histogram} yields a similar picture to what we saw for the
RSVP model: only a negligible number of the state pairs cause the low value
of $\underline{k}(Q)$. The only difference is that $\underline{\kappa}(Q)$ is not
quite as low as it was for the RSVP model. In our experiments, it was still too low
to yield a practically useful error bound for the transient errors when we aggregate
the model. While the actual error $\WD{\widetilde{p}_t}{p_t}$ is only $\approx 0.001$ at time $t = 20$
and thus quite small (for an aggregation with $161$ aggregates, resulting
in $\norm{\Theta A - A Q}_{\mathrm{W}} \approx 0.08$), the error bound using $\underline{\kappa}(Q)$ hits the state space
diameter already around time $0.73$.

\clearpage

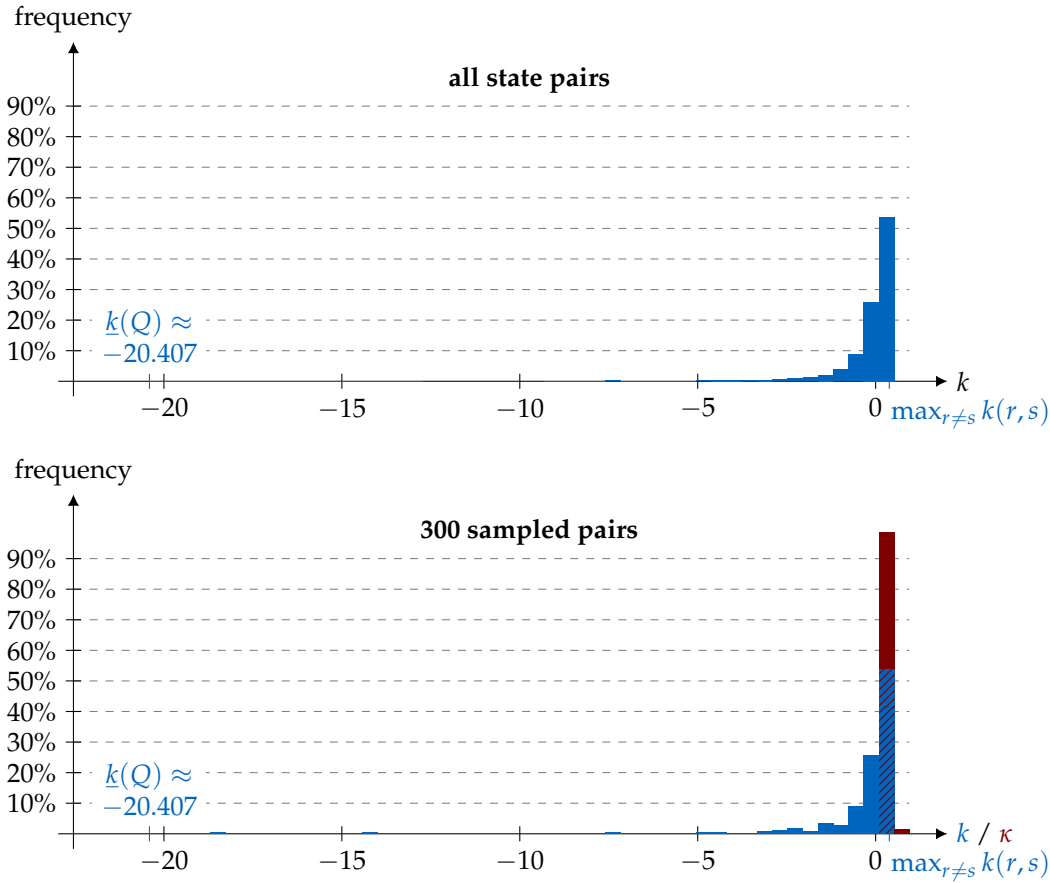
\begin{figure}[htb]
  \begin{center}
    \begin{tikzpicture}[>={Latex[length=1.5mm,width=1.5mm]}]
      \node at (5, 4) {\textbf{all state pairs}};
      \node at (5, -2) {\textbf{300 sampled pairs}};
      \draw[gray,dashed] (-1,0.40541) -- (10,0.40541);
      \draw[gray,dashed] (-1,-5.5946) -- (10,-5.5946);
      \draw[gray,dashed] (-1,0.81081) -- (10,0.81081);
      \draw[gray,dashed] (-1,-5.1892) -- (10,-5.1892);
      \draw[gray,dashed] (-1,1.2162) -- (10,1.2162);
      \draw[gray,dashed] (-1,-4.7838) -- (10,-4.7838);
      \draw[gray,dashed] (-1,1.6216) -- (10,1.6216);
      \draw[gray,dashed] (-1,-4.3784) -- (10,-4.3784);
      \draw[gray,dashed] (-1,2.027) -- (10,2.027);
      \draw[gray,dashed] (-1,-3.973) -- (10,-3.973);
      \draw[gray,dashed] (-1,2.4324) -- (10,2.4324);
      \draw[gray,dashed] (-1,-3.5676) -- (10,-3.5676);
      \draw[gray,dashed] (-1,2.8378) -- (10,2.8378);
      \draw[gray,dashed] (-1,-3.1622) -- (10,-3.1622);
      \draw[gray,dashed] (-1,3.2432) -- (10,3.2432);
      \draw[gray,dashed] (-1,-2.7568) -- (10,-2.7568);
      \draw[gray,dashed] (-1,3.6486) -- (10,3.6486);
      \draw[gray,dashed] (-1,-2.3514) -- (10,-2.3514);
      \draw[->] (-1.2,0) -- (10.5,0) node[right] {$k$};
      \draw[->] (-1,-0.2) -- (-1,4.5) node[above] {frequency};
      \draw[sectionblue] (0,-0.1) -- (0,0.1) node[above,align=center,fill=white] {$\underline{k}(Q) \approx$\\$-20.407$};
      \draw[sectionblue] (9.7352,0.1) -- (9.7352,-0.1) node[below right=0mm and -1mm] {$\max_{r \neq s} k(r,s)$};
      \draw[->] (-1.2,-6) -- (10.5,-6) node[right] {$\textcolor{sectionblue}{k}$ / $\textcolor{linkred}{\kappa}$};
      \draw[->] (-1,-6.2) -- (-1,-1.5) node[above] {frequency};
      \draw[sectionblue] (0,-6.1) -- (0,-5.9) node[above,align=center,fill=white] {$\underline{k}(Q) \approx$\\$-20.407$};
      \draw[sectionblue] (9.7352,-5.9) -- (9.7352,-6.1) node[below right=0mm and -1mm] {$\max_{r \neq s} k(r,s)$};
      \draw (-0.9,0.40541) -- (-1.1,0.40541) node[left] {$10\%$};
      \draw (-0.9,-5.5946) -- (-1.1,-5.5946) node[left] {$10\%$};
      \draw (-0.9,0.81081) -- (-1.1,0.81081) node[left] {$20\%$};
      \draw (-0.9,-5.1892) -- (-1.1,-5.1892) node[left] {$20\%$};
      \draw (-0.9,1.2162) -- (-1.1,1.2162) node[left] {$30\%$};
      \draw (-0.9,-4.7838) -- (-1.1,-4.7838) node[left] {$30\%$};
      \draw (-0.9,1.6216) -- (-1.1,1.6216) node[left] {$40\%$};
      \draw (-0.9,-4.3784) -- (-1.1,-4.3784) node[left] {$40\%$};
      \draw (-0.9,2.027) -- (-1.1,2.027) node[left] {$50\%$};
      \draw (-0.9,-3.973) -- (-1.1,-3.973) node[left] {$50\%$};
      \draw (-0.9,2.4324) -- (-1.1,2.4324) node[left] {$60\%$};
      \draw (-0.9,-3.5676) -- (-1.1,-3.5676) node[left] {$60\%$};
      \draw (-0.9,2.8378) -- (-1.1,2.8378) node[left] {$70\%$};
      \draw (-0.9,-3.1622) -- (-1.1,-3.1622) node[left] {$70\%$};
      \draw (-0.9,3.2432) -- (-1.1,3.2432) node[left] {$80\%$};
      \draw (-0.9,-2.7568) -- (-1.1,-2.7568) node[left] {$80\%$};
      \draw (-0.9,3.6486) -- (-1.1,3.6486) node[left] {$90\%$};
      \draw (-0.9,-2.3514) -- (-1.1,-2.3514) node[left] {$90\%$};
      \draw (0.19074,0.1) -- (0.19074,-0.1) node[below] {$-20$};
      \draw (0.19074,-5.9) -- (0.19074,-6.1) node[below] {$-20$};
      \draw (2.5319,0.1) -- (2.5319,-0.1) node[below] {$-15$};
      \draw (2.5319,-5.9) -- (2.5319,-6.1) node[below] {$-15$};
      \draw (4.873,0.1) -- (4.873,-0.1) node[below] {$-10$};
      \draw (4.873,-5.9) -- (4.873,-6.1) node[below] {$-10$};
      \draw (7.2142,0.1) -- (7.2142,-0.1) node[below] {$-5$};
      \draw (7.2142,-5.9) -- (7.2142,-6.1) node[below] {$-5$};
      \draw (9.5553,0.1) -- (9.5553,-0.1) node[below] {$0$};
      \draw (9.5553,-5.9) -- (9.5553,-6.1) node[below] {$0$};
      \fill[sectionblue] (0,0) -- (0.2,0) -- (0.2,0.00063126) -- (0,0.00063126) -- cycle;
      \fill[sectionblue] (0.2,0) -- (0.4,0) -- (0.4,0.00097558) -- (0.2,0.00097558) -- cycle;
      \fill[sectionblue] (0.4,0) -- (0.6,0) -- (0.6,0.0011248) -- (0.4,0.0011248) -- cycle;
      \fill[sectionblue] (0.6,0) -- (0.8,0) -- (0.8,0.0012855) -- (0.6,0.0012855) -- cycle;
      \fill[sectionblue] (0.8,0) -- (1,0) -- (1,0.0013199) -- (0.8,0.0013199) -- cycle;
      \fill[sectionblue] (0.8,-6) -- (1,-6) -- (1,-5.9865) -- (0.8,-5.9865) -- cycle;
      \fill[sectionblue] (1,0) -- (1.2,0) -- (1.2,0.0014462) -- (1,0.0014462) -- cycle;
      \fill[sectionblue] (1.2,0) -- (1.4,0) -- (1.4,0.0015839) -- (1.2,0.0015839) -- cycle;
      \fill[sectionblue] (1.4,0) -- (1.6,0) -- (1.6,0.0016987) -- (1.4,0.0016987) -- cycle;
      \fill[sectionblue] (1.6,0) -- (1.8,0) -- (1.8,0.0017446) -- (1.6,0.0017446) -- cycle;
      \fill[sectionblue] (1.8,0) -- (2,0) -- (2,0.0015035) -- (1.8,0.0015035) -- cycle;
      \fill[sectionblue] (2,0) -- (2.2,0) -- (2.2,0.0013084) -- (2,0.0013084) -- cycle;
      \fill[sectionblue] (2.2,0) -- (2.4,0) -- (2.4,0.0015954) -- (2.2,0.0015954) -- cycle;
      \fill[sectionblue] (2.4,0) -- (2.6,0) -- (2.6,0.00075751) -- (2.4,0.00075751) -- cycle;
      \fill[sectionblue] (2.6,0) -- (2.8,0) -- (2.8,0.00059683) -- (2.6,0.00059683) -- cycle;
      \fill[sectionblue] (2.8,0) -- (3,0) -- (3,0.00044762) -- (2.8,0.00044762) -- cycle;
      \fill[sectionblue] (2.8,-6) -- (3,-6) -- (3,-5.9865) -- (2.8,-5.9865) -- cycle;
      \fill[sectionblue] (3,0) -- (3.2,0) -- (3.2,0.0003558) -- (3,0.0003558) -- cycle;
      \fill[sectionblue] (3.2,0) -- (3.4,0) -- (3.4,0.00019512) -- (3.2,0.00019512) -- cycle;
      \fill[sectionblue] (3.4,0) -- (3.6,0) -- (3.6,6.8865e-05) -- (3.4,6.8865e-05) -- cycle;
      \fill[sectionblue] (5.2,0) -- (5.4,0) -- (5.4,0.00064274) -- (5.2,0.00064274) -- cycle;
      \fill[sectionblue] (5.4,0) -- (5.6,0) -- (5.6,0.0037875) -- (5.4,0.0037875) -- cycle;
      \fill[sectionblue] (5.6,0) -- (5.8,0) -- (5.8,0.0047172) -- (5.6,0.0047172) -- cycle;
      \fill[sectionblue] (5.8,0) -- (6,0) -- (6,0.0057387) -- (5.8,0.0057387) -- cycle;
      \fill[sectionblue] (6,0) -- (6.2,0) -- (6.2,0.0072193) -- (6,0.0072193) -- cycle;
      \fill[sectionblue] (6,-6) -- (6.2,-6) -- (6.2,-5.9865) -- (6,-5.9865) -- cycle;
      \fill[sectionblue] (6.2,0) -- (6.4,0) -- (6.4,0.006588) -- (6.2,0.006588) -- cycle;
      \fill[sectionblue] (6.4,0) -- (6.6,0) -- (6.6,0.0036383) -- (6.4,0.0036383) -- cycle;
      \fill[sectionblue] (6.6,0) -- (6.8,0) -- (6.8,0.0023529) -- (6.6,0.0023529) -- cycle;
      \fill[sectionblue] (6.8,0) -- (7,0) -- (7,0.0010904) -- (6.8,0.0010904) -- cycle;
      \fill[sectionblue] (7,0) -- (7.2,0) -- (7.2,0.0027431) -- (7,0.0027431) -- cycle;
      \fill[sectionblue] (7.2,0) -- (7.4,0) -- (7.4,0.0099854) -- (7.2,0.0099854) -- cycle;
      \fill[sectionblue] (7.2,-6) -- (7.4,-6) -- (7.4,-5.9865) -- (7.2,-5.9865) -- cycle;
      \fill[sectionblue] (7.4,0) -- (7.6,0) -- (7.6,0.014048) -- (7.4,0.014048) -- cycle;
      \fill[sectionblue] (7.4,-6) -- (7.6,-6) -- (7.6,-5.9865) -- (7.4,-5.9865) -- cycle;
      \fill[sectionblue] (7.6,0) -- (7.8,0) -- (7.8,0.01414) -- (7.6,0.01414) -- cycle;
      \fill[sectionblue] (7.8,0) -- (8,0) -- (8,0.0089065) -- (7.8,0.0089065) -- cycle;
      \fill[sectionblue] (8,0) -- (8.2,0) -- (8.2,0.019431) -- (8,0.019431) -- cycle;
      \fill[sectionblue] (8,-6) -- (8.2,-6) -- (8.2,-5.973) -- (8,-5.973) -- cycle;
      \fill[sectionblue] (8.2,0) -- (8.4,0) -- (8.4,0.025296) -- (8.2,0.025296) -- cycle;
      \fill[sectionblue] (8.2,-6) -- (8.4,-6) -- (8.4,-5.9595) -- (8.2,-5.9595) -- cycle;
      \fill[sectionblue] (8.4,0) -- (8.6,0) -- (8.6,0.033721) -- (8.4,0.033721) -- cycle;
      \fill[sectionblue] (8.4,-6) -- (8.6,-6) -- (8.6,-5.9324) -- (8.4,-5.9324) -- cycle;
      \fill[sectionblue] (8.6,0) -- (8.8,0) -- (8.8,0.050248) -- (8.6,0.050248) -- cycle;
      \fill[sectionblue] (8.6,-6) -- (8.8,-6) -- (8.8,-5.973) -- (8.6,-5.973) -- cycle;
      \fill[sectionblue] (8.8,0) -- (9,0) -- (9,0.082592) -- (8.8,0.082592) -- cycle;
      \fill[sectionblue] (8.8,-6) -- (9,-6) -- (9,-5.8649) -- (8.8,-5.8649) -- cycle;
      \fill[sectionblue] (9,0) -- (9.2,0) -- (9.2,0.15399) -- (9,0.15399) -- cycle;
      \fill[sectionblue] (9,-6) -- (9.2,-6) -- (9.2,-5.8919) -- (9,-5.8919) -- cycle;
      \fill[sectionblue] (9.2,0) -- (9.4,0) -- (9.4,0.35587) -- (9.2,0.35587) -- cycle;
      \fill[sectionblue] (9.2,-6) -- (9.4,-6) -- (9.4,-5.6351) -- (9.2,-5.6351) -- cycle;
      \fill[sectionblue] (9.4,0) -- (9.6,0) -- (9.6,1.0477) -- (9.4,1.0477) -- cycle;
      \fill[sectionblue] (9.4,-6) -- (9.6,-6) -- (9.6,-4.9595) -- (9.4,-4.9595) -- cycle;
      \fill[sectionblue] (9.6,0) -- (9.8,0) -- (9.8,2.181) -- (9.6,2.181) -- cycle;
      \fill[linkred] (9.6,-6) -- (9.8,-6) -- (9.8,-2) -- (9.6,-2) -- cycle;
      \fill[sectionblue] (9.6,-6) -- (9.8,-6) -- (9.8,-3.8243) -- (9.6,-3.8243) -- cycle;
      \fill[pattern=north east lines,pattern color=linkred] (9.6,-6) -- (9.8,-6) -- (9.8,-3.8243) -- (9.6,-3.8243) -- cycle;
      \fill[linkred] (9.8,-6) -- (10,-6) -- (10,-5.9459) -- (9.8,-5.9459) -- cycle;
    \end{tikzpicture}
  \end{center}
  \caption{Histogram showing the frequency of $k$ and $\kappa$ values for the
  translation-invariant CTMC (with $841$ states).
  The upper histogram shows how the values $k(r,s)$ for all state pairs $r,s$ are
  distributed. The lower histogram shows the distribution of both $k(r,s)$ and $\kappa(r,s)$, but
  only for $300$ randomly selected pairs
  ($k$ values in blue, $\kappa$ values in red).}
  \label{fig:tinvctmc_k_histogram}
\end{figure}

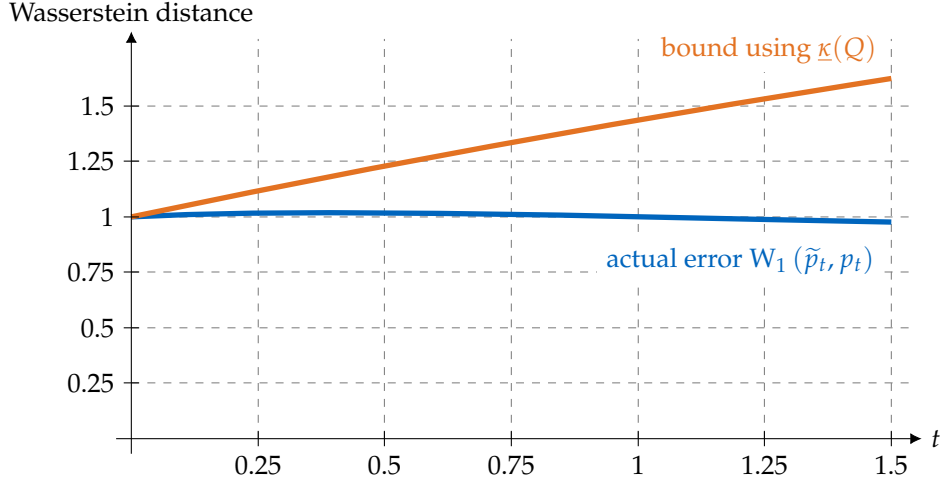
\begin{figure}[htb]
  \begin{center}
    \begin{tikzpicture}[>={Latex[length=1.5mm,width=1.5mm]}]
      \draw[gray,dashed] (1.67, 0.00) -- (1.67, 5.30);
      \draw[gray,dashed] (3.33, 0.00) -- (3.33, 5.30);
      \draw[gray,dashed] (5.00, 0.00) -- (5.00, 5.30);
      \draw[gray,dashed] (6.67, 0.00) -- (6.67, 5.30);
      \draw[gray,dashed] (8.33, 0.00) -- (8.33, 5.30);
      \draw[gray,dashed] (10.00, 0.00) -- (10.00, 5.30);
      \draw[gray,dashed] (0.00, 0.74) -- (10.30, 0.74);
      \draw[gray,dashed] (0.00, 1.47) -- (10.30, 1.47);
      \draw[gray,dashed] (0.00, 2.21) -- (10.30, 2.21);
      \draw[gray,dashed] (0.00, 2.94) -- (10.30, 2.94);
      \draw[gray,dashed] (0.00, 3.68) -- (10.30, 3.68);
      \draw[gray,dashed] (0.00, 4.41) -- (10.30, 4.41);
      \begin{scope}
        \clip (-0.04, -0.04) rectangle (10.04, 5.04);
        \draw[sectionblue,line width=2pt] (0, 2.941) -- (0.725, 2.971) -- (1.594, 2.99) -- (2.609, 2.996) -- (4.058, 2.987) -- (5.652, 2.964) -- (8.406, 2.905) -- (10, 2.872);
        \draw[tumOrange,line width=2pt] (0, 2.941) -- (1.579, 3.268) -- (3.158, 3.581) -- (4.737, 3.878) -- (6.316, 4.163) -- (7.895, 4.434) -- (9.474, 4.692) -- (10, 4.776);
      \end{scope}
      \draw[->] (-0.20, 0.00) -- (10.40, 0.00) node[right] {$t$};
      \draw[->] (-0.00, -0.20) -- (-0.00, 5.40) node[above] {Wasserstein distance};
      \draw (1.67, 0.10) -- (1.67, -0.10) node[below] {$0.25$};
      \draw (3.33, 0.10) -- (3.33, -0.10) node[below] {$0.5$};
      \draw (5.00, 0.10) -- (5.00, -0.10) node[below] {$0.75$};
      \draw (6.67, 0.10) -- (6.67, -0.10) node[below] {$1$};
      \draw (8.33, 0.10) -- (8.33, -0.10) node[below] {$1.25$};
      \draw (10.00, 0.10) -- (10.00, -0.10) node[below] {$1.5$};
      \draw (0.10, 0.74) -- (-0.10, 0.74) node[left] {$0.25$};
      \draw (0.10, 1.47) -- (-0.10, 1.47) node[left] {$0.5$};
      \draw (0.10, 2.21) -- (-0.10, 2.21) node[left] {$0.75$};
      \draw (0.10, 2.94) -- (-0.10, 2.94) node[left] {$1$};
      \draw (0.10, 3.68) -- (-0.10, 3.68) node[left] {$1.25$};
      \draw (0.10, 4.41) -- (-0.10, 4.41) node[left] {$1.5$};
      \node[sectionblue,below left=2mm and 1mm,fill=white] at (10.00, 2.87) {actual error $\WD{\widetilde{p}_t}{p_t}$};
      \node[tumOrange,above left=1mm,fill=white] at (10.00, 4.78) {bound using $\underline{\kappa}(Q)$};
    \end{tikzpicture}
  \end{center}
  \caption{Evolution of the actual error $\WD{\widetilde{p}_t}{p_t}$ and the bound for the translation-invariant CTMC
  (with $841$ states),
  aggregated using a simple coarse gridding approach
  (resulting in $225$ aggregates). The initial distribution $p_0$ was chosen to be a 
  Dirac measure on the state closest to the center of the original state grid.}
  \label{fig:tinvctmc_aggregation_error}
\end{figure}

We also created an example CTMC on a finite subset of $\bbZ^2$ with properties very similar to the ones
used in \autoref{prop:transl_inv_ctmc_nonneg_ricci}, that is, (almost) a translation-invariant CTMC.
However, from every state, we added a transition to one single root state. This
breaks translation-invariance but actually results in a higher Ricci curvature
(and we also wanted to include an example with a high curvature).
As a metric, we used the Manhattan metric or $\norm{\cdot}_1$-norm.
\autoref{fig:tinvctmc_k_histogram} shows that the values of $k(r, s)$ can still be
negative, even though \autoref{prop:transl_inv_ctmc_nonneg_ricci} guarantees
non-negative curvature (we would need to extend the proof slightly to cover
the transitions to the root state, but the conclusion of \autoref{prop:transl_inv_ctmc_nonneg_ricci}
does indeed apply to our example). While $\underline{k}(Q) \approx -20.4$ is higher
than in the previous examples, it is still too low for useful error bounds, even
though we have $\underline{\kappa}(Q) \approx 0.2 > 0$ in this example. Here, we simply
aggregated four neighboring grid points, resulting in $225$ aggregates
and $\norm{\Theta A - A Q}_{\mathrm{W}} \approx 0.68$.
The bound using $\underline{\kappa}(Q)$ is initially of a magnitude comparable to
the actual error, see \autoref{fig:tinvctmc_aggregation_error}.

\begin{figure}[htb]
  \begin{center}
    \begin{tikzpicture}[>={Latex[length=1.5mm,width=1.5mm]}]
      \draw[gray,dashed] (2.50, 0) -- (2.50, 5.30);
      \draw[gray,dashed] (5.00, 0) -- (5.00, 5.30);
      \draw[gray,dashed] (7.50, 0) -- (7.50, 5.30);
      \draw[gray,dashed] (10.00, 0) -- (10.00, 5.30);
      \draw[gray,dashed] (0, 1.41) -- (10.20, 1.41);
      \draw[gray,dashed] (0, 2.81) -- (10.20, 2.81);
      \draw[gray,dashed] (0, 4.22) -- (10.20, 4.22);
      \draw[->] (-0.20, 0) -- (10.30, 0) node[right] {$t$};
      \draw[->] (0, -0.20) -- (0, 5.40) node[above] {Wasserstein distance};
      \draw (2.50, 0.1) -- (2.50, -0.1) node[below] {0.25};
      \draw (5.00, 0.1) -- (5.00, -0.1) node[below] {0.5};
      \draw (7.50, 0.1) -- (7.50, -0.1) node[below] {0.75};
      \draw (10.00, 0.1) -- (10.00, -0.1) node[below] {1};
      \draw (0.1, 1.41) -- (-0.1, 1.41) node[left] {0.05};
      \draw (0.1, 2.81) -- (-0.1, 2.81) node[left] {0.1};
      \draw (0.1, 4.22) -- (-0.1, 4.22) node[left] {0.15};
      \draw[sectionblue,line width=2pt] (0.00, 0) -- (0.14, 0.137) -- (0.29, 0.269) -- (0.43, 0.396) -- (0.58, 0.519) -- (0.72, 0.636) -- (0.87, 0.750) -- (1.01, 0.860) -- (1.16, 0.966) -- (1.30, 1.069) -- (1.45, 1.168) -- (1.59, 1.265) -- (1.74, 1.359) -- (1.88, 1.450) -- (2.03, 1.539) -- (2.17, 1.625) -- (2.32, 1.709) -- (2.46, 1.791) -- (2.61, 1.871) -- (2.75, 1.949) -- (2.90, 2.025) -- (3.04, 2.100) -- (3.19, 2.173) -- (3.33, 2.244) -- (3.48, 2.314) -- (3.62, 2.382) -- (3.77, 2.449) -- (3.91, 2.515) -- (4.06, 2.579) -- (4.20, 2.643) -- (4.35, 2.705) -- (4.49, 2.766) -- (4.64, 2.825) -- (4.78, 2.884) -- (4.93, 2.942) -- (5.07, 2.999) -- (5.22, 3.055) -- (5.36, 3.110) -- (5.51, 3.164) -- (5.65, 3.217) -- (5.80, 3.270) -- (5.94, 3.322) -- (6.09, 3.373) -- (6.23, 3.423) -- (6.38, 3.473) -- (6.52, 3.522) -- (6.67, 3.570) -- (6.81, 3.618) -- (6.96, 3.665) -- (7.10, 3.712) -- (7.25, 3.758) -- (7.39, 3.803) -- (7.54, 3.848) -- (7.68, 3.893) -- (7.83, 3.937) -- (7.97, 3.980) -- (8.12, 4.023) -- (8.26, 4.066) -- (8.41, 4.108) -- (8.55, 4.150) -- (8.70, 4.191) -- (8.84, 4.232) -- (8.99, 4.272) -- (9.13, 4.312) -- (9.28, 4.352) -- (9.42, 4.391) -- (9.57, 4.430) -- (9.71, 4.469) -- (9.86, 4.507) -- (10.00, 4.545) node[above left=0mm and -3mm,fill=white] {actual error $\WD{\widetilde{p}_t}{p_t}$};
      \draw[tumOrange,line width=2pt] (0.00, 0) -- (0.49, 0.475) -- (0.98, 0.950) -- (1.47, 1.425) -- (1.96, 1.901) -- (2.45, 2.377) -- (2.94, 2.854) -- (3.43, 3.331) -- (3.92, 3.809) -- (4.41, 4.287) -- (4.90, 4.765) node[above,fill=white] {bound using $\underline{k}(Q)$};
    \end{tikzpicture}
  \end{center}
  \caption{Evolution of the actual error $\WD{\widetilde{p}_t}{p_t}$ and the bound for the CTMC
  arising from a (state) discretization of the CGMY process (with $800$ states),
  aggregated using a simple coarse gridding approach
  (resulting in $160$ aggregates). The initial distribution $p_0$ was chosen to be a 
  uniform distribution over the aggregate containing the state $0$ (or rather, the
  discretized state which represents the original state $0$ of the CGMY process).}
  \label{fig:cgmy_aggregation_error}
\end{figure}
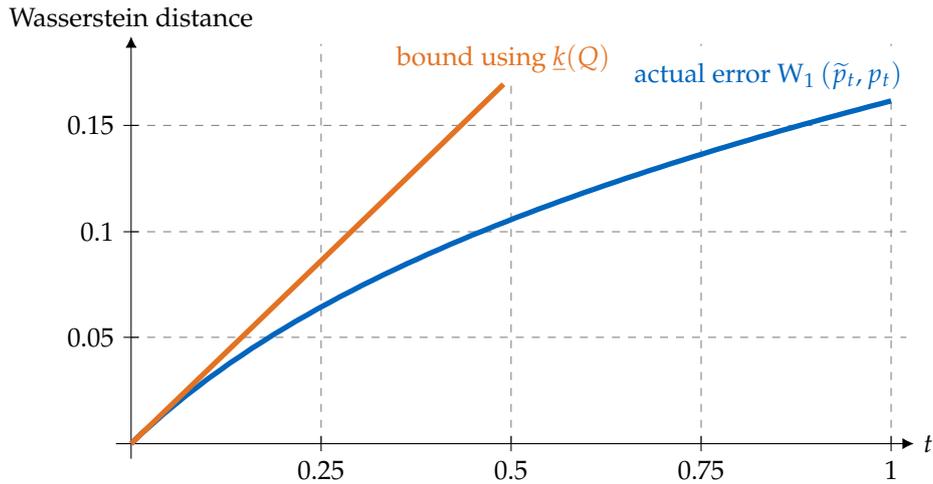

Another example where the error bounds do not explode are discretizations
of Lévy processes or Lévy-driven queues because such discretizations (usually) also
satisfy the assumptions of \autoref{prop:transl_inv_ctmc_nonneg_ricci}.
One area in which Lévy processes are often used as models is finance,
and one particular process used for modelling asset returns is the
CGMY process \cite{cgmy}. We discretized and truncated the original state
space $\bbR$ to obtain a CTMC with $800$ states, resulting in $\underline{\kappa}(Q) \approx 0$,
while $\underline{k}(Q) \approx -0.018$. \autoref{fig:cgmy_aggregation_error}
shows that the error bound using $\underline{k}(Q)$ matches the actual error
almost exactly near $t = 0$ before the distance between the two grows.
In this case, we aggregated five neighboring states on the line, resulting
in $160$ aggregates and $\norm{\Theta A - A Q}_{\mathrm{W}} \approx 0.34$.
In the future, we would like to adapt the error bounds to the Markov process setting
such that the bounds can be used directly for the distance between the original,
continuous transient distribution and the approximation.

\subsection{The DTMC case}
\label{ssec:w_error_bounds_dtmc}

For completeness, we provide a very short overview of how results analogous
to those from the previous section look for DTMCs.
In discrete time, the calculations are simpler, which is why this paper was
focused on CTMCs~-- the more complicated case. For DTMCs, the final inequality
from the proof of \autoref{thm:wasserstein_error_growth_bound} becomes
\begin{align}
  \begin{split}
    \WD{\widetilde{p}_{k + 1}}{p_{k + 1}}
    &\overset{\triangle\textrm{-inequ.}}{\leq} \WD{\widetilde{p}_{k + 1}^{\transp}}{\widetilde{p}_k^{\transp} P} + \WD{\widetilde{p}_k^{\transp} P}{p_{k + 1}^{\transp}} \\
    &\;\;\;\;=\;\;\;\; \WD{\pi_k^{\transp} \Pi A}{\pi_k^{\transp} A P} + \WD{\widetilde{p}_k^{\transp} P}{p_k^{\transp} P}
  \end{split}
  \label{eq:dtmc_step_bound}
\end{align}
Now, on the one hand,
\begin{align*}
  \WD{\pi_k^{\transp} \Pi A}{\pi_k^{\transp} A P}
  \;\overset{\textrm{\eqref{eq:wasserstein_dual_finite}}}{=}\;
  \max_{\substack{f \in \bbR^n \textrm{ is }1\textrm{-Lipschitz w.r.t.\ }\DIST\\\forall s \in S: 0 \leq f(s) \leq d_{\max}}}
  \pi_k^{\transp}(\Pi A - A P) f \; \leq \; \pi_k^{\transp} \dabs{\Pi A - A P}_{\mathrm{W}}
\end{align*}
On the other hand, by \cite[Proposition 20]{riccimarkovmetricspaces}
(we only use one direction of the proposition),
\begin{align*}
  \WD{\widetilde{p}_k^{\transp} P}{p_k^{\transp} P}
  &\leq \big(1 - \underline{\kappa}(P)\big) \cdot \WD{\widetilde{p}_k}{p_k}
\end{align*}
where $\underline{\kappa}(P)$ was defined in \autoref{def:ricci_curvature_dtmc}.
Plugging these two bounds into \eqref{eq:dtmc_step_bound} yields
\begin{align*}
  \WD{\widetilde{p}_{k + 1}}{p_{k + 1}} &\leq \pi_k^{\transp} \dabs{\Pi A - A P}_{\mathrm{W}} + \big(1 - \underline{\kappa}(P)\big) \cdot \WD{\widetilde{p}_k}{p_k}
\end{align*}
which proves a statement analogous to \autoref{thm:wasserstein_error_growth_bound}
for the discrete-time case. We leave it to the reader to derive bounds for $\underline{\kappa}(P)$
similar to $\underline{k}(Q)$ and $K(Q)$ which are easier to calculate than $\underline{\kappa}(P)$ itself.

\section{Conclusion}
\label{sec:conclusion}

We have seen how the bounds presented in \cite{formalbndsstatespaceredmc} can be
extended from measuring the aggregation error in total variation to measuring the error
in the Wasserstein distance w.r.t.\ an arbitrary metric on the finite state
space of a CTMC (or DTMC). The error caused by approximating the model dynamics
with a Markov chain on a lower dimensional state space can be bounded by
a Wasserstein matrix norm on $\Theta A - A Q$ (where $\Theta$ is the aggregated generator,
$A$ the disaggregation matrix and $Q$ the original generator), which is a very similar result
to the one from \cite{formalbndsstatespaceredmc}. The propagation of the accumulated error can be
controlled using the coarse Ricci curvature of the Markov chain.
When the curvature is positive, the bound on the accumulated error will decrease over time;
if it is negative, the bound will grow exponentially. The discrete metric
ensures non-negative curvature and thus a non-increasing accumulated error,
which explains the absence of an additional error term in \cite{formalbndsstatespaceredmc}.
In fact, the curvature can be strictly positive, improving the bounds from
\cite{formalbndsstatespaceredmc} in such settings.

Next to the discrete metric, we also saw that translation-invariant CTMCs
result in a non-negative Ricci curvature, which is desirable to obtain practically
useful error bounds which do not blow up exponentially. However, when applying
the bounds to the examples from \cite{formalbndsstatespaceredmc} equipped
with (more or less) natural metrics, we also saw that negative Ricci curvature
(in different orders of magnitude) can easily render the error bounds useless.
This effect is further aggravated when using the easier-to-calculate $\underline{k}$
instead of the Ricci curvature. The examples demonstrated that only a small
portion of the state pairs can cause the negative curvature while the major part
of pairs is better behaved.

\subsection{Future work}

The Wasserstein error bounds presented in this paper are a first step towards
extending the error bounds for aggregation to general continuous-time Markov
processes with continuous state spaces. Often, the only way to calculate transient
distributions for these is by discretization (which can be seen as aggregation),
and formal bounds for the introduced error in the approximation of the
transient distributions are missing.
The total variation distance is usually not a good choice to measure the
error between an approximated and actual probability distribution on a continuous
state space, but the Wasserstein distance (with the right metric) is more appropriate.
While there is more work to do, many of the results of this paper can
probably be extended to general Markov processes. This is the subject of
ongoing research.

Next to the extension to general Markov processes, a crucial research
question is the practical applicability of the presented error bounds.
The main issue seems to be processes with negative Ricci curvature.
As a negative curvature quickly leads to deteriorating bounds,
it should be investigated whether better bounds can be derived if only
a small part of the state pairs has a negative curvature while the rest
has positive or close-to-zero curvature. Using ``coarse Ricci curvature up to $\delta$''
as defined in \cite[Definition 57]{riccimarkovmetricspaces} instead of the
coarse Ricci curvature could be a way to tackle that problem.
Another research subject would
be to identify processes with non-negative curvature more broadly than
done in this paper, as the error bounds would work well with those.
Are there important examples beyond the discrete metric and translation-invariant
CTMCs (respectively Lévy processes in the more general Markov process setting)?

\section{Preview: error bounds for Markov processes}

Here, we give an overview of how difficult or straightforward an extension
of the theory to the Markov process setting seems to be.

\subsection{Preliminaries}

We consider a Markov process $X_t$ with a continuous state space $S$
and in continuous time. We assume (at least) that $S$ is Polish and
that it is equipped with some lower-semicontinuous metric $\DIST$
(which need not be a metric giving rise to the underlying topology of
$S$). In probability theory, such a Markov process
is also described by a semigroup $P_t$ and by a generator $\calL$. Both
are linear operators on the space of functions from $S$ to $\bbR$.
\begin{itemize}
  \item We have $P_t f(x) = \Exc{x}{f(X_t)}$. If we consider a CTMC $Y_t$ with generator $Q$ on a finite state space, then
  the linear operator $P_t$ is represented by the matrix $e^{tQ}$. Indeed,
  we can represent a function $f$ from the finite state space $\{1,\ldots,n\}$ of $Y_t$ to $\bbR$
  as a vector in $\bbR^n$: $\vec{f} := (f(1), \ldots, f(n))^{\transp}$. We then have:
  \begin{align*}
    \Exc{x}{f(Y_t)} = \1_x^{\transp} \cdot e^{tQ} \cdot \vec{f}
    \qquad\textrm{ where } \1_x \in \bbR^n, \1_x(y) = \begin{cases}
      1 & \textrm{ if } x = y \\
      0 & \textrm{ otherwise}
    \end{cases}
  \end{align*}
  The left multiplication with $\1_x^{\transp}$ amounts to evaluating
  the function $e^{tQ} \cdot \vec{f}$ (which is interpreted as a vector)
  at point $x$.

  It holds that $P_0 = I$ (the identity), and $P_s \circ P_t = P_t \circ P_s = P_{s+t}$.
  $P_t$ is basically a stochastic matrix, but on an infinite-dimensional state space.
  \item We have $\calL f(x) = \left.\frac{\textrm{d}}{\textrm{d}t}\right|_{t=0} \, \Exc{x}{f(X_t)}$.
  If we consider again the CTMC $Y_t$, then $\calL$ corresponds to the generator matrix $Q$. Indeed,
  \begin{align*}
    \left.\frac{\textrm{d}}{\textrm{d}t}\right|_{t=0} \, \Exc{x}{f(Y_t)}
    &= \left.\frac{\textrm{d}}{\textrm{d}t}\right|_{t=0} \, \1_x^{\transp} \cdot e^{tQ} \cdot \vec{f}
    = \1_x^{\transp} \cdot Q \cdot \vec{f}
  \end{align*}
\end{itemize}

In analogy to the CTMC case, the transient distribution $p_t$ of the Markov process,
i.e., the law of $X_t$, is defined as
\begin{align*}
  p_t = p_0 P_t \quad \textrm{ where } \quad
  p_t(A) = \int_S P_t(x, A) \dmux{p_0}{x} \quad \textrm{ for } A \textrm{ measurable}
\end{align*}
where $p_0$ is the initial distribution and $P_t(x, A) = \Prbc{x}{X_t \in A}$.

As detailed in \cite{levymatters} (see, e.g., the preface), basically all Markov
processes which are interesting for applications can be understood as a family of
Lévy processes: these processes are characterized by a state-dependent drift, diffusion coefficient
and jump measure, in contrast to Lévy processes where all three components are
independent of the current state of the process. With the application of
numerically approximating transient laws in mind, a process description via
this so-called state-dependent Lévy triplet is more tractable than the abstract
generator $\calL$ when trying to derive error bounds. Hence, bounds for
continuous-time and -state Markov processes should probably be derived for
such a process description.

In order to be actually able to compute a transient distribution of
a Markov process (which is a non-trivial problem),
two possible forms of discretizations immediately come to mind: discretizing only the
state space to obtain a CTMC, and discretizing both states and time to obtain
a DTMC. We will sketch these two approaches, but it might be even better
to combine them somehow or slightly alter some of the given details.

\subsubsection{CTMC approximation}

In this setting, we approximate $p_t$ by $\widetilde{p}_t$, defined as
\begin{align*}
  \widetilde{p}_t = \sum_{i=1}^n \pi_t(i) \cdot a_i \quad \textrm{ where } \quad
  \pi_t^{\transp} = \pi_0^{\transp} e^{t \Theta}
\end{align*}
with $\Theta \in \bbR^{n \times n}$ the generator matrix of the aggregated
CTMC model, $\pi_0 \in \bbR^n$ the aggregated initial distribution, and
$a_1, \ldots, a_n$ probability measures on the original state space $S$
(``disaggregation measures''). $a_i$ describes how the probability mass
in aggregate $i$, that is $\pi_t(i)$, should be distributed among the original
states in the disaggregation phase. For example, if $S = \bbR$, then
$a_i$ could be a uniform distribution over some interval, which would
imply that the states in that interval are represented by aggregate $i$
in the aggregated model.

\subsubsection{DTMC approximation}

In this setting, we approximate $p_{k\Delta}$ by $\widetilde{p}_k$, defined as
\begin{align*}
  \widetilde{p}_k = \sum_{i=1}^n \pi_k(i) \cdot a_i \quad \textrm{ where } \quad
  \pi_k^{\transp} = \pi_0^{\transp} \Pi^k
\end{align*}
with $\Pi \in \bbR^{n \times n}$ the stochastic transition matrix of the aggregated
DTMC model, $\pi_0 \in \bbR^n$ the aggregated initial distribution,
$a_1, \ldots, a_n$ probability measures on the original state space $S$
(``disaggregation measures''), and $\Delta$ the time discretization parameter / step size.

\subsection{Wasserstein error bounds}

Again, we would like to provide formal error bounds on the distance between
the actual and approximated transient distributions, the latter obtained via
the discretization procedure sketched above. The following two sections give
an overview of the necessary steps to prove bounds similar to the ones
shown in \autoref{thm:wasserstein_error_growth_bound}.

\subsubsection{CTMC approximation}

The basic goal would be to bound $\dufrac{t^+} \WD{\widetilde{p}_t}{p_t}$.
Required steps:
\begin{enumerate}[(1)]
  \item Is $\WD{\widetilde{p}_t}{p_t}$ continuous? Otherwise a derivative bound
  is not useful. $\WD{\widetilde{p}_t}{p_t}$ can be discontinuous, e.g.\ if the
  discrete metric on $S$ is used as $\DIST$.

  It is enough to show
  \begin{align*}
    \WD{\widetilde{p}_t}{\widetilde{p}_{t+u}} \overset{u \to 0}{\longrightarrow} 0
    \quad \textrm{ and } \quad
    \WD{p_t}{p_{t+u}} \overset{u \to 0}{\longrightarrow} 0
  \end{align*}
  The first condition should be true as $\pi_t$ is continuous and even
  differentiable. The second condition can easily be violated if the
  discrete metric is used and $p_t$ is e.g.\ a Dirac measure on $t \in \bbR$.

  We would need to find conditions under which $\WD{\widetilde{p}_t}{p_t}$ is continuous.
  Are error bounds possible if $\WD{\widetilde{p}_t}{p_t}$ is not continuous?
  \item Does $\dufrac{t^+} \WD{\widetilde{p}_t}{p_t}$ exist? Or should we consider
  some $\limsup$ instead?\
  \item Find a bound on
  $\displaystyle \WD{\widetilde{p}_{t + u}}{\widetilde{p}_{t} P_u}$ or directly bound
  $\displaystyle \dufrac{u^+} \WD{\widetilde{p}_{t + u}}{\widetilde{p}_{t} P_u}$
  (i.e., try to find an equivalent of \autoref{cor:wdderiv_approxdyn} / \autoref{cor:wdderiv_approxdyn_bound}).
  
  If we want to bound the derivative directly, we could try to use
  Danskin's Theorem, applied to
  \begin{align*}
    \WD{\widetilde{p}_{t + u}}{\widetilde{p}_{t} P_u}
    = \sup_{f \textrm{ bounded and Lipschitz}} \left(\int_S f \dx{\widetilde{p}_{t + u}} - \int_S f \dx{\widetilde{p}_{t} P_u}\right)
  \end{align*}
  Problem 1: the supremum need not be a maximum. Use the coupling definition
  instead (where the minimum is achieved)?

  Problem 2: the generator $\calL$ cannot be applied to all bounded and
  Lipschitz $f$ to calculate the derivative.
  \item Find a bound on
  $\displaystyle \WD{\widetilde{p}_{t} P_u}{p_{t} P_u} = \WD{\widetilde{p}_{t} P_u}{p_{t + u}}$
  depending on
  $\displaystyle \WD{\widetilde{p}_{t}}{p_{t}}$
  (i.e., find an equivalent of \autoref{lem:wdderiv_error_ricciinfbound}).

  Here, we should be able to use \cite[Theorem 1.9]{markovprocricci}.
  Note that the theorem only applies to processes admitting a left-continuous modification,
  and that we would need to slightly extend the original statement which only
  applies to Dirac initial measures.
  \item Conclude
  \begin{align*}
    \WD{\widetilde{p}_{t + u}}{p_{t + u}}
    \leq \overbrace{\WD{\widetilde{p}_{t + u}}{\widetilde{p}_{t} P_u}}^{=0 \textrm{ for }u = 0}
    + \overbrace{\WD{\widetilde{p}_{t} P_u}{p_{t + u}}}^{= \textrm{l.h.s\ for }u = 0}
  \end{align*}
\end{enumerate}

\subsubsection{DTMC approximation}

The goal here would be to bound $\WD{\widetilde{p}_{k+1}}{p_{(k+1)\Delta}}$, given
$\WD{\widetilde{p}_{k}}{p_{k\Delta}}$. We have
\begin{align*}
  \WD{\widetilde{p}_{k+1}}{p_{(k+1)\Delta}}
  \leq \WD{\widetilde{p}_{k+1}}{\widetilde{p}_{k}P_{\Delta}}
  + \WD{\widetilde{p}_{k}P_{\Delta}}{p_{k\Delta} P_{\Delta}}
\end{align*}
\begin{enumerate}[(1)]
  \item To bound $\WD{\widetilde{p}_{k+1}}{\widetilde{p}_{k}P_{\Delta}}$, it should be enough
  to consider, for all $j \in \{1, \ldots, n\}$,
  \begin{align*}
    b_j := \WD{\sum_{i=1}^n \Pi(j, i) \cdot a_i \,}{\;\; a_j P_{\Delta}}
  \end{align*}
  We should then be able to derive
  \begin{align*}
    \WD{\widetilde{p}_{k+1}}{\widetilde{p}_{k}P_{\Delta}} \leq \sum_{i=1}^n b_i \cdot \pi_k(i)
  \end{align*}
  The main problem here is how to approximate $a_j P_{\Delta}$ in practice, because
  an explicit calculation is typically impossible (if the explicit calculation
  was possible, the discretization procedure would be unnecessary).
  \item To bound $\WD{\widetilde{p}_{k}P_{\Delta}}{p_{k\Delta} P_{\Delta}}$,
  we should be able to use \cite[Theorem 1.9]{markovprocricci}
  (again only if the Markov process admits a left-continuous modification,
  and with a slight extension of the statement of the original theorem).
  A bound along the line
  \begin{align*}
    \WD{\widetilde{p}_{k}P_{\Delta}}{p_{k\Delta} P_{\Delta}} \leq \WD{\widetilde{p}_{k}}{p_{k\Delta}} e^{-\underline{\kappa} \cdot \Delta}
  \end{align*}
  should hold.
\end{enumerate}

\begin{remark}
  In the last two sections, we simply tried to sketch how to transfer the
  Wasserstein error bounds derived in this paper onto discretizations of
  Markov processes. However, one should also consider whether alternative
  discretization procedures which do not fit exactly into the CTMC or DTMC
  approximation framework above might be better suited for numerical approximation
  of transient laws (cf.\ the ideas in \cite{gendifflimmc}, for example).
\end{remark}

\end{document}